\documentclass[fleqn,leqno]{article}
\usepackage{amssymb,latexsym}
\usepackage{amsfonts}
\usepackage[noload]{qtree}
 \usepackage[dvips]{color}

\newcommand{\bi}{\bigskip}
\newcommand{\sm}{\smallskip}

\newcommand{\wh}{\widehat}
\newcommand{\ee}{\end{equation}}
\newcommand{\eea}{\end{eqnarray}}
\newcommand{\bean}{\begin{eqnarray*}}
\newcommand{\eean}{\end{eqnarray*}}
\newif\ifpctex
\newcommand{\noi}{\noindent}

\newcommand{\ve}{\varepsilon}

\newcommand{\wt}{\widetilde}

\newtheorem{theorem}{Theorem}
\newtheorem{remark}{Remark}
\newtheorem{example}{Example}
\newcommand{\D}{\displaystyle}

\newcommand{\qad}{{}\hfill \square}

\newtheorem{proposition}{Proposition}[section]
\newtheorem{definition}[proposition]{Definition}
\newtheorem{lemma}[proposition]{Lemma}
\newtheorem{corollary}[proposition]{Corollary}

\newtheorem{conjecture}[proposition]{Conjecture}
\newenvironment{proof}{\par\noindent{\bf Proof\ }}

 \renewcommand{\theequation}{\thesection.\arabic{equation}}

 \newcommand{\rand}[1]{}

\newcommand{\Rand}[1]{\marginpar{#1}}
     \renewcommand{\Rand}[1]{}

\newcommand{\be}[1]{\Rand{\vspace{0,6cm}\tt #1}\begin{equation}\label{#1}}
\newcommand{\bea}[1]{\Rand{\vspace{0,7cm}\tt
#1\vspace{-0,7cm}}\begin{eqnarray}\label{#1}} \marginparwidth2.5cm

\newcommand{\beL}[1]{\Rand{\vspace{0,6cm}\tt #1}\begin{lemma}\label{#1}}
\newcommand{\beD}[1]{\Rand{\vspace{0,6cm}\tt #1}\begin{definition}\label{#1}}
\newcommand{\beT}[1]{\Rand{\vspace{0,6cm}\tt #1}\begin{theorem}\label{#1}}
\newcommand{\beP}[1]{\Rand{\vspace{0,6cm}\tt #1}\begin{proposition}\label{#1}}
\newcommand{\beC}[1]{\Rand{\vspace{0,6cm}\tt #1}\begin{corollary}\label{#1}}
\newcommand{\beCj}[1]{\Rand{\vspace{0,6cm}\tt #1}\begin{conjecture}\label{#1}}

\newcommand{\ep}{\end{proposition}}

\newcommand{\suml}{\sum\limits}
\newcommand{\intl}{\int\limits}

\newcommand{\prodl}{\prod\limits}
\newcommand{\bigoplusl}{\bigoplus\limits}

\newcommand{\F}{\mathbb{F}}

\newcommand{\R}{\mathbb{R}}
\newcommand{\N}{\mathbb{N}}

\newcommand{\I}{\mathbb{I}}
\newcommand{\Z}{\mathbb{Z}}

\newcommand{\text}{\mbox}

\newcommand{\CA}{{\mathcal A}}
\newcommand{\CB}{{\mathcal B}}

\newcommand{\CF}{{\mathcal F}}

\newcommand{\CG}{{\mathcal G}}
\newcommand{\CI}{{\mathcal I}}
\newcommand{\CL}{{\mathcal L}}
\newcommand{\CT}{{\mathcal T}}

\newcommand{\CM}{{\mathcal M}}

\newcommand{\CP}{{\mathcal P}}

\newcommand{\CS}{{\mathcal S}}

\newcommand{\la}{\longrightarrow}

\begin{document}


\title{Duality for spatially interacting Fleming-Viot processes with mutation and
selection }
\author{{\bf Donald A. Dawson$^{1,2}$ \mbox{ } Andreas Greven$^{2,3}$}}

\date{\small  \today}

\maketitle

\begin{abstract}
Consider a system $X = ((x_\xi(t)), \xi \in \Omega_N)_{t \geq 0}$
of interacting Fleming-Viot diffusions with  mutation and selection
which is a strong Markov process with continuous paths and state space
$(\CP(\I))^{\Omega_N}$, where $\I$ is the type space,  ${\Omega_N}$
the geographic space   is assumed to be a countable group
and $\CP$ denotes the probability measures.

We establish various duality relations for this process. These dualities
are function-valued processes which are driven by a coalescing-branching
random walk, that is, an evolving particle system which in addition
exhibits certain changes in the function-valued part at jump times driven 
by mutation.

In the case of a finite type space $\I$ we construct a set-valued
dual process, which is a Markov jump process, which is very suitable to prove
ergodic theorems which we do here. The set-valued duality contains as special
case a duality relation for any finite state Markov chain.

In the finitely many types case there is also a  further tableau-valued dual
which can be used to study the invasion of fitter types after rare mutation.
This is carried out in \cite{DGsel} and \cite{DGInvasion}.

\end{abstract}

\noi {\bf Keywords:} Interacting Fleming-Viot processes, duality, mutation,
selection,  ergodic theorems.
\vspace*{2cm}

\footnoterule
\noi
\hspace*{0.3cm}{\footnotesize $^1$ {School of Mathematics and Statistics,
Carleton University, Ottawa K1S 5B6, Canada, e-mail:
ddawson@math.carleton.ca}}\\
\hspace*{0.3cm}{\footnotesize $^2$
{Research supported by NSERC, DFG-Schwerpunkt 1033 and later DFG-NWO
Forschergruppe 498, and D. Dawson's Max
Planck Award for International
Cooperation.}}\\
\hspace*{0.3cm}{\footnotesize $^3$
{Department Mathematik, Universit\"at Erlangen-N\"urnberg,
Bismarckstra{\ss}e 1 1/2,
D-91054 Erlangen, Germany, e-mail: greven@mi.uni-erlangen.de}}
\newpage

\tableofcontents

\newpage

\newcounter{secnum}
\setcounter{secnum}{\value{section}}
\setcounter{section}{-1}
\setcounter{equation}{0}
\renewcommand{\theequation}{\mbox{\arabic{secnum}.\arabic{equation}}}

\section{Background and motivation}
\label{s.backmot}

We develop here an extensive duality theory for interacting Fleming-Viot processes
with selection and mutation, which we shall exploit further in
\cite{DGsel} and \cite{DGInvasion}.

Consider the following  spatial multitype population model. The
state of a single colony is described by a probability measure on
some countable type space, the geographic space is modelled by a
set of colonies and the colonies (or demes) are labelled
with the countable hierarchical group $\Omega_N$ or some other countable
Abelian group. The stochastic
dynamics is given by a system of interacting measure-valued
diffusions and the driving mechanisms include resampling (pure
genetic drift) as diffusion term, migration, selection and mutation
as driftterms. Resampling is
modelled in each colony by the usual Fleming-Viot diffusion. The
haploid selection is based on a fitness function on the space of
types and the mutation is type-dependent.
The model belongs to a class of
processes which have been constructed via a well-posed martingale
problem by Dawson and Greven \cite{DG99}.

 More precisely the model
we use arises from the particle model driven by migration of particles
between colonies, and in each site by resampling of types,
mutation and selection. Increase the number of particles as
$\ve^{-1}$ and give them mass $\ve$. Then as $\ve\to 0$, a diffusion
limit of interacting multitype diffusions results if the
resampling rate is proportional to the number of pairs and there
is weak selection, i.e. selection occurs at a rate decreasing with
the inverse of the number of particles per site. Otherwise with
strong selection, i.e. selection at a fixed rate, we get the
deterministic limit, often referred to as the infinite population
limit.

An important technical tool in the analysis of Fleming-Viot processes  are representations
of the marginal distributions of the basic processes in terms of
expectations under the appropriate {\em function-valued}, respectively {\em set-valued}, dual processes.

Recall that two Markov processes $(Z(t))_{t \geq 0}, (Z^\prime(t))_{t \geq 0}$
on Polish spaces $E$, respectively $E^\prime$ are called {\em dual} w.r.t. to
the duality function
$H:E \times E^\prime \to \R$, if for (deterministic) initial states $Z_0$ respectively $Z^\prime_0$
for $Z$ respectively $Z^\prime$, the following identity holds:
\be{add100}
E[H(Z(t), Z^\prime_0)] = E[H(Z_0, Z^\prime(t))] \quad , \quad
\forall \; Z_0 \in E, Z^\prime_0 \in E^\prime.
\ee
This often allows us to draw conclusions on the process $Z$ from a (simpler) process $Z^\prime$ on $E^\prime$.

It has been realized for a long time that
duals play an important role in the analysis of
interacting particle systems \cite{Lig85}
and for measure-valued diffusions \cite{D} and many results are obtained
via this technique but for  population genetics models this is of particular
importance. See also \cite{EK2} for a very general approach.

Some (but not all) of these duality relations can in the case of population models be
interpreted in terms of the genealogical tree of a tagged sample and this will be the
case here. Very often for a given process various different dual processes are available
using  different spaces $E^\prime$ and  different duality functions $H$.
Depending on the application the use of different duals is in fact sometimes necessary.

Dual process representations for Wright-Fisher
processes with selection and migration were first introduced by Shiga and
Uchiyama in \cite{S}, \cite{S1} and \cite{SU}. These dual
representations were extended in Dawson and Greven \cite{DG99}
to spatial models with arbitrary type space, namely interacting
Fleming-Viot diffusion with selection and mutation and lead to a Feynman-Kac duality. The
latter can be used to establish that the martingale problem we formulated below in
(\ref{1.12}) is indeed
well-posed and to derive results on the {\em longtime behaviour} for sufficiently
large state-{\em independent} mutation.
Duals  for particle models of populations with selection have been introduced first by
Krone and Neuhauser in \cite{KN97}.
An interesting class of duality relations is studied by Athreya and Swart in \cite{AS}.

However in order to study the long-time behaviour in general we have to go beyond the
Feynman-Kac duality in \cite{DG99} and we have to develop
finer dual representations for the hierarchically interacting
system on $\Omega_N$ (or any other geographical space). It is often also useful to consider
the  nonlinear Markov process
(McKean-Vlasov process) arising as the mean-field limit $N \to\infty$ of an $N$ site exchangeable
model and also for this case one obtained a duality.

The duality we shall introduce works for
\begin{itemize}
\item
{\em multitype selection},
\item
{\em state-dependent mutation} and
\item
has a {\em historical} and {\em genealogical} interpretation and in fact extension to
tree-valued processes (see \cite{DGP}),
\item it can also be extended to multiple species and multilevel models (see \cite{D2011}).
\end{itemize}

The goal of the dual construction is to construct for a sample of $n$-individuals
from the population at a given time $t$ the probability law of the collection
\be{dgr35}
\left\{\left(\mbox{type of individual $i$, genealogical distance of pairs of individuals } (k,\ell)\right),
\; \forall\; i,k,\ell \in \I\right\}.
\ee

The key point for the analysis is that we have a {\em function-valued} dual process,
that is, in addition to a dual particle process we have a process
which is function-valued but driven
by a particle system.
In the case of a finite type space we derive a refined dual which takes values in sums
of products of indicators respectively a set-valued dual. Both these duals are very
powerful analysing the longtime behaviour.

In case of two-type and state-dependent mutation our  dual is a relative of the ancestral
selection graph of Krone and Neuhauser \cite{KN97}, but note that the latter does
not have the form of a duality with another backward Markov process.

\section{The model}
\label{s.model}

\setcounter{secnum}{\value{section}} \setcounter{equation}{0}

We now introduce  the system of interacting Fleming-Viot diffusions with mutation and selection,
where the interaction is due to  migration between the colonies indexed by the geographic space.

\subsection{Ingredients}
\label{ss.ingredients}

First we introduce the state space of the process and the basic
parameters of the stochastic evolution.\\ \noi $(\mathbf{\alpha})$ \quad
The state space of a single component (describing frequencies of
types) will be \be{1.1} \CP(\I) = \mbox{set of probability
measures on}\, \I. \ee The set of colonies (sites or components)
will be indexed by a set $\Omega_N$, which is countable and
specified in ($\beta$) below. The {\it state space} ${\cal X}$ of
the system is therefore \be{1.2}
 {\cal X} = (\CP(\I))^{\Omega_N},
\ee
with the product topology of the weak topology of probability measures
on the compact discrete set $\CI$. A typical element is written
\bea{1.3}
 X = (x_{\xi})_{\xi\in \Omega_N} \mbox{ with } x_\xi \in \CP(\I).
\eea

\medskip
\noi $(\mathbf{\beta})$   \quad The hierarchical group $\Omega_N$ indexing
the {\it colonies} of the geographic space is defined by:
\bea{1.4} \Omega_N = \Big\{\xi=
(\xi^i)_{i\in \N_0} \, | \, \xi^i \in \Z,\, 0\leq \xi^i \leq N-1,
\; \exists k_0:
  \quad  \xi^j = 0\quad   \forall \, j\geq k_0\Big\},
\eea with group operation  defined as component-wise addition
modulo $N$. A typical element of $\Omega_N$ is denoted by $\xi =
(\xi^0,\xi^1,\dots)$.  Here $N$ is a parameter with values in
$\{2,3,4, ...\}$.

Note that: \be{1.4b} \Omega_N = \bigoplusl^\infty_{i=0} Z_N, \quad
Z_N=\{0,\cdots,N-1\} \mbox{ with addition mod } (N). \ee

We also introduce a metric (actually an ultrametric) on $\Omega_N$
denoted by $d(\cdot , \cdot )$ and  defined as \be{1.5} d(\xi, \xi^\prime
) = \inf \{k \, | \, \xi^j = (\xi^\prime)^j \quad \forall \,
j\geq k\}. \ee

The use of $\Omega_N$ in population genetics in
the mathematical literature goes back to Sawyer and Felsenstein
\cite{SF}.  An element $(\xi^0, \xi^1, \xi^2, \dots)$ of
$\Omega_N$ can be thought of as the $\xi^0$-th village in the
$\xi^1$-th county in the $\xi^2$-th state .... In other words we
think of colonies as grouped and classified according to ``extent
of  neighborhood inclusion''.

Everything we state here remains valid if we replace $\Omega_N$ by a countable
group like $\Z^d$ for example.
\medskip

\noi $(\mathbf{\gamma})$ \quad  The transition kernel $a(\cdot ,\cdot )$ on
$\Omega_N \times \Omega_N$, modelling {\it migration rates}
has some specific properties.
We shall only consider \emph{homogeneous} transition kernels:
 \be{1.6}
a_N(\xi, \xi^\prime)= a_N(0, \xi^\prime - \xi) \quad \forall \xi, \xi^\prime \in \Omega_N.
\ee
\begin{example} Choose $\Omega_N$ the hierarchical group with a multi-level symmetry property:

\be{add1}
a_N(0, \xi)= \sum\limits_{k\geq j} \left(\frac{c_{k-1}}{
N^{(k-1)}}\right)\frac{1}{ N^{k}} \qquad \qquad \mbox{if~}\;
d(0,\xi) = j\geq 1,
 \ee
where
\be{1.6a} c_k> 0 \quad \forall\; k\in \N ,\quad
\sum\limits^\infty_{k=0} \frac{c_k}{ N^{k}} < \infty \quad \textrm{for all
}\; N \ge 2. \ee The kernel $a_N(\cdot ,\cdot )$ should be thought
of as follows. With rate $c_{k-1}/ N^{(k-1)}$ we  choose a
hierarchical distance $k$, and then each point within distance at
most $k$ is picked with equal probability as the new location.
Then (\ref{1.6a}) requires a finite jump rate from a given point.
\end{example}

It is often convenient to write the transition rates as
\be{agre70}
c a_N(\cdot, \cdot) \quad , \quad c \in \R^+ \mbox{ and } a_N
\mbox{ probability transition kernel on }
\I \times \I.
\ee
Then we can say that $c$ is the migration rate and
$a_N(\xi,\xi^\prime)$ the probability that a jump from
$\xi$ to $\xi^\prime$ occurs.
\bi

\noi $(\mathbf{\delta})$ \quad In addition to describe \emph{mutation} and
\emph{selection} we need two further objects. Let \be{1.6b} M
(\cdot,\cdot) \mbox{ be a probability transition kernel on } \I \times \I, \ee
modelling mutation probabilities from one type to another. Furthermore
let
\be{1.6c}
\chi(\cdot) \mbox{ be a bounded function on } \I, 0 \leq \chi (\cdot) \leq 1,
\quad 0=\mbox{ min } \chi, 1= \sup \chi,
\ee
modelling relative fitness of the different types. We have set here the
arbitrary minimum value of $\chi$ equal to zero and then in order to
uniquely specify the parameter $s$ (representing selective
intensity) we have set the maximal value of $\chi$ equal to one.
Using $\chi$ we can embed $\I$ into $[0,1]$ rather naturally
and use on it the relative topology induced by the euclidian
topology on $[0,1]$.

The case we use in \cite{DGsel} is as follows:

\begin{example} We consider a hierarchically structured set of types which allows to model
a fitness function acting on a whole hierarchy of scales. We set:

\be{mg1}
 \mathbb{I}= (\N_0\times \{1,2,\dots,M\}) \cup (0,0) \cup (\infty,1),
 \ee
 where $\N_0 =\{0\}\cup\N$ and $M$ is an integer
 satisfying $1<M<N$.

 Adding the maximal point $(\infty,1)$ allows us to
 work with a {\em compact} set of types and as a consequence a {\em compact Polish state space}
 if we embed the type set in the interval $[0,1]$ appropriately using
 the relative topology.

 Then we can write:
\be{dd1}
 \I =\bigcup_{j=0}^\infty E_j,\ee
where
\be{aGre3} E_j = \{(j,\ell); \ell = 1,\cdots,M\},\; j\in \N, \quad
E_0 = \{(0,\ell), \quad \ell=0,1,\cdots,M\}, E_\infty = \{(\infty,1)\}.
\ee
\end{example}

\subsection{Characterization of the process by a martingale
problem} \label{charact}

 We now proceed  to rigorously specify  the model by formulating the
appropriate {\em martingale problem}.

The  key ingredients for the martingale problem
with state space $E$ (Polish space) are a measure-determining sub-algebra $\CA$ of
$C_b (E,\R)$, the socalled test functions, typically called $F$ here, and a linear
operator on $C_b(E,\R)$, typically called $L$, with domain $\CA$ resulting in values
in $C_b(E,\R)$.

\beD{D.MP} {(Martingale problem)}

(a) The law $P$ on a space of $E$-valued path
for a Polish space $E$, either $D([0,\infty), E)$ or
$C([0,\infty),E)$, is a solution to the
martingale problem for $(L,\nu)$  w.r.t. $\CA$ if and only if
\be{xg103}
\left(F(X(t)) - \intl_0^t (LF)(X(s))ds\right)_{t \geq 0} \mbox{ is a
martingale under } P \mbox{ for all } F \in \CA
\ee
and
\be{xg103b}
 \CL(X(0)) = \nu.
\ee
The martingale problem is called wellposed, if the finite dimensional distributions
of $P$ are uniquely determined by the property (\ref{xg103}) and (\ref{xg103b}).

(b) In our context $E=(\CP(\I))^{\Omega_N}$. $\qquad \square$
\end{definition}

The algebra  of test functions we use here is denoted by
\be{1.11A}
\CA \subseteq C(( \CP(\I))^{\Omega_N},\R)
\ee
and is defined as follows. Given a nonnegative
bounded function $f$   on $(\I)^k$ and $\xi_1,\dots ,\xi_k \in
\Omega_N$ consider the function on $(\CM(\I))^{\Omega_N}$ (here
$\CM$ denotes finite measures) defined by: \be{1.11} F(x) =
\intl_{\I} \dots \intl_{\I} f(u_1,\dots,u_k) x_{\xi_1}(du_1) \dots
x_{\xi_k}(du_k), \ee
\[ \hspace*{3cm}\xi_i \in \Omega_N, \quad x_\xi \in \CM(\I), \quad x \in
(\CM(\I))^{\Omega_N}.
\]
Let $\CA$ be the algebra of functions generated by functions  $F$
of type (\ref{1.11}) restricted to $(\CP(\I))^{\Omega_N}$.

Our operators $L$ will be differential operators.
We define differentiation on the big space $(\CM(\I))^{\Omega_N}$,
even though we deal only with restrictions to
$(\CP(\I))^{\Omega_N}$ later on. This approach facilitates
comparison arguments. Hence we differentiate the function $F$  as
follows:
\be{1.39} \frac{\partial F}{\partial x_\xi}(x)[u] =
\lim_{\ve \to 0} \frac{F(x^{\ve,\xi,u})-F(x)}{\ve}, \ee with
\be{1.40}
x^{\ve,\xi,u} = (x^\prime_\eta)_{\eta \in \Omega_N},\quad
x^\prime_\eta = \left \{
\begin{array}{ll}
x_\eta & \eta \ne \xi\\
x_\xi + \ve \delta_u & \eta = \xi.
\end{array}
\right.
\ee
Correspondingly $\frac{\partial^2 F}{\partial x_\xi \partial
x_\xi}(x)[u,v]$ is defined as $\frac{\partial}{\partial x_\xi} \left (
\frac{\partial F}{\partial x_\xi}(x)[u]\right )(x)[v]$.
\bi

Let
\be{1.40b}
c,s,m, d>0
\ee
be parameters that represent
the relative rates of the migration, selection and mutation,
respectively the resampling rate which in biological language
represents the  inverse effective population size, finally let
$a_N(\cdot,\cdot)$ be as in (\ref{1.6}). Furthermore $\chi$ and $M$
are as in (\ref{1.6b}) and (\ref{1.6c}).
\bi

We now define a linear operator on $C_b((\CP(\I))^{\Omega_N})$
with domain $\CA$, more precisely
\be{xg102}
L : \CA \la C_b((\CP(\I))^{\Omega_N}, \R),
\ee
which is the  generator $L$ of the martingale problem and is modeling
migration, selection, mutation via drift terms and resampling as
the stochastic term (in this order), by (here $d>0$ and $c,s,m
\geq 0)$:
\bea{1.12}
(LF)(x)&  =\sum_{\xi \in \Omega _N} \Bigg[ &c
\suml_{\xi^\prime \in \Omega_N}a_N(\xi ,\xi ^{\prime })\intl_{\I}
\frac{\partial F(x)}{\partial x_\xi }(u)
(x_{\xi ^{\prime }}-x_\xi )(du)\\
&&+\; s \intl_{\I}\left\{  \frac{\partial F(x)}
    {\partial x_\xi }
    (u)\left(\chi(u)-\int_{\I} \chi(w)x_\xi(dw) \right)\right\}x_\xi(du)
\nonumber\\
&&+\; m\intl_{\I}  \left \{ \intl_{\I}
    \frac{\partial F(x)}{\partial x_\xi }
    (v)M(u,dv)
 -  \frac{\partial F(x)}{\partial x_\xi }(u)\right \}
    x_\xi (du)
    \nonumber \\
&&+ \;d\intl_{\I} \intl_{\I}
    \frac{\partial ^2F(x)}{\partial x_\xi \partial x_\xi }
    (u,v)Q_{x_\xi }(du,dv)  \Bigg ], \qquad x \in (\CP(\I))^{\Omega_N}, \nonumber
\eea where
\be{1.13} Q_x(du,dv) = x(du) \delta_u(dv) - x(du)x(dv).
\ee

It has been proved in \cite[99]{DG99} that a model of the type as above is
well defined:

\noi
\beT{T0} \quad (Existence and Uniqueness)\\
Let $\nu $ be a probability measure on $(\CP({\I}))^{\Omega_N}$
specifying the initial state which is independent of the
evolution.

(a) Then the $(L;\nu)$-martingale problem w.r.t. $\CA$, on
the space $C([0,\infty)\; $ (cf. Definition \ref{D.MP}),  $(\CP(\I))^{\Omega_N})$ is
well-posed.
For fixed value of the parameter $N$ the resulting
canonical stochastic process is denoted
\be{xg104}
(X_t^N)_{t\ge 0}.
\ee

(b) The solution defines a strong Markov process with the Feller property.
$\qquad \square$
\end{theorem}

In the sequel we shall consider  often initial states which are
either deterministic or  random with law  $\nu$ in $\CP
((\CP(\I))^{\Omega_N})$ which satisfies the following two properties:
\be{GIC1} \nu \mbox{ is
invariant w.r.t. spatial shift on  } \Omega_N, \ee \be{GIC2} \nu
\mbox{ is shift ergodic}. \ee

\section{Function-valued dual}
\label{basictool} \setcounter{secnum}{\value{section}
\setcounter{equation}{0}
\renewcommand{\theequation}{\mbox{\arabic{secnum}.\arabic{equation}}}}

In this section we first recall the Feynman-Kac duality from \cite{DG99}, next
we develop the {\em new} dual representations for the interacting system.
This allows us in particular to establish
the ergodic theorem for the mean-field limit with
state-\emph{independent} mutation component.

We then develop a further modified dual that
covers state-\emph{dependent} mutation in order  to prove for example ergodic theorems or results
describing the transition from one quasi-equilibrium to the next (see \cite{DGsel}).
This latter dual representation is very flexible and allows for various
refinements useful for specific purposes. In particular the multi-scale
analysis requires some new arguments due to the interplay between selection
and mutation.
Each of these arguments requires some specific features of the
dual representation which requires some slight modifications of the
dual mechanism that are developed here in this section systematically.

In \cite{GLW} and \cite{GPWmp} it is shown that for the neutral model with migration the
duality for the marginal distribution arises from a dual representation of
the genealogical  process of the underlying particle process. Roughly speaking the
dual dynamics generates the law of the genealogical tree of a randomly sampled finite tagged
population among the total population at a given time $t$ in terms of the genealogical tree
associated with a coalescent. This raises the question to what extent
this can be extended to processes with mutation and selection. This remains to
be investigated and will be discussed later on in Subsection \ref{histint}.
On the level of genealogies such a dual is used in \cite{DGP}.

{\bf Outline of Section \ref{basictool}}
We begin in Subsection \ref{dualp} to introduce the duality functions and
the dual process and give in Subsection \ref{FKd} the basic duality
relations on which one can base all the modifications and refinements
we present later on.
In Subsection \ref{modified} we give
special versions of the dual representation of the mutation
 useful for the different possible applications of the dual we
have a need for.
In Subsection \ref{histint} we discuss the
historical meaning of the dual representation.
In Subsection \ref{ss.refine} we derive a very important refinement of the duality relation for models
with finitely many types.
This allows also to obtain a historical interpretation of the dual
which we discuss in Subsubsection \ref{sss.histint}.

\subsection{The dual process and duality function for the interacting system}
\label{dualp}

In this section we introduce the basic ingredients of the duality theory.
In order to obtain a  dual representation of a
$(\mathcal{P}(\mathbb{I}))^{\Omega_N}$-valued process $(X_t)_{t
\ge 0}$ we need two ingredients, (1) a \emph{family of functions} on the state space that is
measure-determining from which we construct the duality function $H(\cdot, \cdot)$
and (2) the stochastic dynamics of a \emph{dual process}
whose state space is the set which labels the family in the previous point.

{\bf (1) Duality functions}

We begin with the first point.
Recall that in order to determine the marginal distribution at
time $t$ of the law $P_{X_0}$ of the system $X(t) =
(x_\xi(t))_{\xi \in \Omega_N}$, with initial point $X(0)=X_0 \in
(\CP(\I))^{\Omega_N}$,  it suffices to calculate for every fixed time $t$ a suitable form of {\em mixed
moments}, i.e. we have to determine the following
quantities parametrized by $(\xi_0,f)$
(the initial state of the intended dual dynamics) with
$\xi_0$ a vector of geographic positions and a basic function $f$
on finite products of the type space $\I$:
 \be{3.2}
  F_t((\xi_0, f), X_0) : = E_{X_0} \left [
  \intl_{\I} \dots \intl_{\I} f(u_1,\dots,u_n )
  x_{\xi_0(1)}(t)(du_1) \dots x_{\xi_0(n)}(t)(du_n)\right ], \ee
  for
  all $n \in \N,\; \xi_0 = (\xi_0(1),\dots,\xi_0(n)) \in
  (\Omega_N)^n$ and  $f \in L_\infty((\I)^n, \R)$.
This means we want to consider the bivariate function
\be{add105}
H(X,(\xi,f)) = \intl_{I^n} f(u_{1}, \cdots, u_{n})
x_{\xi(1)} (du_1) \cdots x_{\xi(n)}) (du_n)
\ee
to obtain a duality. This would however only allow for a dual for the migration part.
Namely introducing the dual process we shall see
that we have to enrich the second argument further, to be able to obtain a Markov
process as a dual, namely resampling and selection require to introduce
partitions of the variables.

  The class given above contains the functions which we will use
  for our dual representations.
  In fact it suffices to take smaller, more convenient sets of test functions $f$:
  \be{productform} f(u_1,\dots,u_n )=\prod_{i=1}^n f_{i}(u_i),
  \quad f_i \in L_\infty (\I). \ee
  Since in our case $\mathbb{I}$ is
  countable, it would even suffice to take functions
  \be{3.2b}
  f_{j}(u)=1_{j}(u) \mbox{ or } f_j(u)=1_{A_j}(u),
\ee
where $1_{j}$ is the indicator function of $j\in\mathbb{I}$ and
in case of finite type sets, $1_A$ with $A$ some sets of types
for example like $\{u| \chi(u) \geq e_k\}$ if $0=e_0 < \cdots < e_\ell=1$
are the fitness values, which will
  be a very useful choice later on. \bi

\begin{remark} \quad (Notational convention).

In this section we
write
\be{coll1}
X_t \mbox{ instead of } X(t)
\ee
on the level of the collection (but
we still write $x_\xi(t)$ for components). \end{remark} \bi

{\bf (2) Dual process: Spatial coalescent with births driving a function-valued process}

Now we come to the second point, the dual process on a suitable state space.
We observe  here that due to the interaction  of selection and mutation the state space of the
dual process is more complicated than for most  particle systems or interacting diffusions.
We proceed in four steps.
\sm

{\bf Step 1} {\em (The state of the dual process)}

We begin this point by reviewing the  dual process used
in \cite{DG99} which involves a {\em function-valued process}
$(\CF_t)_{t \geq 0}$ driven by a {\em finite ordered particle system}
$(\eta_t)_{t \geq 0}$. Indeed in
\cite{DG99}  the quantities on the r.h.s. of (\ref{3.2}) were
expressed in terms of the evolution semigroup corresponding to a
\be{e.markov}
\mbox{ bivariate Markovian stochastic process } (\eta_t,\mathcal{F}_t)_{t \geq 0}
\ee
consisting of
\begin{itemize}
\item a process $\eta$ involving a random set of finitely many ordered individuals,
\item a function-valued process $\CF$.
\end{itemize}
This is the dual process to be discussed in the sequel where $\eta_t$ takes
values in the set of {\em marked partitions} of a {\em random set of individuals} and $\CF_t$ is a
{\em function} defined on {\em products of the type set $\I$}.
The joint evolution of these   is such that the dynamic of $\CF_t$ is driven by
the process $\eta_t$ which evolves autonomously.
\sm

We will use two main types of such processes which are denoted by
\be{an3}
( \eta_t, \CF_t), \mbox{ respectively } ( \eta_t, \CF^+_t) \mbox{ or its variant }
(\eta_t, \CG^+_t)_{t \geq 0},
\ee
where the first component, $\eta_t$ is
a particle system described formally below and the second component $\CF_t, \CF^+_t$,
takes values in the  space of functions on $\I^\N$ which are measureable, bounded and depend
only on finitely many components which we can identify with
\be{3.3}
\bigcup^\infty_{m=1} L_\infty((\I)^m),
\mbox{ resp. } \bigcup^\infty_{m=1} L_\infty^+ ((\I)^m).
\ee

This  process, called the
\emph{function-valued dual process}, will allow us in Subsections
\ref{FKd} and \ref{modified}
to represent
for every given $X_0$ the function $F_t$ from (\ref{3.2}) as the expectation of an
appropriate functional over the evolution of $( \eta_t, \CF_t)$ or
$( \eta_t, \CF^+_t)$ respectively $(\eta_t,\CF^+_t)$. We construct it here in such a way (for example its order structure)
that it can be refined later better for some special applications.
\bi

The process $(\eta_t, \CF_t)_{t \geq 0}$ is constructed from the following four ingredients:

\begin{itemize}
\item $N_t$ is a non-decreasing $\N$-valued process. $N_t$ is the number of {\em individuals} present in the dual process  with $N_0=n$, the number of initially tagged individuals,
and with $N_t-N_0\geq 0$ given by the number of individuals born during the interval $(0,t]$.

\item $\zeta_t =\{1,\cdots,N_t\}$ is an {\em ordered particle system} where the individuals are given an assigned order and the remaining particles are ordered by time of birth.

\item \begin{tabular}{ll} \label{dualcomp}
\hspace{-0.2cm}$\eta_t=$ & $(\zeta_t, \pi_t, \xi_t)$ is a trivariate process consisting of the above $\zeta$ and\\
& $- \pi_t$:  {\em partition} $(\pi^1_t,\cdots,\pi^{|\pi_t|}_t)$ of $
\zeta_t$, i.e. an ordered family of subsets, where the \\&\qquad
{\em index} of a partition element is the smallest element (in the ordering of individuals)\\&\qquad of the partition element.\\
& $- \xi_t: \pi_t \longrightarrow \Omega^{|\pi_t|}_N$,  giving {\em locations} of the
partition elements. Here $\xi_0$ is the vector \\&\qquad of the prescribed space points
where the initially tagged particles sit.
\end{tabular}
\item $\CF_t$ is for given  $\eta_t=( \zeta_t, \pi_t,\xi_t)$ is a function in
$L_\infty (\I^{|\pi_t|})$.
\end{itemize}

Therefore the state of $\eta$ is an element in the set
\be{dg30.1}
\CS = \bigcup\limits^\infty_{m=1} \{1,2, \cdots, m\}
\times \mbox{Part}^<(\{1,\cdots,m\}) \times \{\xi: \mbox{Part}^< (\{1,\cdots,m\})
\la \Omega_N\},
\ee
where Part$^< (A)$ denotes the set of ordered partitions of a set $A \subseteq \N$
and the set $\{1,\cdots, m\}$ is equipped with the natural order.

It is often convenient to order the partition elements according to their
indices
\be{dg30}
\mbox{{\em index} } (\pi^i_t)= min (k \in \N : k \in \pi^i_t)
\ee
and then assign them ordered
\be{l1}
\mbox{ {\em labels} } 1,2,3, \cdots, |\pi_t|.
\ee

\begin{remark} \quad
Note that
\be{3.2b2}
N_t = |\pi^1_t| +\cdots + |\pi_t^{|\pi_t|}|.
\ee
\end{remark}

\begin{remark} \label{R.eta2} \quad
The process $\eta$ can actually be defined starting with countably many individuals located in $\Omega_N$. This is due to the quadratic death rate mechanism at each site which implies that the number of individuals at any site will have jumped down into $\N$ by time $t$ for any $t >0$.
\end{remark}

\begin{remark} \quad
We view the configuration of $(\eta_s)_{s \leq t}$ as the
description of the genealogy of the sample drawn at time $t$
observed back at time $t-s$. We interpret $N_0=n$ as the sample size at time $t$ and
$N_s$ as the population at time $t-s$ which can potentially influence the type of the
individuals in the sample at time $t$. Then $(\xi_s)_{s \leq t}$ give the
ancestral path of the individuals which are potentially relevant
for the $n$-sample at time $t$. The partition $\pi_s$ describes the groups of
individuals at time $t$ which have a common ancestor at time $t-s$.
See Subsection \ref{histint} and Subsubsection \ref{sss.histint} for more information.
\end{remark}
\sm

\begin{remark}\label{R.classic} \quad
We can now formulate our dual process  into the classical scheme of duality theory (see \ref{add100}). We
write  $X^\prime$ for the pair $(\eta,\CF)$,
\be{dg30.2}
E = (\CP(\I))^{\Omega_N} \quad , \quad E^\prime = \{X^\prime \in \CS \times (\bigcup^\infty_{m=1}
L_\infty (\I^m))\; | m =|\pi|\}.
\ee
If we use the standard topology on the set $E$, it  is a Polish space.

$E^\prime$ can be embedded in a complete metric space which is a subset of $\wt{\CS}$
based on ordered partitions of $\N$ (rather than finite sets) and the set  of functions in $L_
\infty(\I^\N)$ which can be approximated by functions of finitely many variables.
If we want to achieve separability of $E^\prime$ one needs to
assume more on $\chi$ and $H$ to be able to restrict the function space to continuous functions
on $\I$. However one should note that the measures appearing as states of $X$ are
always at most countably supported for positive time.

Then define a map
$H:E\times E^\prime \to \R$ by
\be{dg30.3}
H(X, (\eta,\CF)) = \intl_{\I} \cdots  \intl_{\I} \CF (u_{1}, \cdots, u_{|\pi|})
x_{\xi(1)} (du_1) \cdots x_{\xi(|\pi_t|)} (d u_m).
\ee
Then the collections of bounded measureable functions
\be{dg30.4}
\{H(\cdot,(\eta,\CF)),  (\eta,\CF) \in E^\prime\}, \quad
(\{H(X, \cdot), \quad X \in (\CP(\I))^{\Omega_N} \})
\ee
are measure-determining on $(E,\CB(E))$ respectively $(E^\prime, \CB (E^\prime))$.
The second property is not always satisfied in the duality theory and typically
not needed in applications. Then $H(\cdot,\cdot)$ is called a duality function
for the pair $(E,E^\prime)$ and $E$ and $E^\prime$ are legitimate state spaces
for Markov processes, recall Remark \ref{R.eta2}. (See \cite{Lig85} and \cite{EK2}
for detailed expositions of duality theory).
\end{remark}

We describe the dynamics of the process in two steps, first we
give the autonomous dynamics of $(\eta_t)_{t \geq 0}$ and then construct
$(\CF_t)_{t \geq 0}$ respectively $(\CF^+_t)_{t \geq 0}$ or $(\CG^+_t)_{t \geq 0}$ given a
realization of the process $\eta$. In subsequent sections we
introduce various modifications of this tailored for specific
purposes.\bi

\noi
{\bf Step 2} {\em (Dynamics of dual particle system)}

The dynamics of $\{\eta_t\}$ is that of a pure Markov jump process with the following transition mechanisms which correspond to resampling,
migration and selection in the original model
(in this order):
\begin{itemize}
\item \emph{Coalescence} of partition elements:  any pair of  partition
elements which are at the same site  coalesce at a constant
rate; the resulting enlarged partition element is given the
smaller of the indices of the coalescing elements as index and the
remaining elements are reordered to preserve their original order.
\item \emph{Migration} of partition elements: for $j\in
\{1,2,\dots,|\pi_t|\}$, the partition element $j$ can migrate
after an exponential waiting time to another site in $\Omega_N$ which means we have a jump of
$\xi_t(j)$.
\item \emph{Birth} of new individuals by each partition element
after an exponential waiting time: if there are currently $N_t$ individuals then a newly
born individual is given index $N_{t-}+1$ and it forms a {\em partition
element consisting of one individual}  and this partition element is given a label
$|\pi_{t-}|+1$. Its location mark is the same as that of its parent
individual  at the time of birth,
\item {\em Independence}: all the transitions and waiting times described in the previous
points occur independently of each other.
\end{itemize}

Thus the  process $( \eta_t)_{t \ge 0}$ has the form
\be{3.4}
\eta_t = (\zeta_t, \; \pi_t, (\xi_t(1), \xi_t(2), \cdots, \xi_t(|\pi_t|))),
\ee
where $\xi_t(j) \in\Omega_N,\;j=1,\dots,|\pi_t|,$ and $\pi_t$ is an
ordered partition (ordered tuple of subsets) of the current basic set of the form
$\{1,2,\dots,N_t\}$ where $N_t$ also grows  as a random process,
more precisely as a pure birth process.

We order partition elements by their smallest elements and then
{\em label} them by $1,2,3,\dots$ in increasing order. Every partition
element (i.e. the one with label $i$) has at every time
a location in $\Omega_N$, namely the i-th partition
element has location $\xi_t(i)$. (The interpretation is that we
have $N_t$ individuals which are grouped in $|\pi_t|$-subsets and
each of these subsets has a position in $\Omega_N$).

Furthermore denote by
\be{3.2c}
\pi_t (1), \cdots, \pi_t (|\pi_t|)\ee
the index of the first, second etc. partition element.
In other words the map gives the index of a partition element
(the smallest individual number it contains) of a label
which specifies its current rank in the order.

For our concrete situation we need to specify in addition the
parameters appearing in the above description, this means that we define
the dual as follows.

\beD{D.Could}{(First component of the dual process: $\eta_t$)}

(a) The initial state $\eta_0$ is of the form $(\zeta_0, \pi_0,\xi_0)$ with:
\be{dg31}
\zeta_0=\{1,2,\cdots,n\}, \quad \pi_0 =\{\{1\}, \cdots,\{n\}\}
\ee
\be{dg32}
\xi_0 \in (\Omega_N)^n.
\ee
(b)
The evolution of $( \eta_t)$ is defined as follows:
\begin{itemize}
\item[(i)] each pair of partition elements which
 occupy the same site in $\Omega_N$ \emph{coalesces} during the
 joint occupancy of a location into one partition element after a \emph{rate $d$}
 exponential waiting time,
 \item[(ii)] every partition element performs, independent of the
 other partition elements, a \emph{continuous time random walk} on
 $\Omega_N$ with transition \emph{rates $c \bar a(\cdot, \cdot)$} (see
 (\ref{1.6})), with $\bar a(\xi,\xi^\prime)=a(\xi^\prime,\xi)$,
 \item[(iii)] after a rate $s$ exponential waiting time each partition element
 gives \emph{birth} to a new
 particle   which forms a new (single particle)
 partition element at its location, and this new partition element
 is given as label $|\pi_{t-}|+1$ and the new particle the index $N_{t-}+1$.
 \end{itemize}
 All the above exponential waiting times are \emph{independent} of each other.
 $\qquad \square$
 \end{definition}
 Note that this process is well-defined since the total number of
 partition elements is stochastically dominated by an ordinary linear birth
 process which is well-defined for all times. \bi

\noi
{\bf Step 3} {\em (Dynamics of function-valued part of dual process: $\CF_t, \CF^+_t, \CG^+_t$)}

To complete the specification of the dual dynamics we
want to define a bivariate process $(\eta_t, \CF_t)$ resp.
   $(\eta_t, \CF^+_t)$ associated with (\ref{1.12}).
   Hence we now describe the evolution of the process$(\mathcal{F}_t)_{t\geq 0}$
   conditioned on a realization of the process $(\eta_t)_{t\geq 0}$.
   For a given path of the first component  $( \eta_t)_{t \ge 0}$ the second component
   is a \emph{function-valued process} $(\CF_t)_{t \ge 0}$,
   respectively $(\CF^+_t)_{t \geq 0}, (\CG^+_t)_{t \geq 0}$ with values in
   $\bigcup^\infty_1 L_\infty(\I^m)$, respectively $\bigcup^\infty_1
   L^+_\infty(\I^m)$.

The evolution of $\CF_t$, respectively $\CF^+_t$, starts in the state
\be{evol1}
   \CF_0 = \CF^+_0 = f,
\ee
with a bounded function $f$ of
$|\pi_0|$-variables, a parameter we denoted by $n$, the variables running in $\I$.

The evolution of $(\CF_t)_{t \ge 0}$ involves three mechanisms, corresponding to the
\emph{resampling}, the  \emph{selection} mechanism and the \emph{mutation} mechanism.
We now describe separately these mechanisms.

\beD{D.Fuvadure} {(Conditioned evolution of $\CF, \CF^+$: the coalescence mechanism)}

 If a coalescence of two partition elements occurs,
 then the corresponding variables of $\CF_t$ are set equal to the
 variable indexing the partition element, (here $\wh u_j$ denotes an omitted
 variable), i.e. for $\CF_t=g$ we have the transition
 \bea{ga3}
 g(u_1,\cdots,u_i,\cdots,u_j,\cdots,u_m)
 & \longrightarrow & \wh g(u_1,\cdots,\wh u_j,\cdots,u_{m}) \\
 && =  g(u_1,\cdots,u_i,\cdots,u_{j-1}, u_i, u_{j+1},\cdots,u_{m}), \nonumber
 \eea
 so that the function changes from an element of $L_\infty(\I^m)$ to one of
 $L_\infty(\I^{m-1}). \qquad \square$
 \end{definition}

 \begin{remark}
  \quad We could also work with
 \be{ga4}
 \wt \CF_t (u_{\wt \pi_t(1)}, \cdots, u_{\wt \pi_t(N_t)}) \mbox{ instead of }
 \CF_t (u_{1},\cdots,u_{(|\pi_t|)}).
 \ee
This means in (\ref{ga3}) we view the r.h.s. as element of $L_\infty(\I^m)$.
Then the state space $\wt E^\prime$ is obtained by replacing in (\ref{dg30.2}) the
restriction by $m=|\zeta|$ and using the duality function $\wt H(\cdot,\cdot)$ given by
\be{add127}
\wt H (X,(\eta,\wt{\CF})) = \intl_{I^{|\pi|}} \wt{\CF}(u_{\wt \pi(1)}, \cdots, u_{\wt \pi (N_t)}))
x_{\xi(1)} (du_1) \cdots x_{\xi(| \pi |)} (du_{|\pi|}).
\ee
The form in (\ref{ga4}) codes some historical information which is lost in the form in (\ref{dg30.3}).
 On the other hand expressions in the duality relation become often simpler
 in the reduced description.
 \end{remark}

 We  continue with the selection mechanism and its dual.
 There are several alternate ways to incorporate the effects of
 selection into the dual.  In particular some versions involve
 signed function-valued processes and others use only
 non-negative-valued functions.  Some versions involve a
 Feynman-Kac factor in the representation and others do not. We now
 describe those versions of the dual process that will be elaborated
 on as required in later sections.

\beD{D.Funkvdu} {(Conditioned evolution of $\CF$ and $\CF^+, \CG^+$: the selection
mechanisms)}

\begin{itemize}

\item[(i)] If a birth occurs in the process $(\eta_t)$ due to the
partition element  $\pi_i,\;i\in\{1,\dots,|\pi_t|\}$ , then
for $\CF_{t-}=g$ the following transition occurs from an element in
$L_\infty((\I)^m)$ to elements in $L_\infty((\I)^{m+1})$:
\be{3.6}
g(u_1,\dots,u_m) \longrightarrow \chi(u_i)g(u_1,\dots,u_m) -
\chi(u_{m+1}) g(u_1,\dots,u_m),
\ee
where the new variable is associated with the partition element of the newly born individual.
\item[($i^\prime$)] For $\CF^+_t$ the transition (\ref{3.6}) is replaced by (provided that
$\chi$ satisfies $0 \leq \chi \leq 1$):
\be{3.6b}
g(u_1,\cdots,u_m) \longrightarrow \wh g(u_1,\dots,u_{m+1})=(\chi
(u_i) + (1-\chi (u_{m+1}))) g(u_1,\cdots,u_m),
\ee
where the new variable is associated with the newly born individual.
In particular
\be{ga2}
 \Vert \wh g\Vert_{\infty}\leq 2\Vert g\Vert_\infty, \mbox{ and } \wh g \geq 0,
  \mbox{ if } g \geq 0.
   \ee
\item[($i^{\prime\prime}$)] For $(\CG^+_t)_{t \geq 0}$ we use the following transition:
\be{ga2b}\begin{array}{l}
g(u_1, \cdots, u_m) \la (\chi (u_i)1(u_{m+1}) g(u_1, \cdots, u_m)\\
\hspace{5cm}
+ (1- \chi(u_i)) g(u_1, \cdots, u_{i-1}, u_{m+1}, u_{i+1}, \cdots, u_m),
\end{array}
\ee
in which case
\be{ga2c}
\|\wh g\|_\infty \leq \|g\|_\infty. \qquad \square
\ee
\end{itemize}
\end{definition}

     \begin{remark} \quad The mechanism (i) will lead to a
     signed-function-valued process and requires a Feynman-Kac factor,
     while the mechanism ($i^\prime$) leads to a non-negative
     function-valued process and in this case NO Feynman-Kac factor is
     needed.
The mechanism ($i^{\prime\prime}$) induces the same change in the duality function $H$ as ($i^\prime$) and
therefore leads also to a duality function.
     \end{remark}\sm

\begin{remark} \quad  We note that the dynamic ($i^{\prime\prime}$) disconnects the strict parallel
between the variables in the function and the individuals in the process $\eta$ which corresponds to the
fact that an individual in the sample may be replaced by one from outside the
sample or vice versa.
\end{remark}

We now introduce some basic objects used to incorporate the
mutation mechanism into the dual. Let  $\CG^{\rm{mut}}_t$ denote the
semigroup on $\bigoplus\limits_{n \in \N} L_\infty(\I^n)$ induced by independent
copies of   the mutation
process acting on $L_\infty(\I)$, i.e. the process we obtain from
(\ref{1.12}) by letting the geographic space (i.e. the index set for the components consist
of one point and putting all other coefficients $d, s, r, c$ in
(\ref{1.12}) equal to 0).  We can represent $\CG^{\rm{mut}}_t$ in terms of
the process
$\{M^\ast_t\}_{t\geq 0}$  which is a Markov pure jump process with jumps $L_\infty(\I^n)\to L_\infty(\I^n)$:
\be{} f\to M_jf\quad\text{at rate }m\;\text{ for each }j=1,\dots,n
\ee
with $M_j$ given by (\ref{an4}) below.
Then
\be{an4ext}
\CG^{\rm{mut}}_t(g) = E[ M^\ast_t(g)].
\ee
 is a (deterministic) linear semigroup on
$\bigoplus_{n \in \N}  L_\infty (\I^n)$ (with exactly that domain) with generator
\be{ang5}
\mathcal{M}^\ast=m \sum_j (M_j-I),
\ee where $M_j$ acting on $L_\infty(\I^n)$ is for $j \in \{1,\cdots,n\}$ given by
\be{an4}
M_jg(u_1,\dots,u_j,\dots,u_n)= \intl_{\I}
g(u_1,\dots, u_{j-1},v,u_{j+1},\dots,u_n)M(u_j,dv).\ee
This allows to define:

\beD{D.Funkvdu2} {(Conditioned evolution of $\CF,\CF^+, \CG^+_t$: the
mutation mechanism)}

If no transition occurs in  $(\eta_t)_{t \ge 0}$ then $\CF_t$,
respectively $\CF^+_t, \CG^+_t$ follows the deterministic evolution given
by the semigroup $ \CG^{\rm{mut}}_t$ acting on $L_\infty(\mathbb{I}^m)$
provided that currently $\CF_t$, respectively ($\CF^+_t$ or $\CG^+_t$) is
a function of $m$ variables, meaning at the last transition in
$(\eta_s)$ up to time $t$ became an element of $L_\infty((\I)^m) \qquad \square$. \sm

  \end{definition}

   This completes the evolution mechanism driving the process
   $(\CF_t)_{t \geq 0}$ respectively $(\CF^+_t)_{t \geq 0}$ or
   $(\CG^+_t)_{t \geq 0}$ for given $\eta$.
   In this fashion, the function-valued processes starting with an element $f \in
   L_\infty(\I^n)$ for some $n$ are uniquely defined for every given
   path for the process $\eta$.
   \bi

\noi {\bf Step 4} {\em (Definition bivariate dual process)}

Combining Step 1 and Step 2 we have defined  uniquely for every (ordered) set of $n$
individuals with $n$ locations in $\Omega_N$ and a function $f_0$ in
$L^\infty(\I^n)$
a bivariate Markov process $(\eta_t, \CF_t)_{t \ge 0}$,
respectively $ (\eta_t,\CF_t^+)_{t \ge 0}$ or $(\eta_t,\CG^+_t)_{t \geq 0}$  in which the jumps
occur only at the times of jumps in the pure Markov jump process
$(\eta_t)_{t \ge 0}$ and changes in the function-valued part
occur deterministically in a continuous way (action of
mutation semigroup).

\begin{remark} \quad
One advantage of working with the process $(\eta_t,\mathcal{F}_t)$ (or
$(\eta_t, \CF^+_t)$) over working with the original process is, that it
can be explicitly constructed and only a finite number of
random transitions occur in a finite time interval whereas  changes occur
instantaneously in countably many components of the process $X$  in any
finite time interval and are furthermore diffusive changes.
\end{remark}

\subsection{Duality relation for interacting systems}
\label{FKd}

The following is the basic {\em Feynman-Kac duality } relation
that was established in \cite{DG99} for the general class of interacting
Fleming-Viot processes with selection and mutation from which
we derive the {\em new duality relation} in Theorem \ref{T.orddual} below
without Feynman-Kac term which we use in this paper.

\beP{P3.1} \quad (Duality relation - signed  with Feynman-Kac dual)\\
Let $(X_t)_{t \ge 0}$ be a solution of the $(L,X_0)$-martingale
problem with $L$ as in (\ref{1.12}),
$X_0=(x_\xi)_{\xi \in \Omega_N}$. Choose $\xi_0 \in
(\Omega_N)^n$ and $f \in L_\infty((\I)^n)$ for some $n \in \N$.
(Recall (\ref{3.3}), (\ref{3.4}) for notation).

Assume that $t_0$ is such that:
\be{3.1.1}
\label{expcondition} E\left(\exp\left(s\int_0^{t_0}|\pi_r|dr\right
)\right)<\infty.\end{equation}
Then for $0\leq t\leq t_0$, $(\eta_t, \CF_t)$ is the Feynman-Kac
dual of $(X_t)$, that is:
\bea{3.8}
\lefteqn{F_t(( \eta_0,f),
X_0) =  E_{( \eta_0, \CF_0)}\Bigg \{ \Bigg [ \exp(s \intl^t_0
|\pi_r|dr) \Bigg ] \cdot  }\\
[1ex] & & \hspace*{3cm}  \cdot \Bigg [\intl_{\I} \dots \intl_{\I}
\CF_t(u_1,\dots,u_{|\pi_t|}) x_{\xi_t(1)} (du_1)\dots
x_{\xi_t(|\pi_t|)}(du_{|\pi_t|}) \Bigg ] \Bigg\}, \nonumber \eea
where the initial state $(\eta_0, \CF_0)$ is given by

\bea{3.9}
 \eta_0 & =&  [\{1,\dots,n\},\xi(1), \cdots, \xi
 (n);(\{1\},\{2\},\dots,\{n\})],\\
 &&\quad\qquad \xi (i) \in \Omega_N \mbox{ for } i=1,\cdots,n; \quad n \in \N, \nonumber \\
 \CF_0 & =&  f. \qquad \square \nonumber \eea
 \end{proposition}

 \begin{remark} In \cite{DG99} it was shown that there exists $t_0 > 0$ for
 which (\ref{3.1.1}) is satisfied.
 \end{remark}
 \bi

The unpleasant features of the above dual are the exponential term and the signed function, which
make it  difficult to analyse for $t \to \infty$. However
 if we have a fitness function $\chi$ satisfying $0 \leq \chi
 \leq 1$ (which because of scaling properties of the dynamic really only means
 that we do not allow {\em unbounded} fitness but cover all other cases)
 then we can obtain the following duality relation that
 does {\em not} involve a Feynman-Kac factor and preserves the {\em positivity} of functions:

\beT{T.orddual}(Duality relation - non-negative)

With the notation and assumptions as in Proposition \ref{P3.1}
(except \ref{3.1.1}) we get for $\chi$ with
\be{3.9.0}
0 \leq \chi \leq 1
\ee
that for all $t \in [0,\infty)$ we have:
\bea{3.9.1}
F_t (( \eta_0,f)),X_0) &=& E_{(\eta_0,\CF_0^+)}
\left [ \intl_\I \cdots
\intl_\I  \CF^+_t (u_1, \cdots,u_{|\pi_t|}) x_{\xi_t(1)}
(du_1) \cdots x_{\xi_t (|\pi_t|)} (du_{|\pi_t|}) \right ]\\[1ex]
&=& E_{(\eta_0,\CG_0^+)}
  \left [ \intl_\I \cdots \intl_\I \CG_t^+ (u_1, \cdots, u_{|\pi_t|})
	  x_{\xi_t(1)} (du_1) \cdots x_{\xi_t (|\pi_t|)} (du_{|\pi_t|}) \right]. \nonumber
\eea
 Moreover, $\CF^+_t$ is always non-negative if $\CF^+_0 =f\geq 0$ , similarly for
 $\CG_t^+$ if $\CG^+_0 =f\geq 0$.\hfill $\square$
 \end{theorem}

\begin{remark} \quad The duality relations we develop here for $\Omega_N$ work
in fact for every geographic space which forms a countable group and for every
migration mechanism induced by a random walk on that group. With some
care we can even pass to migration given by a general Markov chain.
\end{remark}

\begin{remark} \quad The duality relations
work for general type space for example a continuum $[0,1]$ as well,
the special structure of $\I$ has not been used.
 \end{remark}

We note that the duality relation given above is of the general form considered
in duality theory but the state space is more complex than in the classical
cases for Fisher-Wright diffusions or measure-valued processes as demonstrated by the
following remark.

\begin{remark} \quad
In the language of classical duality theory as described in Remark \ref{R.classic}
we have abbreviating $X^\prime_t =(\eta_t,\CF_t)$, the duality relation
\be{dg30.5}
E_{X_0} [H(X_t, X^\prime_0)] = E_{X^\prime_0} [H(X_0,X^\prime_t)], \quad \forall \; t \geq 0
\ee
whenever $(X_0,X^\prime_0) \in E \times E^\prime $.
Therefore the l.h.s. which is the object of interest can be calculated in terms
of the r.h.s. involving a process of a simpler nature.

For  Feller processes
on compact state spaces with generators $G_X$ respectively $G_{X^\prime}$ , respectively,  then the  duality relation follows from the generator relation
\be{dg30.6}
(G_X H(\cdot, X^\prime_0)) (X_0) = (G_Y H(X_0, \cdot))(X^\prime_0), \mbox{ for all } (X_0, X^\prime_0)
\mbox{ in } E \times E^\prime
\ee
(see e.g. \cite{Lig85}, Chapt. 2).
On non-compact spaces some integrability properties have to be verified in addition
(see \cite{EK2}, Chapt. 4).

\end{remark}

The duality relation above can be extended as follows.

 \begin{remark} \quad It is often necessary to consider moments of the
 time-space process and then the process is observed say at a time $t$ and
 $t+s$. In this case one considers the dual particle system where the
 particles corresponding to the moment at the earlier time $t$ are
 activated in the dual process only at times $s$ and the dual is
 evaluated after evolving till time $t+s$. This observation has first
 been used in the voter model and is quite useful for us later on in
 some calculations.
 \end{remark}

The duality relation in Theorem \ref{T.orddual} can be related to the  previous FK-duality
 in \cite{DG99} quoted in Proposition \ref{P3.1} if we observe that the action of the selection part of the
 generator of $X$ can be written in two ways, namely as (assuming that $\chi$ is
 bounded by 1 and in the expression $\chi g$ below taking $ \chi$
 to be a function of \emph{one} of the variables of $g$):
 \be{3.9.2}
 \{(\chi g - \chi \otimes g) - g\} + g=(\chi g +
 (1-\chi)\otimes g)-g. \ee
 The first expression, that is, $\{(\chi g
 - \chi \otimes g) - g\} + g$ corresponds to a jump defining the
 transition of $\CF_t$ (the term in $\{\}$) plus a Feynman-Kac term $g$ as we have in
 the duality in (\ref{3.8}).
 The second expression, that is, $(\chi g + (1-\chi)\otimes g)-g$
 corresponds to a jump describing $\CF^+_t$ which now yields as new
 function-valued state or function consisting of two
 summands but does not require a Feynman-Kac factor.
Since integrability issues are involved we now give a detailed proof below.

\begin{proof}{\bf of Theorem \ref{T.orddual}}

We have to do two things, (1) verify the generator relation and
(2) an integrability condition to guarantee that the r.h.s. in (\ref{3.9.1})
is always finite.
\sm

(1) {\em Generator relation.}
Denote the generator of the process by $L$ and of the dual process
$(\eta, \CF^+_t)_{t \geq 0}$ by $L^\prime$. Then we have to verify
for the duality function $H$ that
\be{df1}
(LH(\cdot,(f,\eta))) (X) = (L^\prime H(X, \cdot))(\eta,f).
\ee
In order to verify the generator relation (\ref{dg30.6}),
we have to calculate $LF$, recall (\ref{1.12}) and for that we evaluate
first the first and second order differential operators acting on $F$.

Fix $n \in \N$ and a map $\xi:\{1, \cdots,, n\} \la \Omega_N$.
We need to calculate for functions $G$ with
\be{fu1}
F(X) =\intl_{\I^n} f(u_{1}, \cdots, u_{n})
x_{\xi(1)} (du_1) \cdots, x_{\xi(n)} (du_n), \quad X \in (\CP(\I))^{\Omega_N},
\quad X= \otimes^{\Omega_N}_i x_i
\ee
the first and second order derivatives (recall (\ref{1.39})-
(\ref{1.40})). Note that the  $\{\xi(j):j=1,\dots,n\}$ are not necessarily distinct and that repetitions  can occur.

For $\xi\in\Omega_N$ (cf. (\ref{1.39}))
\be{fu2}
\frac{\partial F(X)}{\partial \mu_{\xi}} [v]
= \sum_{\ell\in\{1,\dots,n\}\atop\xi(\ell)= \xi}\left(\;\intl_{\I^{n-1}} f(u_1, \cdots, u_{\ell-1}, v, u_{\ell+1}, \cdots, u_n) \bigotimes\limits^n_{j=1 \atop j \neq \ell} x_{\xi(j)}\right),
\ee
	

\bea{fu3}
&& \frac{\partial^2 F(X)}{\partial \mu_{\xi} \partial \mu_{\xi}} [v,{v^\prime}]\\
&& =\sum_{\ell,\ell^\prime\in\{1,\dots,n\},\ell\ne\ell^\prime\atop\xi(\ell)=\xi(\ell^\prime)= \xi}\left(\intl_{\I^{n-2}} f(u_1, \cdots, u_{\ell-1}, v, u_{\ell+1}, \cdots, u_{\ell^\prime-1},
 v^\prime, u_{\ell^\prime+1}, \cdots, u_n) \bigotimes\limits^n_{\ell=1 \atop \ell \neq i,j}
 x_{\xi(\ell)} (du_\ell)\right).\nonumber
\eea
Note that if $|\{\ell:\xi(\ell)=\xi\}|=k$, then there are $\left(\begin{array}{c}k \\ 2\end{array}\right)$
summands.

We can now apply this formula to our mixed moment $F$, recall (\ref{3.2}).
Recall (\ref{1.12})

\bea{fu4}
(LF)(X)&& =(L^{\rm{mig}}+L^{\rm{sel}}+L^{\rm{mut}}+L^{\rm{sam}})F(X)\\&& =\sum_{\xi \in \Omega _N} \Bigg[ c
\suml_{\xi^\prime \in \Omega_N}a_N(\xi ,\xi ^{\prime })\intl_{\I}
\frac{\partial F(x)}{\partial x_\xi }(u)
(x_{\xi ^{\prime }}-x_\xi )(du)\nonumber\\
&&+\; s \intl_{\I}\left\{  \frac{\partial F(x)}
    {\partial x_\xi }
    (u)\left(\chi(u)-\int_{\I} \chi(w)x_\xi(dw) \right)\right\}x_\xi(du)
\nonumber\\
&&+\; m\intl_{\I}  \left \{ \intl_{\I}
    \frac{\partial F(x)}{\partial x_\xi }
    (v)M(u,dv)
 -  \frac{\partial F(x)}{\partial x_\xi }(u)\right \}
    x_\xi (du)
    \nonumber \\
&&+ \;d\intl_{\I} \intl_{\I}
    \frac{\partial ^2F(x)}{\partial x_\xi \partial x_\xi }
    (u,v)Q_{x_\xi }(du,dv)  \Bigg ], \qquad x \in (\CP(\I))^{\Omega_N}, \nonumber
\eea where $Q_x(du,dv)$ is given by (\ref{1.13}).

We now apply this to the function $F(\cdot)=H(\cdot,(\eta,\CF))$ (cf. (\ref{dg30.3})).

We first consider the action of $L^{\rm{mut}}$ using (\ref{fu1}),
\bea{add180}
 &&L^{\rm{mut}} H(X,(\eta,\CF))\\
 && = \hspace{1cm}
       m\intl_{\I^n} \Big(\sum_{\ell\in\{1,\dots,n\}}\nonumber\\
&& \qquad
\Big(\;\big[\int_{\I}f(u_1, \cdots, u_{\ell-1}, v, u_{\ell+1}, \cdots, u_n)M(u,dv)-f(u_1, \cdots, u_{\ell-1}, u, u_{\ell+1}, \cdots, u_n)\big]\Big)\nonumber\\
&& \hspace{11cm} \bigotimes\limits^n_{j=1} x_{\xi(j)}(du_j)\Big)\nonumber\\
 &&= \hspace{1cm}
     H(X,(\eta,\mathcal{M}^\ast f)).\nonumber
\eea

We next consider the action of $L^{\rm{mig}}$ using (\ref{fu1}),
\bea{fu60}
&&L^{\rm{mig}} H(X,(\eta,\CF))\\
&& \hspace{1cm} = \quad c\intl_{\I^n} \Big(\sum_{\ell\in\{1,\dots,n\}}
  \Big(\; \int_{\I}f(u_1, \cdots, u_n)\nonumber\\
&& \hspace{2cm}  \Big[\sum_{\ell=1}^n\sum_{\xi^\prime}{a_N({\xi(\ell),\xi^\prime})}\bigotimes\limits^n_{j=1\atop j\ne \ell} x_{\xi(j)}(du_j)\otimes x_{\xi^\prime}(du_\ell)-\bigotimes\limits^n_{j=1} x_{\xi(j)}(du_j)\Big]\Big)\nonumber\\
&& \hspace{1cm}
= \quad H(X,\sum_{\ell=1}^n\sum_{\eta^\prime} \wh a_\ell(\eta,\eta^\prime)(\eta^\prime, f))=L^{\prime,\rm{mig}} H(X,(\eta,\CF))\nonumber
\eea
where $\eta_\ell^\prime =(\zeta,\pi,\xi_\ell^\prime)$, $\xi_\ell^\prime(j)=\xi(j),\; j\ne \ell,\;\xi_\ell^\prime(\ell)=\xi^\prime$, $\wh a_\ell(\eta,\eta_\ell^\prime)=a_N(\xi(\ell),\xi^\prime)$.

\bigskip

Consider the following function on the state space of the dual process.
Fix an element $X \in (\CP(\I))^{\Omega_N}$ and define
for $n \in \N$,  partition $\pi$ of $\{1,\cdots,n\}$,  map
$\xi : \{1,\cdots,|\pi|\} \to \Omega_N$ and  $f:\I^{|\pi|} \to \R^+$,
define

\bea{fu5}
&&G((\eta,f)) = \intl_{\I^{|\pi|}} f(u_{1}, \cdots, u_{|\pi|})
x_{\xi(1)} (du_1) \cdots x_{\xi(|\pi|)} (du_{|\pi|})\\&&
=\intl_{\I^{|\pi|}} f(u_{\xi(1)}, \cdots, u_{\xi(|\pi|)})
\bigotimes x_{\xi} (du_{\xi}) =H(X,(\eta,f)).\nonumber
\eea
We have to calculate the action of the generator $L^\prime$, i.e.
we have to determine $(L^\prime G)((\eta,f))$, where $L^\prime$
is the generator of the pure Markov jump process with a piecewise deterministic
part (for the mutation) which we defined as dual process. We obtain by explicit calculation
\bea{fu61}
&&L^{\prime,\rm{mig}} H(X,(\eta,f))=(L^\prime G)((f,\eta)) = \\
&&\hspace{2.7cm}
\intl_{\I^{|\pi|}} \big[\sum_{\ell=1}^n\sum_{\xi^\prime} a_N(\xi(\ell),\xi^\prime)(f(u_{\xi(1)},\cdots,u_{\xi(\ell-1)},u_{\xi^\prime},u_{\xi(\ell+1)},\dots u_{\xi(|\pi|)})
\nonumber\\
&& \hspace{7.6cm} -f(u_{\xi(1)}, \cdots, u_{\xi(|\pi|)}))\big]
\bigotimes x_{\xi} (du_{\xi}) \nonumber\\
&& =L^{\rm{mig}}H(X,(\eta,f)). \nonumber
\eea


Combining both (\ref{fu60}) and (\ref{fu61})  we obtain
the claim  (\ref{df1}).  The generator calculations corresponding to selection and resampling (coalescence) follow
in the same way (compare \cite{DG99}) and this completes the required  generator calculation.

{\em (2) Finite expectation}

Here we have to show that
\be{add126}
E_{(\eta_0, \CF_0)} [ H(X, (\eta_t, \CF_t))] < \infty \mbox{ for all } t \geq 0,
\ee
in order to use Theorem 4.11 in
\cite{EK2} to conclude the duality relation from the
generator relation we established in the previous point.
We begin by verifying this for some $t_0 >0$ which is
sufficiently small.

     Namely we realize first that all transitions occuring in the
     dual preserve the $\|\cdot \|_\infty$-norm of $\CF_t$ except
     the selection transition. Here we have $g \to \chi g +(1-\chi)
     \otimes g$ which satisfies
\be{transi1}
     \|\chi g +(1-\chi) \otimes g \|_\infty \leq 2 \|g\|_\infty.
\ee
     Since the number of selection events is given by a pure birth
     process with birth rate $s$ we have
\be{transi2}
     \|\CF_t \|_\infty \leq 2^{N_t} \cdot \| \CF_0\|_\infty.
\ee
     Hence by explicit calculations of the Laplace-transform of $N_t$
     we have:
\be{transi3}
     E[\|\CF_t\|_\infty] \leq (1-p) (\frac{2}{1-2p})^{N_0}\|\CF_0\|
     < \infty, \mbox{ if } p=(1-e^{-st}) < \frac{1}{2}.
\ee
     Therefore
\be{transi4}
     E[\|\CF_t\|_\infty] < \infty \quad , \quad \forall \; t<\frac{\log 2}{s}
     \mbox{ for all } N_0 \in \N.
\ee
We now have to extend (\ref{dg30.5})  to all $t \geq 0$. We argue as follows.

First using Theorem 4.11 in \cite{EK2} with $\Gamma_T = 2^{N_T}$
and  $T< \frac{\log 2}{s}$, so that $E[2^{N_T}]<\infty$, we get
\be{transi5}
E_{X_0} [H(X_t, (\eta_0, \CF_0))]
= E_{\eta_0, \CF_0} [H(X_0, (\eta_t, \CF_t))], \mbox{ if } t\leq T,
\ee
and therefore the martingale problem for $X$ is wellposed and has as solution a Markov process
(in fact a Feller process). Furthermore we know that
we can write (this is implied by the form of the selection transition
which at each birth creates as new state a sum of two functions derived from the old function),
\be{transi6}
\CF_t = F(t, (\eta_0, \CF_0)) = \suml^{N_t}_{i=1} \CF_{t,i}.
\ee
Then we know from the Markov property of the dual process and its form that
\be{transi7}
\CF_{2t} = F(t, (\eta_t, \CF_t)) = \suml^{N_t}_{i=1} F(t, (\eta_t, \CF_{t,i})).
\ee
Then observe that the Markov property of $(X_t)_{t \geq 0}$
allows to calculate as follows.
Using the duality for time $t$ to time $2t$ (here $\wt E$ denotes expectation over the
dual dynamic)

\be{expe1}\begin{array}{l}
E[H (X_{2t}, (\eta_0, \CF_0)) |X_0] = E[E[H(X_{2t}, (\eta_0, \CF_0)) | X_t]|X_0]\\ [2ex]
\hspace{4cm} = E[\wt E[H(X_t, (\eta_t, \CF_t)) | (\eta_0, \CF_0)]X_0]\\[2ex]
\hspace{4cm}= E[\wt E[H(X_t, (\eta_t, \suml^{N_t}_{i=1} \CF_{t,i})) | (\eta_0, \CF_0) ] |X_0].
\end{array}
\ee
Now we calculate the r.h.s. of the above equation as:
\be{expe2}\begin{array}{l}
\wt E[ E[H(X_t,  |  (\eta_t, \suml^{N_t}_{i=1} \CF_{t,i})]|X_0] (\eta_0, \CF_0))
= \wt E \left[\suml^{N_t}_{i=1} [H(X_0, (\eta_{2t}, F(t, \CF_{t,i}))) |(\eta_0, \CF_0)]\right]\\[2ex]
\hspace{5cm} = [\wt E \left[ H(X_0, (\eta_{2t}, \suml^{N_t}_{i=1} F(t, (\eta_t,  \CF_{t,i})))|(\eta_0, \CF_0)  \right]\\[2ex]
\hspace{5.4cm} = \wt E[\wt E[ H(X_0, (\eta_{2t}, \CF_{2t})) | (\eta_t, \CF_t)]|(\eta_0, \CF_0)]\\[2ex]
\hspace{7cm} = \wt E[ H(X_0, (\eta_{2t}, \CF_{2t}))|(\eta_0, \CF_0)],
\end{array}
\ee
where we used duality between  times $t$ and $2t$, the construction of the dual and its Markov property. Hence
we get the duality relation up to time $2T$. Iteration gives the claim  for all positive times.

\end{proof}

\subsection{Two alternative duals for the mutation component}\label{modified}

In order to handle certain applications  where mutation is a key it is sometimes useful to
modify the dual representation. Here we give two modifications of
the mutation induced part of the dual as well as their combined version
which will be crucial later on also for the refined version of the
dual which is the main tool for the renormalization analysis.

\subsubsection{Modified dual for  state-independent mutation component}
\label{mutcomp1}
In the case in which there is a non-zero state independent component
 of the mutation mechanism, that is, $\bar m >0$
and $mM \geq \bar m 1 \otimes \rho, \quad \rho$ strictly positive on $\I$,
 we can obtain a modified dual which is particularly useful when $
 mM-\bar m (1 {\D \otimes} \rho) =0$.

 \beD{D.moddu1} {(Modified mutation dual)}

 We define the {\it modified dual process} $(\wh \eta_t, \CF_t)_{t
 \geq 0}$ resp. $(\wh \eta_t,\CF^+_t)_{t \geq 0} (\wh \eta_t, \CG^+_t)_{t \geq 0}$ for a mutation matrix
 $M$ satisfying $m M \geq \bar m 1 \otimes \rho$ as follows.
 First, we enlarge the  space $\Omega_N$ to
 \be{Z1}
 \Omega_N \cup \{\ast\}. \ee
 In $\{\ast\}$ all the transition rates of the dual particle process $\eta$
 are 0. In addition we set for the function-valued part $M(\{\ast\},\{\ast\})=1$ and
 for consistency extend $\rho$ by $\rho (\{\ast\})=0$.

  The dynamics of the process $(\wh \eta_t, \CF_t)$ resp.
  $(\wh \eta_t,\CF^+_t)_{t \geq 0}, (\wh \eta_t, \CG^+_t)_{t \geq 0}$ is now obtained by adding
  to the mechanism of $\eta$ for partition elements still located on $\Omega_N$ (rather
  than $\ast$)
  \be{Z2} \mbox{ jumps of the partition elements to } \{\ast\}
  \mbox{ after exponential waiting times at rate } \bar m
  \ee
  independently of each other (and independent of those of all other
  transitions)
  and by changing the dynamics of $\CF_t, \CF^+_t, \CG^+_t$ by
  \bea{an5}
   &\mbox{replacing } m M^\ast_t \mbox{ generated by the mutation kernel $M$ by the semigroup }\\
    &M^\ast_t \mbox{ corresponding to the transition rates }
     (m M-\bar m 1 \otimes \rho). \nonumber \qquad \square
      \eea
      \end{definition}

      \begin{remark} \quad Note that $\wh \eta_t = (\wh \xi_t, \wh \pi_t)$ and
      both components have a different law than $(\xi_t, \pi_t)$, since
      in particular on $\ast$ there is no coalescence taking place.
      This observation has  important consequences. In the case $ M-\bar
      m (1 \otimes \rho) =0$ once a partition element reaches $\{\ast\}$
      this element does not undergo any further change and the dual process
      can be easily analysed since it is eventually trapped with
      all locations on $\{\ast\}$.
      \end{remark} \bi

With this modification the same duality relation between
$(X_t)_{t \geq 0}$ and $(\wh \eta_t, \CF_t)_{t \geq 0}$ and its variants holds
if we enrich by an additional component associated with the site
$\ast$, in exactly the form given in Proposition \ref{P3.1} respectively Theorem
\ref{T.orddual}. Precisely:

\beP{P.dual} \quad (Modified Duality for state-independent mutation jumps)

Let $(X_t)_{t \geq}$ be as in Proposition \ref{P3.1} satisfying $\bar m>0$ and
extend it to a process on $\Omega_N \cup \{\ast\}$.
Here the state of the original process $X(t)$ is defined in the additional
state $\ast$ for all times as
the probability measure $\rho$ on the type space which is the measure giving the
state-independent part of the mutation rates.

Let now
$(\wh \eta_t, \CF_t)$, respectively $ (\wh \eta_t,\CF^+_t), (\wh \eta_t, \CG^+_t)_{t \geq 0}$ denote
the modified dual processes defined in (\ref{Z1})-(\ref{an5}). Then
the analogues of (\ref{3.8}),(\ref{3.9})  and (\ref{3.9.1}) hold.
$\qquad \square$
\end{proposition}

\begin{proof}
First we note that we can decompose the mutation semigroup in the independent superposition of
the state-independent and the stat-dependent part.

Next note that the expectation of the newly introduced jump occurring
in our test function is exactly given by the action of the state-independent part
of the mutation semigroup. In particular we have not changed the expected value
of the test function switching to the new dynamic. Since the duality relation
involves on both sides the expectations we get the claim. Alternatively apply
(\ref{dg30.6}) and note that the state-independent part of the generator and the
part for the jump to $\{\ast\}$ in the dual satisfy this relation.
\end{proof}

\subsubsection{A random representation of the (state-dependent) mutation term of the dual}
\label{mutterm}

The modification we now describe  is useful in dealing  with mutation
which is not state-independent, for example in our context when we have also rare
mutations from one level to the next. In that case it is useful to change the dynamics of
$\CF_t; \CF^+_t, \CG^+_t$ by replacing the \emph{deterministic} function-valued evolution driven
by the mutation semigroup as specified in Definition \ref{D.Funkvdu2}
or in Definition \ref{D.moddu1} by a \emph{random} and {\em function-valued}
jump process. The process $\eta$ remains untouched.

Namely remove the deterministic evolution from (\ref{an4ext}) and
replace it by the following Markov pure jump process in function
space. Introduce the following jumps in function space
for $g \in L_\infty (\I^k)$:
\be{dy1}
g(u_1,u_2,\cdots,u_k) \longrightarrow ( M_ig)(u_1,u_2,\cdots,u_k),
\quad k \in \N,
\ee
where $M_i$ denotes the application of the
operator $M$ to the ith variable of the function.
These jumps occur for a function $f \in L_\infty (\I^k)$ at rate $m \cdot k$.
This defines (uniquely) a Markov pure jump
process with values in $(\bigcup_n L_\infty (\I^n)$ starting from
any initial state.

Consider now the situation that $\CF_t, \CF^+_t, \CG^+_t$
depends on $k$ variables. For each $i \in \{1,\cdots,k\}$
the jumps in (\ref{dy1})  occur at rate $m$ and this jump time is independent of
everything else. Therefore as long as the state of $\CF_t$, respectively
$\CF_t^+$ is a function of $k$ variables these random transition
from (\ref{dy1}) occur after  exponential waiting times. We denote by
\be{dy2}
 \wh\CF_t,\wh\CF^+_t, \wh \CG^+_t \ee
the resulting modified function-valued processes. Note that in particular
we can assign then the mutation jump occuring always with a particular
partition element, which will be of some relevance for the historical
interpretation.

Since the expectation of the jumps occuring in our test function is by construction
given by the action of the mutation semigroup and since the duality relation only
claims the identity of two expectations, we have not changed the r.h.s. of the
duality relation and conclude the following.

\beP{P.dumoM} All the previous duality relations remain valid if we replace
$(\CF_t)_{t \geq 0}$ by $(\wh \CF_t)_{t \geq 0}$ or $(\CF^+_t)_{t
\geq 0}$ by $(\wh \CF^+_t)_{t \geq 0}$, similarly $(\CG^+_t)_{t \geq 0}$
by $(\wh \CG^+_t)_{t \geq 0}$ and leave $\eta$ untouched.
The same holds for the process
$(\wh \eta_t, \CF_t)_{t \geq 0}, (\wh \eta_t, \CF_t^+), (\wh \eta_t, \CG^+_t)_{t \geq 0}$
from Definition \ref{D.moddu1}.
   $\qquad \square $
   \end{proposition}

   \begin{remark} Note how here the interplay between selection and mutation is reflected
   in this form of the duality.
   For example assume that  we have rare mutation events, i.e. with rates $<<s$, the
   selection rate. Then: the rate of mutation events in the dual process is proportional to the number
   of
   partitions elements we have in the dual process and therefore rare
   mutations become visible in the dual process as soon as  many births have
   occurred in the dual process due to a much higher selection rate,
   which then compensates via the large number of individuals a small
   mutation rate.

\end{remark}

\subsubsection{Pure jump process dual}
\label{sss.purejump}

We can use both constructions presented in the previous two subsubsections
at once:
\beD{D.pure} {(Pure jump process dual)}

In particular in combining the new representations of the state-dependent
part $m M-\bar m 1 \otimes \rho$ and the state-independent part in $\bar m 1 \otimes \rho$
we get a pure Markov jump process
\be{dg33}
(\wh \eta_t,\wh \CF_t)_{t \geq 0}, \quad (\wh \eta_t,\wh \CF^+_t)_{t \geq 0},
(\wh \eta_t, \wh \CG^+_t)_{t \geq 0}.
\qquad \square
\ee
\end{definition}

It is the dual process on which we base our refined dual in Subsection
\ref{ss.refine} and it is also the version best suited for the historical interpretation
since it generates a marked (locations, mutation events) random graph in
which we can find the marked ancestral tree of a tagged subpopulation.
(Compare Subsubsection \ref{sss.histint}).

\subsection{Historical interpretation of the dual process}
\label{histint}

Does this analytical construction above have a heuristic meaning
or is it simply an analytical trick revealing an explicitly
solvable model? Both these alternative situations occur in the theory of
stochastic processes. In our model, however, the duality can be
interpreted in a nice way if one considers what is sometimes
called the historical process.
In fact this duality can be extended to reach processes including more
information about the individuals past (see \cite{DGP}).

In order to understand this, one has
to remember that the diffusions we work with arise as small mass -
many individual limits of  particle models. In such a particle
model it is possible to follow through time, the fate of
individuals and their descendants and define for each individual currently
alive its ancestral path leading back to its father then its grandfather and so on
where this ancestral path gives the location and the type of the ancestor
at every time.
In fact we can from a time horizon $t$
backward generate the collections of ancestral path and their genealogical
relation which define a tree or rather a forest of trees corresponding to
different founding fathers. This way we obtain a random marked (with types and locations) forest.
A more modest question would be to determine the types and the
time back to the first common ancestor of the sample taken from time $t$
(this means we only get the types at time $t$ and not at all earlier
times along the ancestral path.)

How can we study this complicated object? For each given time
horizon $t$ we can zoom in on a \emph{finite subpopulation of
tagged individuals}, tagged in the time $t$ population and trace
 their history \emph{backwards in time}.
 Then we can ask whether this object can be generated by a suitable
 stochastic process, the {\em backward process}.
 Can we hope for this process to again be Markov and time-homogeneous?

  To understand this better, take the
  non-spatial case first. Ask the question: What are the types of a
  $k$-sample of individuals from the time $t$ population.
  To answer this question generate the
  law of the historical evolution of this subpopulation by running a
  stochastic process backward. Similarly in a spatial model one can
  take samples at different locations. A nice case arises if this backward
  dynamic of the sample alone is already a Markov process, which is time-homogeneous.
It is in this situation that one traditionally speaks of the existence of
  a dual process.

  Resort first to a simpler, the {\em neutral} case. Due to exchangeability of individuals
  for the neutral case (i.e. no selection and mutation) such a Markovian
  backward process can be given based on a coalescent generating the family structure
  of the sample. A key tool to establish this in our context is the
  representation in terms of the lookdown process of Donnelly and
  Kurtz adapted to this spatial situation, (compare \cite{GLW}).
  This then allows to rigorously establish and identify the dynamics
  of the backward process in terms of a spatial coalescent (\cite{GLW}).
This can be extended and one can show that the genealogical trees of the
neutral model evolve in such a way that their state at a fixed time $t$
is given by the genealogical tree associated with the coalescent, see
\cite{GPWmetric} and \cite{GPWmp}.

In the case in which selection and mutation both depend on the type
of the individual  exchangeability is no longer preserved and
a complicated interaction between the tagged subpopulation in the sample and the
remaining population arises, which is not anymore in law equivalent
to an autonomous evolution of the sample.
Hence in order to generate the genealogy by
a backward process we have to add a richer device which
consists of a reservoir of possible histories in order to still obtain  a Markov
process driving the backward picture.

In the literature this problem has been treated first in models
where the process and its dual can be specified via a random graph
(graphical representation). In the case of the existence of
graphical representations, such as the voter model or stochastic
Lotka-Volterra models which are of that type and which appear in
the models and work of Krone and Neuhauser (\cite{KN97}) there are
typically arrows between points in the random graph which might or
might not be used which describe the possible action of selective
forces.

In our context selection
requires, as explained above, adding new individuals to our dual (sampled) population
 as we move backward from the time horizon which represent the
 \emph{potential} insertion of a fitter type in the tagged population.
 Whether such a potential insertion actually takes place depends on
 the fitness of the ''victim'' in the tagged population and of the
 fitness of the potential intruder. This means that in
 order to decide the types of the $k$-sample of the population we
 have to consider for selection a growing dual population
representing typical individuals drawn from the population at
a certain time $s$ back. Then we have
 to assign weights to these possibilities representing the
 probabilities with which they are realised.
Therefore the ancestral lines of our sample form a subtree
of the graph generated by the ordered particle process $\eta$.
Among the possible choices for subtrees in this graph the
various possibilities have probabilities which we can read
of from the function-valued part.

 In addition to the complication arising from selection,
 along the ancestral path mutation events have to
 be taken into account. For the state-independent part such an event
 decouples the final type from everything happening earlier, but
 if the mutation is state-dependent, we have again to consider
 the complete system of backward path but since the tagged
 population is small compared to the basic population a law of
 large number effect occurs and this can be represented through
 a functional dual reflecting the fact that mutations occur based
 on the current type and do not depend on the overall populations.
However with the mechanism for the dual as described in Subsubsection
\ref{mutterm} we can even associate with every individual a chain of
mutation events. But we have not yet a rich enough system to
associated with the ancestral path of the individual, a path
in type space. We will see in Subsection \ref{ss.refine} how this possibility
arises.

 In the spatial context we have to sample $k_1$ particles from a
 site $\xi(1)$, $k_2$ from site $\xi(2), \cdots, k_m$ from site
 $\xi(m)$ and the ancestral paths migrate in space. Altogether we
 therefore get a spatial coalescent with birth and with mutation
 operations associated with each ancestral path which represents
 the historical evolution of a randomly drawn $k$-tuple from the
 time $t$-population of the original process. Then we  can use the dual
 to calculate the probabilities that $k$-sampled individuals have specific
 types, a specific genealogy and paths in geographic space. We shall discuss more
 of this as we go along, see also the explanation of the refined dual
 from a historical process perspective in Subsubsection \ref{sss.histint}.

\subsection{Refinements of the dual for finite type space: $(\eta_t, \CF^{++}_t)_{t \geq 0},
(\eta_t, \CG^{++}_t)$}
\label{ss.refine}

In this section we focus on the case of a {\em finite} type space and we consider refinements of the dual $(\eta_t,
\CF^+_t)$, (resp. $(\wt \eta_t, \CF^+_t)$ or $(\eta_t, \wh \CF^+_t)$) and similarly with the
versions using $\CG^+, \wh \CG^+, (\eta_t$, which we denote $\CG^{++}_t)_{t \geq 0}$ which  mainly require an enrichment of
the mechanism corresponding to selection and mutation in the function-valued
part of the dual and which works (only) for a smaller set of functions
in which the function-valued part can start,
namely the function $f$ in (\ref{3.2}) has to be a {\em sum of products of
indicators}. This set of functions however is still
generating the full set of functions we consider and is therfore in particular
distribution-determining for finite type space and hence suffices for a duality theory.
However to guarantee that this subset of functions is {\em preserved} under the
part of the dual dynamics corresponding to selection and mutation dynamics, we
have to change the dynamics of $(\eta, \CF^+)$ or $(\eta, \CG^+)$.

In particular this duality allows for a nice historical interpretation, since
it generates a marked graph (marked with types and locations) of which the
ancestral marked tree of a finite sample is a subtree (see Subsubsection
\ref{sss.histint}).
      \bi

This constructions we describe in this subsection are written for the case of a \emph{finite} type
space
\be{fts1}
      \mathbb{I}=\{1,\dots,K\}.\ee

\begin{remark} \quad
An extension to countable type space is possible, but not needed in
the sequel. In particular
if we assume more about the fitness functions, for example that fitness
values have 1 as the only accumulation point, then we can handle also
immediately the case of countably many types.
\end{remark}

\subsubsection{Ingredients} \label{sss.ingredients}

For this refined duality we use a suitably chosen smaller set of
functions as state space for the function-valued process.
However the key point is to also use a modified  dynamic of the function-valued
process $\CF^+_t, \CG^+_t$ (later on called $\CF^{++}_t, \CG^{++}_t$) which is generated by
a {\em refinement} of the {\em birth process} in
$\eta_t$, the particle process driving the function-valued process
together with a {\em refinement} of the {\em mutation part} of the function-valued part of the dual process.
The point of these modifications (leading in a way to a more complicated dynamic), will be (1) the
smaller subset of test functions is preserved and (2) this is suitable for development in \cite{DGsel} to  consider the time-space process
of this dual and use the additional information built into the dynamic to define a process on
a richer state space (space of tableaus representing decompositions of the set $\I^\N$) which
can then be analysed better if we are concerned with the longtime behaviour.

In points (0) -(iv) we now explain step by step the needed changes in the state space and the four
mechanisms in the dual to obtain the refined version.

\emph{(0) State space of function-valued part}

The state space for the dual process better its function-valued
part is the subset of functions on $ \mathbb{I}^{\mathbb{N}}$ arising as follows.  For
$k\in \mathbb{N}$, define first the space of certain finite sums of products of at
most $k$-indicator functions of a variable in $I$ each:

\be{set1}
\mathbb{F}_k:=
 \Big\{f=\sum_{i=1}^n\prod_{j=1}^{k_i} 1_{B_{i,j}}(u_j),B_{i,j}
 \subseteq \{1,2,\cdots,K\}, \quad n \in \N,  k = \max_{i=1,\dots,n}k_i\}.
\nonumber \ee
Here $f$ above is viewed as a function of $\I^k$.

\begin{remark}\quad
Note that the states in this set of functions need not satisfy
$\int f d \mu^{\otimes k} \leq 1$, as would be the case if $f$ defines
a decomposition of $(\I)^k$. Therefore we have two cases, the version of the dual
based on $\CF^+$ where this is not the case and the dual based on $\CG^+$, where
in the definition above we can impose the condition
\be{condi1}
\intl_{\I^k} f d\mu^{\otimes k} \leq 1
\ee
and still obtain a set of states preserved under the dynamics.
\end{remark}

   The \emph{state space} of the function-valued component of the refined dual process
   is given by: \be{Gre2}
       \mathbb{F}:= \bigcup^\infty_{k=1} \mathbb{F}_k. \ee

 This function space is
   associated with the situation in which  we have $k$ particles in
   the dual particle process $\eta$. For consistency the product running from $k_i$ to $k$ could be
   filled up with $k-k_i$ indicators of the whole type
   space $\mathbb{I}$. The parameter $n$ allows us to consider
   the evolution of functions in which there is a mechanism that replaces
   a  product of indicators by a sum of products of indicator functions.

       With each variable we associate a position in space. However this is
       not changed under selection, mutation or resampling. Hence in order
       to explain the dual mechanism for each of those we can ignore the
       spatial aspect for the moment.
       \bi

       We now explain the dual mechanisms corresponding to selection,
       mutation, resampling and migration in the original process step by step.

       \noi \emph{(i) Selection.}

In order to describe the change of the transition  in $\eta$ at a birth event and
the new transition in the function-valued
component of the dual upon a birth event we use the following structure.
Assume that there are $\ell$ fitness levels reaching from 0 to 1
and denote them by
\be{Gre3}
0=e_1<e_2<\dots<e_{\ell}=1. \ee
(Note without loss of generality to assume for bounded fitness functions in assuming that $e_1=0,\;e_{\ell}=1$.)

First change the dynamic of the first component, i.e. $ \eta$,  so that the
  birth process is replaced by a \emph{multitype} birth process,
  more precisely a $(\ell-1)$-type birth process.
  Births occur at total rate $s$ but now when a birth occurs the
{\em  type of birth} ($\in \{2,\dots,\ell\}$)  is chosen i.i.d. for this event with probabilities
\be{Gre13c}
(e_i -e_{i-1})_{i=2,\cdots, \ell}. \ee

We shall now introduce the jumps in the function-valued part occuring upon
a birth in $\eta$ of $\wh{}$ specific type.
For each tagged individual of a certain level of fitness say $e_j$, all those
individuals with a strictly larger level of fitness can be a
candidate to take the place of the tagged individual.
Define therefore for each level of fitness $i$, with
$i=2, \cdots,\ell$ the set
\be{Gre4}
A_i:=\{j\in\mathbb{I}:\chi(j)\geq e_i\}.
\ee
We note that $A_{i+1}\subset A_i$.

 We define for a subset $C \subseteq \{1,2,\cdots,K\}$ of the type
 space the operator $\psi_C$ on $\F$  as follows. Let $f$ be a
 function of $k$-variables and define
 \be{Gre13} \Psi_{C}
 f:=\sum_{m=1}^k  [1_C(u_m)f(u_1,\dots ,u_k)
 +f(u_1,\dots,u_{k})(1-1_C(u_{k+1})].
 \ee

 Observe that $(\psi_C -Id)$ is the generator of a rate 1 jump
 process on $\F$. Namely denote $1^\ell_C = 1_C(u_\ell)$ and
 introduce for $f \in \F_k$ for every variable $u_j, j=1, \cdots,k$
 at rate 1 transitions:
 \be{Gre10}
 f \longrightarrow 1^j_C f +
    f \otimes (1 - 1_C^{j}), \quad \quad j=1,\cdots,k. \ee
 Then each jump from some $f$ in $\F_k$ ends in
   $\F_{k+1}$.

We apply this to our context for $C$ being replaced by the sets
defined in (\ref{Gre4}), we use the notation
\be{Gre10b}
\chi^j_i = 1_{A_i} (u_j) \ee
and set for $f \in \F_k$,

\be{Gre11} (\psi^{j,k}_i f) (u_1,
\cdots, u_{k+1}) = ((\chi^j_i f) \otimes 1^{k+1}_\I) +
   (1-\chi_i^{j})\otimes f))(u_1, \cdots, u_{k+1}). \ee
Again we can define an alternative transition as follows
\be{Gre11b}\begin{array}{l}
(\wh \psi^{j,k}_i f) (u_1, \cdots, u_{k+1})
= \chi^j_i (u_i) f(u_1,\cdots, u_n) 1(u_{k+1} \\[1ex]
\hspace{6cm} +(1-\chi^j_i (u_i)) f (u_1, \cdots, u_{i-1},u_{k+i}, u_{i+1}, \cdots, u_k).
\end{array}
\ee

   \begin{remark} \quad It is easy to verify that if $0\leq \int f
   d\mu^{\otimes k}\leq 1$, then $0\leq \int \psi^{j,k}_if
   d\mu^{\otimes (k+1)}\leq 2 \|f\|_\infty$ for any probability measure $\mu$,
   respectively $0 \leq \int \wh \psi^{j,k}_i f d \mu^{\otimes k+1} \leq \|f\|_\infty$.
   \end{remark}

\begin{remark} Note that this definition means that if $f$ is
a sum of products that the jump occurs simultaneously in the
corresponding factor in all summands.\end{remark}

 \begin{remark}
 A different ordering of the factors and variables in different summands and a modified
 dynamics will be used later to  couple
 the dynamics of different summands arising with every birth event.
 \end{remark}

\beD{D.seljump}{(Selection jump)}

Introduce transitions (jumps in $\F$) for the function-valued part of the form:
\be{Gre13c2} f
\longrightarrow \psi^{j,k}_i f \mbox{ from } \F_k \to \F_{k+1}, \ee
whenever a {\em birth of type $i$}
due to partition element $j$ in $\eta$ occurred.
(Recall here $i$ is chosen with probability $(e_i-e_{i-1})$).
This  will replace the transition we had in $\CF^+_t$ before. For
$\CG^{++}_t$ we use $\wh \psi^{jmk}_i$. $\qad$
\end{definition}

Now the r.h.s. of (\ref{3.6b}) is interpreted as rate $s$ births
being of type $i$ with probability $(e_i-e_{i-1})$ and leading
to a transition given by (\ref{Gre13c2}).

The type of birth does
  not influence the further evolution of $\eta$, it only will change
  the function-valued part occurring at this transition of $\eta_t$;
  we therefore do not enlarge the state space of the process $\eta$
  to store the type assigned to the birth event. \sm

\noi \emph{(ii) Mutation.}

We have to specify the action of
mutations on functions $f \in \F_k$, i.e. functions of
$k$-variables in  the special case in which  they are certain sums of
products of indicators. We use here a refinement of the random
representation of the mutation semigroup as used in Subsection
\ref{mutterm} since the latter does not necessarily preserve
indicator functions.  We now construct a random
{\em indicator-function-valued jump process} that represents the
mutation semigroup $(M^\ast_t)_{t \geq 0}$.

\begin{remark}

(Set-valued dual for Markov chains)

The construction we give below produces in fact for every Markov jump
process on a finite state space a set-valued dual process. Let $E$ be its
state space, then we can calculate
$p[Z_t = i |Z_0=j] = E[i \in \CA_t | \CA_0=j]$
for a set-valued process, i.e. values in $2^E$ called $(\CA_t)_{t \geq 0}$
with jumps and rates in Definition \ref{D.svmutj}.
In particular we can calculate also its equilibrium distribution.
\end{remark}

 This process acts independently on each
 variable of the function $\CF^{++}_t$ (corresponding to partition elements of
 the dual process). That is, at random times the function $\CF^{++}_t$
 or $\wh \CF^{++}_t$ for our new function-valued process
 in the set of sums of products of indicators
 is changed by a jump from a product of indicators
 to a new product of indicators, where
 \be{grev29}
 \mbox{all factors for a {\em given variable } change in {\em every} summand at once}.
 \ee

 We next describe the action of this modified mutation dual acting
 on one variable in the argument of $f$  (to keep the
 notation simple we think of $f$ as a function of one variable with
 the others fixed).

The mutation semigroup driving the
function-valued part of the dual process when acting on indicator
functions can also be represented by an indicator-function-valued dual (random) process
whose jumps are specified next. Later we extend this to the sums of
indicators.

We  specify for each pair $(i,j)$ of types the jumps $f(\cdot)\to \wt f (\cdot)$
corresponding to the mutation transition
$i \to j$, by the following prescription.
Recall $M = (m_{i,j})_{i,j, \cdots, K}$.

\beD{D.svmutj}{(Set-valued mutation jumps)}

 This transition from type $i$ to type $j$ occurs at rate
 \be{rate1} m \cdot m_{i,j} \quad ; \quad
 i,j \in \{ 1,2,\cdots,K\} \ee and results in a jump
 (depending on whether $\ell \in \{i,j\}$ or not):
 \bea{rate2}
 1_{\{j\}}& \longrightarrow & 1_{\{i\} \cup \{j\}} \\
 1_{\{i\}}& \longrightarrow & 0\nonumber \\
 1_{\{\ell\}}& \longrightarrow & 1_{\{\ell\}} \quad \ell \notin \{i,j\}\nonumber
 \eea
 for $i,j \in \{1,2,\cdots,K\}$ and $j \neq i$. $\qad$
\end{definition}

Next extend this to $\F$. For this purpose
let $B \subseteq \{1,\cdots,K\}$ and consider the indicator
$1_B(u)$. Since this indicator can be written as
\be{ind1} 1_B(u)
 = \suml_{k \in B} 1_{\{k\}} (u), \ee
we continue the transition in (\ref{rate2})
as a linear map acting on the indicators $f=1_{\{\ell\}}(\cdot)$, with $\ell \in\{1,\dots,K\}$.

More generally proceed as follows.
The transition $f \to \wt f$ associated with the parameter $(i,j)$
    is obtained by applying to $f$ the matrix $\overline{M}$ (in
    $(k,\ell)$) which we define as \be{ma1} \overline M(i,j) [k,\ell]
    = \left \{ \begin{array}{r@{\quad} l}
    1 & k=\ell\neq i\quad , (k,\ell)=(i,j)\\
    0 & \mbox{otherwise}.\end{array} \right. \ee This specifies the
    transition occurring in one of the several variables of $f$.

\beD{D.refmutj}{(Refined mutation jump)}

We now obtain  a function-valued dual process $\CF^{mut}_t$  for
the mutation semigroup acting on indicators or sums of indicators
if we introduce the following collection of jumps for $f \in \F_k$
for each $\ell \in \{1,2,\cdots,k\}$:
\be{rate3}
f(u_1, \cdots,u_\ell,\dots, u_k) \longrightarrow \suml_v \overline M_\ell
(i,j)[u_\ell,v] f(u_1, \cdots, v,\cdots, u_k); \quad i,j \in
\{1,\cdots,K\},
\ee at rate \be{rate4} m \, m_{i,j}, \ee  where
$\overline M_\ell$ indicates that $\overline M$ is applied to the
$\ell$th variable. $\qad$
\end{definition}

For every variable $\ell$ these jumps
preserve the sets $\F_k$ for every $k \in \N$ since the number of
variables remains fixed.  Hence the set
$\F$ is preserved under a dynamic consisting of jumps like in (\ref{rate3}).
Therefore if $\CF^{mut}_0=f\in \F_k$ for some $k$, then
$\CF^{mut}_t$ is a $\mathbb{F}$-valued process.
\bi

\noi \emph{(iii) Resampling-Coalescence.}

Here as before we identify two variables located at the same site in the
 corresponding factor of all summands of the element in $\F_k$ we
 are dealing with. When coalescence occurs the resulting
 (coalesced) element is given the lower of the indices of the
 coalescing elements and the remaining elements are reordered to
 eliminate gaps but preserve the original order.
 Note that this operation turns two factors
$1_A (u_i), 1_B(u_j), \quad i<j$
 each related to one of the
 coalescing variables into the indicator of the intersection of the
 two sets as a function of the ''new'' merged variable.
That is upon coalescence:
\be{vco1}
1_A (u_i) 1_B(u_j) \la 1_{A\cap B}(u_i),\quad
1_C (u_\ell) \la 1_C (u_\ell), \ell \notin \{i,j\}.
\ee
Therefore again $\F$ is preserved since this operation sends
 \be{op1}
 f \in \F_k \longrightarrow \wt f \in \F_{k-1}.
 \ee

\noi \emph{(iv) Migration}

The action of the migration is as before and has only to do with the
 process $\eta_t$ assigning locations to the variables but not with
 the jumps in the function-valued part.
 \bi

 \begin{remark} We can now proceed as in the definition of $(\wh \eta_t, \wh \CF^+_t)$ and add a
 state $\{\ast\}$ to the geographic space where all rates are set
 equal to zero. Note that under this dynamic $\F$ is preserved.
 \end{remark}

 \subsubsection{The refined dual $\{\eta_t, \mathcal{F}_t^{++}\}$ or $(\eta_t, \CG_t^{++})_{t \geq 0}$}
 \label{sss.dualform}

 We construct the full refined dual dynamics of the process which is
 denoted $( \eta_t, \CF^{++}_t)$ respectively $(\eta_t, \CG_t^{++})_{t})$ by setting:
 \beD{D.Gre13c}
 {(Refined dual $(\eta_t, \CF^{++}_t), \eta_t, \CG_t^{++}$)}

 (a) The process $( \eta_t, \CF^{++}_t)_{t \geq 0}$ is defined as
 $(\eta_t, \CF^+_t)_{t \geq 0}$  except that births now have  a
 random type as specified in (\ref{Gre13d}) which changes the transition
 of $\CF^{++}_t$ at a birth event as specified below and the mutation transition is replaced by a jump
 process for the indicators in the representation of $\CF^{++}_t$. Similarly we proceed for
 $(\eta_t, \CF^{++}_t)_{t \geq 0}$. Precisely the two changes are:

(i) The transition occurring in the second component at a birth event is modified as follows:

  if $\CF^{++}_t\in \mathbb{F}_k$ and a birth of a type $i$
  occurs at time $t$ from the partition element with index $j$ we have
  the transition for the factor $f$ of the $i$-th variable (recall (\ref{Gre11})):
  \be{Gre13d}
  f \to \psi^{j,k}_if=\chi^{j}_i f \otimes 1 + f \otimes (1-\chi^{j}_i),
  (\mbox{ resp. for } \CG^{++} :f \la \chi^j_i f \otimes 1 + (1-\chi^j_i) \otimes f),
  \ee
where the factor $(1-\chi^j_i)$ and its (new) variable corresponds to the newborn individual.

(ii) The mutation transition in the function-valued part is replaced by jumps according to (\ref{rate3})-(\ref{rate4}).

  (b) Similarly we define $(\wh \eta_t, \wh \CF^{++}_t)_{t \geq 0}$ or
  $(\wh \eta_t, \wh \CG^{++}_t)$ if we have
  state-independent mutation at rate $\bar m$ as the corresponding modification
  of $(\wh \eta_t,\wh \CF^+_t)$ respectively $(\wh \eta_t, \wh \CG^{+}_t)$.  $\qquad \square$
  \end{definition}

The key feature of the refined dual $(\eta_t, \CG^{++})$ is the fact that the
terms $1_A$ and $(1-1_A)$ which are created by births, evolve by further selection
and mutation to terms being identically 0 or 1, a state which we call {\em resolution}.
At resolution time one of the two summands generated at the birth event disappears.

\subsubsection{The refined duality relation}
\label{sss.refdualrel}

Now we can state and prove a duality relation between the process
$(\eta_t, \CF^{++}_t)_{t \geq 0}$ and the interacting Fleming-Viot diffusion
with mutation-selection.  \sm

  \beT{T.refdual}{(Refined duality)}

Under the assumption of finite type space we have the following duality
relation between the process $(X_t)_{t \geq 0}$ and the two refined dual processes:
if $|\xi|=k$ and
$f\in \mathbb{F}_k$, then
\be{Gre13bb}\begin{array}{l}
F_t ((\xi,f),X_0)\\[1ex]
= E_{(\xi,f)} \Big[\intl_\mathbb{I} \cdots \intl_\mathbb{I}
\mathcal{F}^{++}_t(u_{1}, \cdots, u_{|\pi_t|})
x_{\xi_t(1)} (du_1) \cdots x_{\xi_t (|\pi_t|)} (du_{|\pi_t|})\Big]\\[1ex]
= E_{(\xi,f)} \Big[\intl_\mathbb{I} \cdots \intl_\mathbb{I}
\CG^{++}_t(u_{1}, \cdots, u_{|\pi_t|})
x_{\xi_t(1)} (du_1) \cdots x_{\xi_t (|\pi_t|)} (du_{|\pi_t|})\Big]
\end{array}\ee
and analogously for the state-independent mutation and $(\wh
\eta_t, \wh \CF^{++}_t), (\wh \eta, \wh \CG^{++}_t)_{t \geq 0}. \qad$
\end{theorem}

{\bf Proof of Theorem \ref{T.refdual}}. This follows from the
previous duality relations by observing that we changed two things (I)  we use the dual we
had previously now on a {\em restricted} class of test functions $f$,
namely those in $\F$ which is a set of functions preserved under
the dynamics. (II) we have reinterpreted this dynamic on
$\F$ in an autonomous way, see (\ref{Gre12}) and (\ref{rate5}), such that the
{\em expected} jump induced in the duality function remains the same,
which we now explain in detail for the selection and mutation transition
where the changes in the dynamic occured since the duality claims that
certain expectations are equal it is preserved under the change.
\sm

{\bf (1)} First consider the {\em selection} transition.
We can  rewrite the fitness function $\chi$ as follows. Namely

 \be{Gre9b} \chi=\sum_{i=2}^{\ell}  (e_i-e_{i-1})
 1_{A_i}.
 \ee
 At the same time (since $e_1=0, \;e_{\ell}=1$):
 \be{Gre9c}
 1-\chi = \suml^{\ell}_{i=2} (e_i- e_{i-1}) (1-1_{A_i}).
 \ee

Now observe that the new state of $\CF^+_t$ after a birth event
increasing the number of basic particles from $k$ to $k+1$ can be
rewritten as follows:
\be{Gre12}
( \chi f\otimes 1 + (1-\chi)\otimes f) =
 \sum_{i=2}^{\ell}[(e_i-e_{i-1})\sum_{j=1}^k \psi^{j,k}_i f]
\ee
and
\be{Gre12b}
(e_i -e_{i-1})_{i=2, \cdots, \ell} \in \CP(\{1,\cdots,\ell\}).
  \ee
We therefore can interpret the transition $f \longrightarrow
(\chi f\otimes 1+ (1-\chi)\otimes f)$ in the process $\CF^+_t$ now
differently as superposition independent transitions defined by
$\Psi^{j,k}$ which gives then the transition in
the refined version denoted $\CF^{++}_t$.

{\bf (2)}
 To prove this alternative representation of the {\em mutation} part of
 the dual process first note that we can write the
 generator  of the mutation process, i.e. the jump process on type space
 for a fixed individual in the form of the independent superposition
 of mutation processes where each summand corresponds to the transition
 from type $i$ to type $j$. This reads in formulas:
 \be{aGre6}
   {\mathcal{M}} f(i)=\suml^K_{j=1}  M(i,j)(f(j)-f(i))
   =\sum_{k,\ell=1}^K\sum_{j=1}^K \wt M(k,\ell) [i,j] f (j),\ee where
   \be{aGre6a} \{\wt M(i,j), i=1, \cdots, K ; \quad j=1,\cdots,K; i
   \neq j\} \ee is a collection of kernels $\wt M(i,j) [\cdot,\cdot]$
   given by

\be{aGre6b}\begin{array}{ll}
&\wt M(i,j) [i,j] = m_{i,j}\\
&\wt M(i,j) [i,i]= -m_{i,j}\\
&\wt M(i,j) [k,\ell]=0,\quad \text{for } (k,\ell) \mbox{ different from } (i,i)
\mbox{ or } (i,j).
\end{array}\ee
Note that each $\wt M(i,j)$ is the generator of a mutation semigroup
with only mutations from type $i$ to type $j$ at rate  $m_{i,j}$.
Hence the collection above represents the collection of all possible mutation
jumps at their rates.

    Note that: \be{ma2} \wt M(i,j) [\cdot,\cdot] = m_{i,j} (\overline
    M (i,j) [\cdot,\cdot] - \delta_{[i,i]}), \ee which means that this
    jump in (\ref{rate2}) represents the generator $\wt
    M(i,j)[\cdot,\cdot]$.
This completes the proof of  Theorem \ref{T.refdual}.

\begin{remark}
Suppose only the mutation transitions occur.
One calculates that  the generator of this process is
$\sum_\ell \CM_\ell$ where $\CM_\ell$ is $\CM$
acting on the $\ell$-th component.
Hence the duality representation of the mutation semigroup $M^\ast_t$
acting on indicator functions is given in terms of the
indicator-function-valued dual as follows:
\be{rate5}
M^\ast_t \left(\prod_{\ell=1}^L 1_{B_\ell}(x_\ell)\right) =E[\CF^{mut}_t
(x_1,\dots,x_L)|\CF^{mut}_0 =\prod_{\ell=1}^L
1_{B_\ell}(x_\ell)].\ee
This will serve as the mutation component
in the construction of an indicator-function-valued dual for the
mutation-selection-migration system. As before, in the case of a
sum of products of indicator functions the jumps are made
simultaneously in the corresponding factors in all summands.
\end{remark}

\subsubsection{ Historical interpretation revisited and outlook on a modified dual}
\label{sss.histint}

We can use the representation of $\CF^{++}_t$ given in the
previous points to clarify the relation between the refined
duality and the selection graph of Krone and Neuhauser
\cite{KN97} and concepts in population genetics referring to the calculations
of probabilities for certain genealogical relationships.
Of course there are differences and specifics of  our model, (1) we
have multiple occupation of sites, (2) resampling only in one
colony and in addition (3) we have taken a diffusion limit of many
small mass particles. All this generates different features.

However if we consider an $n$-sample and their ancestry we are
back at a discrete model and instead of arrows attached to sites
specifying potential insertion of fitter types in the graphical
representation we generate potential insertions attached to
individuals in the sample. However our backward process is a Markov process
despite the presence of selection or mutation. This we explain now in detail.

Given sets $C_1,\dots,C_m$ the dual allows us to determine the
probability that $m$ individuals chosen randomly from the
population have types in these sets and also the probability that
they have a given genealogy.
Already in the neutral model with mutation we have to use a function-valued
dual (or for finite type space a set-valued dual) which allows to determine
the probability that a given genealogy of the sample generated by the
coalescent can result in a certain marking with types at time $t$.

Selection complicates the picture even further.
The calculations involving the new
factors created by selection provide tests comparing individuals
with the type of other randomly chosen individuals from the
population.  The corresponding probabilities are obtained by
summing over potential histories. To be more precise consider
the following.

We can now associate with our dual particle a marked graph as follows.
Starting with the tagged sample we draw from every individual an edge,
till the first selection or coalescence event occurs. If a selection
event occurs we place a vertex from which two new edges start, one
associated with $\chi^j_i$, the other with $(1-\chi^j_i)$ if the
birth occurs with the $i$-th particle and is of type $j$. If a
coalescence occurs we join the two edges and continue with one edge.
On each edge we record the current location of the corresponding
individual. Therefore edges are now marked with functions $\chi^j$
or $1-\chi^j$ and with locations.

At each edge we place at rate $m \cdot m_{i,j}$ mutation markers
indicating a mutation event on the variable of the form $i \to j$.

With this construction we have generated a {\em random marked graph} in which we consider all paths
leading from the roots (the tagged individuals) to the leaves
at time $t$. Every such path generates a product of indicators
by multiplying the function at the root and the $\chi$-functions
attached along the way. Then act with the mutation markers on the factors
attached with an edge. Finally
sum over all these products corresponding to path connecting the
leaves to a root. This is the state of $\CF^{++}_t$. Note
that if we put $f_0\equiv 1$, then we see that the expectation of the
sum of factors arising from the birth events and undergo mutation
have expectation 1.

Can we associate with these objects, i.e.
the marked random  graph the marked ancestral tree for the tagged sample
meaning we know the joint law of the genealogical distances
between the tagged sample together with the type and geographical location
at time $t$?
This would require that based on the random graph explained above
we can decide at each branching point
which way to go. If we use $\CF^{++}$ the obstacle is that we
have $\chi(u_i)$ and $1-\chi (u_{m+1})$ which since they
belong to different variables do not define a decomposition
in the potential set of  marked ancestral path into disjoint sets.

For the purpose  of achieving this we use $\CG^{++}$ instead.
Suppose we have a birth due to particle number $i$ of type $j$. Then we replace this transition by
\be{ag40-}
f(u_i) \la \chi^j (u_i) f(u_i) 1 (u_{m+1}) + (1- \chi^j)(u_i) f(u_{m + 1}).
\ee
The key feature is that selection introduces a
{\em decision tree} with corresponding probabilities for the different genealogies which
are possible after the interaction of the sample with the rest of the population.
The selection transition
$f \la \chi f +(1-\chi) \otimes f$
we defined for $\CF^{++}$ (and for $\CF^+$ of course as well)
is changed in the new picture corresponding to $\CG^{++}$ into
(suppose here $f=1_B$)
\be{ag40--}\begin{array}{lcl}
f &\la&  1_{A_j} (u_i) f(u_i) 1 (u_{m+1}) + 1_{A_j} (u_i) f(u_{m+1})\\
&& = 1_{A_j \cap B} (u_i) 1 (u_{m+1}) + 1_{A_j} (u_i) 1_B (u_{m+1}),
\end{array}
\ee
where now the two summands are alternative possibilities the
ancestral path can take and their respective probabilities associated
to the two new particles in each summand where one represents the preservation of the sample in
the other the insertion of a superior individual in the sample.

Recall now  the coupling
rule of the summands (the factors follow the transitions of the
individuals they are associated with). We now see that the
marked ancestral graph associated with this model has now the
property there is at most one path from the initial particle ''root'' to the
''leave'' at the other end, since now $\chi^j \cdot (1-\chi^j)=0$
and $\chi^j +(1-\chi^j) f \leq 1$.
The $\leq $ sign appears since it  might happen that not all
$n$-path end at the leaves, which means that the tagged sample
with the given type configuration and the generated genealogy is
not consistent and has probability zero.

In any case we can read each summand as a
possible line of ancestry and the factors generating the probabilities.
In the sequel we shall carry out and apply this construction at great length
and introduce {\em sums of ordered factors} to parallel the
genealogical relations.

Here we give an explanation how the genealogy-type structure is generated.
This means that we draw $n$ individuals from the population and
record their respective type and location and the time we have to go back for
every pair to find the most recent common ancestor. The joint law of this
statistics we want to derive from the dual process.
We start with the model (1) only resampling, (2) adding types and mutation,
(3) adding selection.
\sm

(1) The genealogy of a sample of $n$-individuals is obtained by considering
the time when pairs of two individuals of the sample first land in the same
partition element. This $n \choose 2$ different  random times are in
distribution equal to the entries in the matrix of genealogical distances
with the sample.

(2) Now the individuals carry a type and for a given genealogical tree
we have to calculate the probabilities that we find specific types on
the individual of the sample. For a realisation of the coalescent we
can, with the function-valued dual, determine the probabilities of a
specific type configuration at time $t$ given the one at time 0.
Namely with every partition element at every time $s \leq t$,
there is a connected factor a function on type space which changes
according to the mutation dynamic. Testing with the initial state
the factor turn 0 if the type is not possible and 1 if it is possible.
Note that if the set-valued process has reached the full or the empty set,
the process contains no further information about the time further back in the
original model.

(3) Finally selection enters into the picture. The selection transition
in the dual accounts for the possibility of an interaction of the sample
with the rest of the population.
In the neutral case this can be suppressed without changing the law
since all individuals are exchangeable if they are at the same site
at the moment resampling occurs. This is not true of course if we have type-based selection.
Each action of the selection operation
introduces two {\em alternative} forms of the genealogical tree of the
sample by the interaction with a new randomly chosen individual from
the population at the time the selection operator acts. Each realisation
of all mutation transitions, coalescence events for an $N$-coalescent
(in which the sample is embedded) $k$-selection transition generate
$2^k$-realizations of a marked random tree which are the potential
genealogy-type configurations for the sample at time $t$. For each of
those we can calculate the probability  from a decision tree
(whose leaves are possible genealogical trees)
and where the edges carry certain factors which change
according to mutation and selection.

A key point is now that all the possibilities are alternative, i.e. only one
such tree is the actually realized one for the sample. Furthermore the
structure is such that a finite random time back it is resolved which
of the different possibilities actually occurs.

The following diagram illustrates the decision tree represented in the dual.
Here the subscripts refer
to the order of the factors and the selection operator has acted
twice on the first factor and then the first and second factors
have coalesced. We also used  a coupling in which the operations
are simultaneously applied to the same factor (according to the
order) in each summand.  The result of this history has two
summands $(1-\chi ) \otimes f$ and $\chi f \otimes 1$.
The indices 1,2,3 refer to the dual particle involved.
This object can be viewed as a decision tree to decide which
ancestral paths are possible for the tagged sample represented
by $f$ and every leaf representing a possibility.

The following
picture indicates this but note depicted is the {\em decision} tree and {\em not} the genealogical tree.
In that decision tree we have first two selection events (the second applies to newborn particle)
and then coalescence of particle 1 and 2:
\bigskip

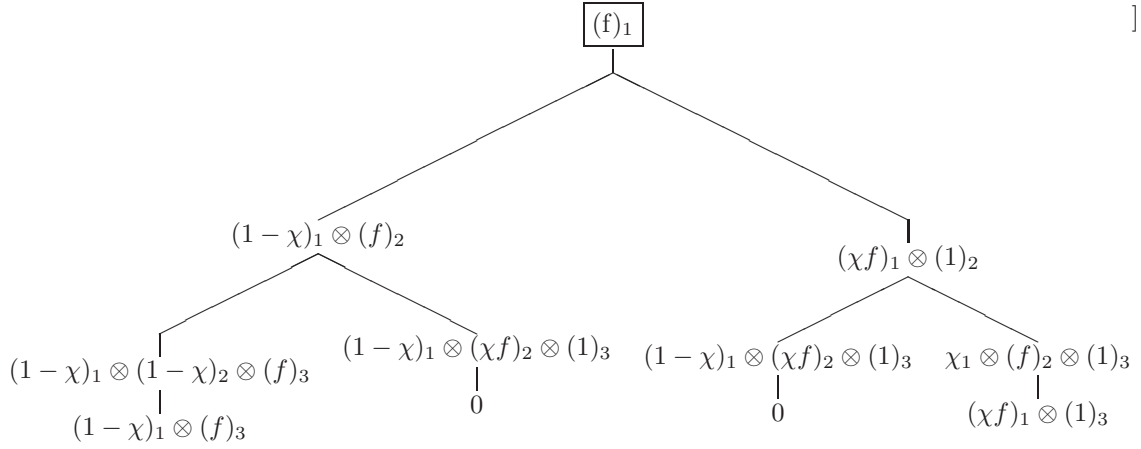
\Tree
 [.\fbox{(f)_1}
   [. [.{{$(1- \chi)_1\otimes (f)_2$}}
  [. [.{$(1-\chi)_1\otimes (1-\chi)_2\otimes (f)_3$}  {$(1- \chi)_1\otimes (f)_3 $} ]
    ]
    [.{$(1-\chi)_1 \otimes (\chi f)_2 \otimes (1)_3$} {$0$}
]]
   [. [.{$(\chi f)_1 \otimes (1)_2$}
   [.{$(1-\chi)_1 \otimes (\chi f)_2 \otimes (1)_3$} {$0$}
   ]
    [.{$\chi_1 \otimes (f)_2 \otimes (1)_3$} {$(\chi f)_1 \otimes (1)_3$}
   ]]
  ]]
]]
\bigskip

\section{Set-valued dual}
\label{s.svd}

\setcounter{secnum}{\value{section}} \setcounter{equation}{0}

We consider now the case where $\I=\{1,\cdots,M,M+1\}$
where $M+1$ is the type of maximal fitness. Our starting point will be the
function-valued dual processes $(\eta, \CF^{++}), (\eta, \CG^{++})$.

\subsection{Ordered function-valued dual}
\label{ss.ordfvd}

We introduce now an order into our dual process, based on the historical information
on the successive action of the selection operator. We treat $\CF^{++}$ and
subsequently $\CG^{++}$.

\subsection{Modification of dual 1: $ \wt{\CF}^{++,<}_t$ with ordered factors}
\label{ss.moddual1}

In this subsubsection we carry out the modification of the dual
$\CF^{++}_t$ constructed in Subsection \ref{ss.refine} involving an enrichment (see Remark \ref{R.dualf}).
In the case of our $(M,1)$-model we introduce here a simplified version applying for
the case of $M$ lower level types and one higher level type $M+1$ with only
up-mutation to the higher level.

\begin{remark}\label{R.dualf}
That enrichments of dual processes give again duality relations is a general principle.
For this we recall the general form of a duality (see (\ref{add100}) and
note that if we have for $X$ a dual process $Y$ with a duality function
$H(\cdot, \cdot)$ and if we can introduce a new process $Y^\ast$ on a new state
space $E^{\prime,\ast}$ such that there exists a map
\be{add118}
\kappa : E^{\prime, \ast} \to E^{\prime} \mbox{ such that }
\CL [\kappa (Y^\ast_t)] = \CL[Y_t];
\ee
then we have a duality with duality function
$H^\ast : E \times E^{\prime,\ast} \to \R$ by setting
\be{add120}
H^\ast (\cdot, \cdot) = H(\cdot,\kappa(\cdot)).
\ee

\end{remark}

\begin{remark}
Typically such enrichments arise from a duality  relation of the new dual
to an enriched original process. We no not use this fact here and
hence do not prove such a relation here.
\end{remark}

We proceed in six steps.
\bi

{\bf Step 1:} \quad {\em Preparation: further observations on the nature of
$\mathcal{F}^{++}_t$ .}
\smallskip

As a preliminary step to the construction of the modified dual process
that we shall use in the case $M\geq 2$ we now look in more detail
at the structure of the function-valued part $\mathcal{F}_{t}^{++}$ of
the dual process.
 For that purpose we view $\CF^{++}_t$ not only as a certain function
 but we explicitly work with its form as a sum of products of
 indicators, in other words this form will become part of the state,
 see (\ref{ga8b1}) below.

 This information is based on historical information about the process $\eta$,
 namely the information about the complete ancestral relations in the
 particle system $\eta$.

 First consider the special case $M=2$. This can be expressed as follows.
Define the abbreviations:
\bea{ga7b}&&f_{1}={(110)},f_{2}={(100)},f_{3}=(010),f_{4}={(000)},\\&& f_{5}={(111)}, f_6=(001), f_7=(011),f_8=(101).\nonumber
\eea
Let $\mathbb{F}=(f_1,\dots f_8)$  denote the set of possible factors.

Then starting with $\mathcal{F}_{0}^{++}= (110)$ and applying selection and the lower level mutation operations we have
\be{ga8b}
\mathcal{F}_{t}^{++} = \suml^{\wt N_t}_{i}
\bigotimes^N_{j=1}
f_{1}^{\otimes k_{1}(i,j,t)}f_{2}^{\otimes k_{2}(i,j,t)}f_{3}^{\otimes k_{3}(i,j,t)}
f_{4}^{\otimes k_{4}(i,j,t)} := \suml^{\wt N_t}_i \CF^{++}_{i,t}, \ee where
$\wt N_t$ denotes the (random) number of summands and
\be{ga8b1}
k_1,k_2,k_3,k_4: \{1,\cdots,\wt N_t\} \times \{1,2,\cdots,N\}\times
[0,\infty) \longrightarrow \N_0, \ee are appropriate (random)
functions which depend on $t$, for example, $k_1(i,j,t)$ is the
number of $f_1=(110)$ factors in the ith summand at time $t$ located
at $j\in\{1,\dots,N\}$. We can encode all information by  coding
$\CF^{++}_t$ in the form of a different type of object (misusing notation since
we refrain from using a new letter, since we go even further below introducing
$\CF^{++,<}$, where we explicitly spell out the dynamics:
\be{ga14aa}
\CF^{++}_t
 = (\wt N_t,\{(k_1(i,j,t),\cdots,k_4(i,j,t)); i \in
   \{1,\cdots,\wt N_t\},  j \in \{1,\cdots,N\}\}).
\ee
We note here that whenever a birth event in $\eta$ occurs and the selection
operator acts on $\CF^{++}_t$, we can get new
summands.

\begin{remark}
Note we could have if $1_A f \equiv 0$ that
\be{ga8b0}
f \la 1_A f + (1-1_A) \otimes f = (1-1_A) \otimes f
\ee
and even so a birth occurs we get no new summand.
\end{remark}

We can perform  calculations by expressing the dual
expectation in terms of a suitably chosen population of factors, which
gives the dynamic of the $(k_1(i,j,t),\cdots, k_4(i,j,t))$. In
particular we now consider the population of summands and the
dynamic of the exponents in each summand.
\bi

{\bf Step 2:} {\em (The basic idea: passing to tableaus)}

The key trick is to order factors and introduce factors 1. This is  based on the following observation.
We can also keep track of some historical information in order to
deal with the sum  on the r.h.s. of (\ref{ga8b}).
In particular we introduce an {\em ordering of
the factors} that allows us to {\em associate factors in different
summands}. First note that a product of indicators can change into
a {\em sum} of products of indicators only by the selection
mechanism and all other transitions preserve the product structure
of the indicators of a summand in the $\CF^{++}_t$ respectively $\CG^{++}_t$ dynamic.

Recall next that selection operates by
$1_A \cdot f + (1-1_A) \otimes f$ with $1-1_A$
corresponding to a new variable. Observe first that
we can write without changing integrals but having the same
number of variables for both summands as
\be{ga14aaa}
f \to 1_A \cdot f\otimes 1 + f \otimes (1-1_A),
\ee
where with  $M=2$ types we get $1_A=(011)$. Then it is useful to place the
new factors 1 and $(1-1_A)$ not at the end but next to the factor giving
birth.
\begin{remark}
Recall that in the $\CG^{++}$
we have  {\em reordered} the factors so that  in all
summands  $1_A\cdot f$
and $(1-1_A)$  always are associated with the {\em same} variable, we call this
having  the same {\em rank}. Then we can couple
the operations on different summands in a different way so that both are operated on
simultaneously by the mutation, selection and migration operators
acting at a given rank (and not particle).
This will be exploited again in the next subsubsection.
\end{remark}

From (\ref{ga14aaa}) we see that
we get a sum of two terms now involving a new
variable corresponding  to the particle added by the birth event
in the dual particle system. As time goes on this produces a {\em
binary tree} and each path in this tree starting from the root and
ending in a leaf corresponds uniquely to a summand.

In general and more precisely we make the resulting two summands  comparable in the sense that they
have  the same number of variables in both summands in
(\ref{ga16c}), but without changing the value. To do this we write this action of the selection
operator as follows. Take a function $f = f^3 \otimes f^1 \otimes f^2$, where
$f^1$ is a function of {\em one} variable, the one on which selection acts
 and $f^2, f^3$ of the remaining variable ordered in the order of the corresponding
 partition elements in the ordered particle system. Write:
\be{ga16d}
 f\to f^3 \otimes (1_A\cdot f^1) \otimes 1 \otimes f^2
 + f^3 \otimes  f^1 \otimes (1-1_A)\otimes   f^2,
\ee
where $1$ stands for $(111)$, in general $(1 1 \cdots 1)$ with $M$ ones,  which integrates to $1$ for any
probability measure. This means that we have inserted a new
variable in {\em both} summands. This allows us then to associate
factors in the two summands after the selection jump one to one starting from the initial
factor and ending with the last born variable.

To formalize the above ideas we need two changes.

(1) We want to view the function $\CF^{++}_t$ not just as a function but to add as part of its description  some additional
information on its {\em form}. Namely, we want to consider a sum of products of indicators where each factor is associated with a particular particle in $\eta$
and the particles are ordered in a certain way related to their appearance as selection acts.

(2) This means we can enrich the $\CF^{++}_t$ to a marked {\em tableau}
whose rows correspond to summands and columns to factors
(which are indicators) each of which is a function of one variable.

For this purpose we need  to {\em order the factors} and to
introduce {\em new factors of 1} where necessary in order that
all the summands in the expression for $\CF^{++}_t$ consist of an equal number of factors.
This leads to a collection of factors which are ordered, carry a location and
are organized in summands. This object we will call $\CF^{++,<}_t$.
In the next two steps we introduce this new process rigorously.

{\bf Step 3} \quad {\em State space of the ordered dual $\mathcal{F}^{++,<}_t$}

Next we introduce the ingredients necessary to construct
the {\em modified dual}, $\mathcal{F}^{++,<}_t$ formally.
This means that we add to the particle system a further function
assigning each particle a rank, an element in $\N$ which defines
automatically an additional order relation.
(Note that we use here a {\em different} order than in
Subsections \ref{dualp}-\ref{modified}.)

Start by observing  how selection acts in order to see how the order
must be set up.
An initial set of factors $f^1\otimes f^2\dots\otimes f^n$
each corresponding to {\em one} single variable is given a linear order (from left to right).
The new particles appear  by birth naturally ordered  in time.
However we now associate with a factor a {\em rank} which introduces then
a {\em new order} by rank. The rank is defined  dynamically as follows.
\sm

\beD{D.order} {(Order relation among factors, ranks, transition under selection)}

(a)
We start with giving rank $1,\cdots,N_0$ to the $N_0$ initial particles and the associated factors.

New factors are created by the selection operator as follows.

Suppose there is a birth by selection operating on the $j$-th particle in $\eta$ and
the corresponding factor is in the $\ell$-th rank and $m$ denotes the current number
of partition elements. That is,  selection hits
the factor with rank $\ell $  and produces the transition of
$f^1 \otimes \cdots \otimes f^m$ to:
\be{ga8b2}\begin{array}{l}
f^1 \otimes \cdots f^{\ell-1} \otimes (f^\ell 1_A) \otimes 1
\otimes f^{\ell+1} \otimes \cdots \otimes  f^m\\[1ex]
\hspace{4cm}
+ f^1 \otimes \cdots f^{\ell-1} \otimes f^\ell  \otimes (1-1_A)
\otimes f^{\ell+1} \otimes \cdots \otimes f^m.
\end{array}
\ee
In order to work with an ordered object we use the rule that
the additional factor 1 is
placed directly to the right of the factor on which the operator
is acting and the remaining factors are shifted one unit to the
right while $1$ and $1-1_A$ are put at the $(\ell+1)$-th position.

(b)
We denote this {\em order} of factors of one variable by
\be{xg43}
<, \mbox{ that is, } f^1<f^2
\ee
means that the factor  $f^1$ lies to the left of the factor $f^2$ in this order.

(c)
We shall relabel all factors counting left to right and
we assign the factors the
\be{fa8b3}
\mbox{{\em ranks} } 1,2,3,\cdots
\ee
with the natural order on $\N. \qad$
\end{definition}
Note that in this labelling the {\em order in rank is inherited}.
Therefore in (\ref{ga8b}) each of the $\sum_{\ell=1}^{4}
k_\ell(i,j,t)$ many factors has a rank in the linear order. This means that for each summand the factors are assigned ordered by their ranks.
This leads then to the following description.

\beD{D.orddual} {(State description of the ordered dual $(\eta^\prime,{\mathcal{F}}_t^{++,<})) $}

(a)
The state of the ordered dual has the form of a {\em marked  tableau}, where the tableau is
given by
\be{agd1}
\{\varphi_i(k); \quad i=1, \cdots, \wt N_t; \quad k=1,\cdots,N_t\},
\ee
where  rows corresponding to the index $i$ represent summands, and  columns corresponding to
the index $k$ represent  the factors of a given rank in different summands
 furthermore every column is assigned a location in
$\{1, \cdots, N\}$ by the $\varphi_i$ as well as the type of factor in position $(i,k)$ in
the array and assigning the rank to dual particles.

A row of the tableau can be split into a product of a set of sub-products of factors,
one sub-product for each occupied location.

We denote (a factor is a function $\I \to \R^+$)
\be{agd2}
\F = \mbox{ the set of possible factors}.
\ee

(b)Then the state at time $t$ can be uniquely described
by a collection of
\be{adg3}
\wt N_t \mbox{  summands (rows) of } N_t \mbox{ factors (columns)}
\ee
and the state is denoted
\be{ga14-a}
\{\varphi_i, \quad i \in \{1,\cdots,\wt N_t\}\},
\ee
where each summand is associated to a map
\be{ga14}
  \varphi_i: \qquad \{1,\cdots,N_t\} \la \F \times \{1,\cdots,N\}\times\mathbb{N}
\ee
specifying  separately for each summand
 for each factor  the {\em type} of the factor, the  {\em location } in geographic space and finally
{\em rank} in the order called $<$. We write
\be{agd4}
 \varphi_i(k)=(f_k,j_k,\ell_k)\quad, \quad k=1, \cdots, N_t,
\ee
with the constraints that
\be{agd5} f_k=f_{k^\prime}\text{  and  } j_k=j_{k^\prime}\text{  if  }\ell_k=\ell_{k^\prime}.\ee

We set
\be{tab1}
\CF^{++,<}_t =
\{\varphi_i(k)=(f_k,j_k,\ell_k): k\leq N_t, \;i=1,\dots,\wt N_t\}.
\ee

(c)
At time $0$ we set $\wt N_0=1$ and the particles $\{1,\dots,N_0\}$ are assigned an initial index and the particles $\{1,\dots,N_t-N_0\}$ are indexed in order of their birth times, the first birth assigned index $N_0+1$, etc. The rank assigned to a particle  changes dynamically  due to coalescence or other births to be described by the dynamics below
after (\ref{ga8c2}).

(d)
The set of particles assigned the same rank
corresponds to the partition elements defined in Subsection \ref{dualcomp}.
As in Subsection \ref{dualp} the set of partition elements at time $t$ is denoted
\be{agd6a}
\pi_t = (\pi_t(1), \cdots, \pi_t(|\pi_t|))
\ee
and the partition elements have the form
 \be{adg6}
 \pi_t(\ell)=\{k:j_k=j,\ell_k=\ell\},\quad \ell=1,\dots,|\pi_t|,
\ee
where the partition element $\ell$ is located at site
\be{adg7}
 \xi_t(\ell)= j.
\ee
Then $\pi_t$ defines a mapping
\be{adg8}
\pi_t:\{1,\dots,N_t\}\to \{1,\dots,|\pi_t|\}
\ee
where the partition elements are now {\em indexed} by the smallest index of the particles they contain and $\xi_t$ defines a mapping
(giving the locations of partition elements)
\be{adg9}
\xi_t:\{1,\dots,|\pi_t|\}\to\{1,\dots,N\}.
\ee

We define the enriched (from $\eta$) particle system
\be{b5b}\eta^\prime_t=(N_t,\pi_t,\xi_t, \Re_t),
\ee
that is,  $\eta^\prime_t$ given by the set of particles, their partition structure,
locations together with the
list of current ranks $\Re_t$ of individuals.
$\qquad \square$
\end{definition}

\bi

{\bf Step 4:}\quad {\em Dynamics of ordered dual process $(\eta^\prime, \CF^{++,<}_t)$.}

In order to clarify the ideas we start in {\bf (i)} with an informal description and then
we formalize this in {\bf (ii)} writing down a Markov jump process by specifying
the state, transitions and their rates.

{\bf (i)} The point will be to place the factors in the summands of
$\CF^{++,<}_t$ in a particular order useful for analysing the resulting
expressions. We describe now informally the transitions of the dual ordered in this way.

The dual evolution can be described as follows. First the particle system $\eta$.
Starting with $N_0$ initial dual particles each located at one of the sites in $\{1,\dots,N\}$  at time 0,  new particles are created by a pure birth process with birth rate $sm$ when the current number of dual particles is $m$. More precisely, independently each particle  gives birth at rate $s$. Each particle is located at one of the points in $\{1,\dots,N\}$ and is assigned in
addition a {\em rank}.
Then the enriched evolution with the transitions by coalescence birth (selection)
and mutation is as follows.

{\em Coalescence} (two particles at the same location) {\em effectively}
decreases the number of factors even though here we keep formally the number
of factors and introduce an additional factor 1.  We make the convention that upon
coalescence of two factors we replace the one with lower rank by
the product of the two factors  and the one with upper rank
becomes $(111)$ (this is the analogue of the wellknown
``look-down'' process). In formulas, at the coalescence of individuals
$i^\prime $ and $k^\prime $ in the particle system $\eta$ which corresponds
in a particular summand to factors with ranks
$i$ and $k$ (in the labelling (\ref{fa8b3})) we get for
the function-valued part the transition given by
\be{ga8c2}
f^1 \otimes \cdots
\otimes f^i \otimes \cdots \otimes f^k \otimes \cdots \la f^1
\otimes \cdots \otimes (f^i f^k) \otimes \cdots \otimes 1 \otimes
\cdots
\ee
with each factor we associate a rank and the new factor 1
gets the same rank which introduces then a new order by rank. The rank is defined
as in a  selection operation the newly created factor.

{\em Mutations} act as before on factors  corresponding to each variable
independently and we note that in particular  acting on $(111)$
or (000) has no effect.

{\em Selection} creates a new particle in the dual particle system $\eta$
and at the same time the transition given in (\ref{ga14aaa}) occurs.

The rank allows us to consider for a selection transition
the new factors now next to each other, i.e. we have
if the factor $f$ gives birth now after the transition the
sum of the two factors:
\be{ga8c3}
1_A \cdot f \otimes 1, \quad f \otimes (1-1_A).
\ee
(Recall $1_A$ and $1-1_A$ are connected with different variables).

{\em Ranks} change at the same time as follows.
The offspring of a particle of rank $k$ is assigned the rank $k+1$ and the
ranks of all particles  with ranks $\ell\geq k+1$ are reassigned rank $\ell+1$.  We denote the resulting number of particles in the dual system $\eta_t$ at time $t$ by $N_t$ where $N_t-N_0$ denotes the number of births.
Moreover at a birth time the offspring of a particle is located at the same site as the parent.

To complete the description of the dual we must define the corresponding function $\mathcal{F}_t^{++,<}$
which is a function of $N_t$ variables but has the form of a sum of {\em ordered}
products of factors.  The transitions of $\mathcal{F}_t^{++,<}$
are as listed above.

We can now represent the dual expression $\CF^{++}_t$ via our tableau. By construction we have:
\beL{L.reprtab}{(Representation via tableau)}

Consider the dual process $(\eta^\prime_t, \CF^{++,<}_t)_{t \geq 0}$.
Let $pr$ denote the map associating with $\eta^\prime$ the triple arising by ignoring the rank.
Given the tableau associated with $(\eta^\prime, \CF^{++})$ as
$\{\varphi_i(k)=(f_k,j_k,\ell_k): k\leq N, \;i=1,\dots,\wt N\}$,
we define for this state
\be{b5bf}
(\eta, \CF^{++}) =  \big(pr(\eta^\prime), \sum_{i=1}^{\wt N}
\prod_{j\in S}\prod_{\ell=1}^{|\pi|} (1_{\xi(\ell)=j}
\left( f_\ell(u_{\ell}))\right)\big).
\ee
Then we obtain a version of $(\eta_t, \CF^{++}_t)_{t \geq 0}.  \qad$
\end{lemma}

\begin{remark}
Note that we can also construct a version of
$(\eta^\prime_t, \CF^{++,<}_t)$
from the path $(\eta_s, \CF^{++}_s)_{s \leq t}. $
\end{remark}

It is often convenient to order the summands in the tableau in a specific way,
i.e. instead of a set as in (\ref{ga14-a}) we want to use a tupel.
Recall that the action of selection (corresponding to a birth event for $\eta_t$) on a function $f$ at a given location and rank is defined by (\ref{ga8c3}).
This replaces a summand by two summands (one of which can be $0$) in  $\mathcal{F}^{++,<}$  and we must keep track of the summands.

In order to keep track of the summands and produce a convenient tableau  we adopt the following convention for ordering the summands.
Starting with a single factor $f$ operated on by selection as in (\ref{ga8c3}) we produce the ordered rows
\be{add121a}\begin{array}{cc}
       \chi_A\cdot f & 1 \\
       (1-\chi_A) & f
     \end{array}
\ee
In general
\begin{itemize}
\item The operation of selection at a rank is applied successively to each row starting at the top and then moving down to each original row.
When it acts on the rank at a row it produces an additional row immediately below the row on which it acts. If it produces a
row with a zero factor this row is removed.
\item Mutation and coalescence can produce a row with a zero factor which is then removed but otherwise they do not change the order of the rows.
\end{itemize}

\sm

This means that we want to think of the
\be{add121}\begin{array}{l}
\mbox{state of $(\eta^\prime, \CF^{++,<})$ as a {\em tableau} where {\em rows} correspond to summands and {\em columns} to}\\
\hspace{3cm}
\mbox{variables and {\em marks} on the entries indicate the location of the factor}.
\end{array}
\ee
The ordered summands allow to put the entries 1 in convenient positions as we shall see later on
working with this object.

\begin{remark}
We will denote this object, misusing notation a bit, again with
$( \eta^\prime_t, \CF^{++,<}_t)_{t \geq 0}$
and use whatever appearance is most convenient.
\end{remark}

\begin{remark} A natural way to keep track of the summands is as the set of leaves in  a tree.  When selection event occurs at a vertex  two additional edges and vertices are produced in the tree. The dual representation will then involve an expectation over the summands indexed by the set of leaves in the tree.  This means that only the marginal distributions of the summands are involved, not the joint distribution of the summands.
One possible choice for the dynamic is that different summands evolve independently. But then we must have a separate birth process (associated to the selection events) for each summand. Instead in the construction below we have a single birth process but this process is applied simultaneously to all summands following the dynamics described below.
\end{remark}

{\bf (ii)} Based on the ideas introduced above we now formally define the dual process
$(\eta^\prime,{\mathcal{F}}_t^{++,<})) $ in the  case of $M$ types at the lower level and one type
at the higher level by specifying transitions and transition rates.

\bi

We next define the  {\em dynamic of the ordered dual} first in $\eta^\prime$ and
then in $\CF^{++,<}$, however this is only the formal version of what we described
in point (i).

We specify the transitions
that result in an increase or decrease of the number of factors, namely, selection and coalescence at one site.   This will produce an  enriched version of the  birth and death process  used in the case $M=1$, $\wt N_t$.   We now describe the transitions in
the dynamic step by step, first for $\eta^\prime$ in (1) and (2) and then for
$\CF^{++,<}$ in (3).

{\em (1)}  The pair $(N_t,\pi_t)$ of processes is defined as follows.  $N_t$ is a pure birth process with birth rate at time $t$ \be{b1x}s\cdot|\pi_t|\ee and $|\pi_t|$ is a birth and death process.

We label  the particles in the order
explained above, that is,  the offspring which results when selection acts on rank $\ell$ is
placed in the linear order at rank $\ell+1$ and the ranks of particles to
the right are shifted one place to the right, that is, their rank is increased by $1$.
Each new particle produced when there are $n$ particles is
the offspring of  a randomly chosen rank  $\ell$ which is chosen
with the uniform measure on $1,\dots,n$.

{\em (2)} The dynamics of the ranks is as follows.
Deaths of ranks at a site are due to coalescence and migration and occur at rate
\be{b2}
d\left(\begin{array}{c}
  m \\
  2 \\
\end{array}\right)1_{m\geq 2} +c m 1_{m\geq 1},\quad (1_{m>1}\;\text{for }N=\infty \mbox{ in the latter}).\ee

\bigskip

On coalescence of ranks $\ell$ and $\ell^\prime>\ell$ which are located at the same site
the particles with rank $\ell^\prime$ are  removed and the
other partition element at rank $\ell$ now includes all particles that had been in the two corresponding partition elements.

When a death occurs, that is a rank is removed from a site,  it is a migration with probability
\be{b3}
\frac{c1_{m\geq 1}}{cm 1_{m \geq 1} +d\left(\begin{array}{c}
  m \\
  2 \\
\end{array}\right)}, \quad (1_{m>1} \mbox{ for } N=\infty \text{ with the above convention}).
\ee
The rank  of a migrant is chosen with uniform
measure on $1,\dots,|\pi_t|$.

\begin{remark}
(a) The ranks $(\ell,\ell^\prime),\;\ell^\prime
>\ell$ involved in a {\em coalescence} are chosen with uniform measure on the
set of $m(m-1)/2$ pairs of particles. Note that when $m=1$ deaths due to coalescence do not
occur.

(b) Also we can adopt the convention that deaths at a site where $m=1$  due to migration do not occur in the case of infinitely many sites and initial measures of the form $\mu^{\otimes \N}$, that is, case we adopt the   convention that ranks at singly occupied sites do not migrate. Note that we cannot use this convention in the collision regime.
\end{remark}

{\em (3)} Turn now to the transition of $\CF^{++,<}$.

We first give the description of the transition that results from selection operator associated to $1_A$ acting on a rank $\ell$ at a
site $j$. We first consider the action on just the $i$-th summand, namely that this summand
is modified and a new summand is produced we denote by $V$ and we describe below how to label
the summands and the newly created ones.  The  original summand is modified as follows:

\be{b3b}\begin{array}{ll}
&\{\varphi_i(k)=(f_k,j_k,\ell_k)\}_{k=1,\dots,N_\ast}, \mbox{ as original state}.\\
&\text{number of particles changes:}\\
&N_\ast\to N_\ast+1,\\
&\text{function change to:}\\
&\varphi_i(k)=(f_k,j_k,\ell_k),\quad \text{ if } \ell_k< \ell, \\
&\varphi_{i}(N_\ast+1)=(1,j_k,\ell+1),\\
&\varphi_i(k)=(f_k,j_k,\ell_k+1),\quad \text{ if } \ell_k> \ell, \\
&\varphi_i(k)=({1}_Af_k,j_k,\ell_k),\; \text{  if  } \ell_k=\ell.
\end{array}
\ee

\noi
Moreover for every summand a new offspring  summand $\varphi_{i^\prime}$ (if we use the special order of
summands from (\ref{add121a}), $i^\prime = i+1$ and the number of all other higher summands is shifted
by one) is produced which is defined as follows:

\be{b3b2}\begin{array}{ll}
& \{\varphi_{i^\prime}(k)=(f^\prime_k,j^\prime_k,\ell^\prime_k)\}_{k=1,\dots,N_\ast+1}\\
&\text{is defined by}\\
&\varphi_{i^\prime}(k)=(f_k,j_k,\ell_k),\quad \text{ if } \ell_k< \ell,
\\
&\varphi_{i^\prime}(k)=(f_k,j_k,\ell_k+1),\quad \text{ if } \ell_k> \ell,
\\
&
\varphi_{i^\prime}(k)=(f_k,j_k,\ell_k),\; \text{  if  } \ell_k=\ell\\&
\varphi_{i^\prime}(N_\ast+1)=(1-1_A,j_k,\ell+1),\; \text{  if  } \ell_k=\ell.
\end{array}
\ee

The transition due to coalescence of ranks $\ell$ and $\ell^\prime>\ell$ at a site $j$ is given as follows:
\be{b3bcc}\begin{array}{ll}
&\{\varphi_i(k)=(f_k,j_k,\ell_k)\}_{k=1,\dots,N_\ast}\\&
\text{changes to }\\&
\varphi_i(k)=(f_k,j_k,\ell_k),\quad \text{ if } \ell_k< \ell,\\&
\varphi_i(k)=(f_\ell\cdot f_{\ell^\prime},j_k,\ell_k),\quad \text{ if } \ell_k= \ell,\\&
\varphi_i(k)=(f_k,j_k,\ell_k),\quad \text{ if } \ell< \ell_k<\ell^\prime,\\&
\varphi_i(k)=(f_k,j_k,\ell_k-1),\quad \text{ if } \ell_k> \ell^\prime,
\\&\varphi_i(k)=(f_\ell\cdot f_{\ell^\prime},j_k,\ell),\quad \text{ if } \ell_k= \ell^\prime.
\end{array}
\ee

Next we consider the transition that occurs if rank $\ell$ migrates to a new  site
say $j^\prime$:

\be{b3bcm}\begin{array}{ll}
&\{\varphi_i(k)=(f_k,j_k,\ell_k)\}_{k=1,\dots,N_\ast}\\&
\text{changes to }\\&
\varphi_i(k)=(f_k,j_k,\ell_k),\quad \text{ if } \ell_k\ne  \ell,\\&
\varphi_i(k)=(f_k,j^\prime,\ell),\quad \text{ if } \ell_k= \ell.
\end{array}
\ee

The collection of summands has a natural tree structure which is described as follows.   We begin with
\be{agd9} \wt N_0=1\text{   and give this index }1.\ee
We index
\be{agd9b}
\mbox{the offspring of a summand of index $i$ by $i1,i2,i3,\dots $}
\ee
which then defines a natural tree structure.  In this way we see that after $N_t-N_0=n$ births, there are $2^n$ summands.  However we note that some summands can be $0$ so that this provides an upper bound on the set of non-zero summands.

\bi

\begin{remark} \quad Assume we can work with the collision-free
regime $(N=\infty)$.
We can  construct a richer labelling system that incorporates both
the {\em birth order} and {\em rank} in the linear order, as
follows.

It can be encoded starting with the founder at $0$. The
history of transitions up to time $t$ and therefore the state at
time $t$ is determined by a sequence $\{0\rightarrow,\dots\}$
of transitions, where each transition is given by one of
\be{b4}
k\rightarrow,
k\leftarrow \ell \,(\mbox{ where } \ell>k) \mbox{ and  } k\uparrow.
\ee
The first
corresponds to a selection event and produces an offspring at
$k+1$, the second a coalescence (of $\ell$ with $k$) and the third an emigration, i.e.
migration to an new unoccupied site.  Note
that the founder position $0$ never migrates. We denote the
position of the $kth$ particle born in the linear order at time
$t$ by $v(k,t)$. We can have $k\leftarrow \ell $ only if
$v(\ell)>k$, that we, have a {\em lookdown}.
The relative order of the particles does not change but when a lookdown occurs we denote
the outcome as the ``compound particle'' $[k,\ell]$, etc.

For example

\be{b5}
0,0\rightarrow,0\rightarrow,1\rightarrow,0\rightarrow,1\rightarrow,2\leftarrow
3 \ee would result in (the labels according to birth times)
\be{b6}
0,01,021,0213,04213,042153,04[2,3]15.
\ee

Note that in
the special case with no coalescence this labelling above identifies
the tree associated to a branching process and in the special case
with no births identifies the coalescent graph (via the lookdown
process). Therefore this can be viewed as a coding of the
branching coalescing (or selection-mutation) graph.

The following observation for the enriched process is important.
Since the birth and death process is recurrent, when the number of
particles at a given location returns to one, the sequence can be
deleted and relabelled $0$ and the process begins again. \bi
\end{remark}

{\bf Step 5} {\em (Summarizing the construction)}

We can summarize the construction so far as follows.
We have modified the $(\eta_t, \CF^{++}_t)_{t \geq 0}$
to a new process $(\eta^\prime_t,  \CF^{++,<}_t)_{t \geq 0}$
where $(\eta^\prime_t)_{t \geq 0}$ is described by the enriched
(by information on the {\em rank}) birth
and death process (branching coalescing particle system) and the function-valued part is given by
\be{AAG5}
( \CF^{++,<}_t)_{t \geq 0}, \quad
\CF^{++,<}_t = \suml_{i=1}^{\wt N_t}  \CF^{++,<}_{i,t}
\ee
where $\wt N_t$ denotes the number of  summands and the
(possibly also ordered) summands are ordered products of factors.

The  {\em tableau} provides a nice way to to represent the state of $\wt \CF^{++,<}_t$
which is an array of strings of factors called {\em rows}.  The {\em columns} in this array have a nice structure
and correspond to factors with the same rank.
For that purpose we had written at each selection event the new rows
directly under the parent row. (See (\ref{mtd0}) for example).

For convenience we denote the tableau associated to the  object
$\CF^{++,<}_t$ as
\be{ga16caa}
\mathfrak{t}((\eta^\prime_t, \CF^{++,<}_t)),\ee
where a {\em column} corresponds to a {\em variable} and a {\em row} to a particular
{\em summand}.

\begin{definition} (Tableau-valued process)

The evolution of the Markov jump process $(\eta^\prime_t, \CF^{++,<}_t)_{t \geq 0}$
induces a {\em marked tableau-valued}  pure jump process driven by an
autonomous particle process $\eta^\prime_t$:
\be{ga16cb}
(\mathfrak{t}((\eta^\prime_t, \CF^{++,<}_t)))_{t \geq 0}.
\qquad \square
\ee\textbf{}
\end{definition}

We have for the modified process according to the Remark \ref{} the duality relation:
\beP{P.moddual} {(Modified duality )}

We have $x^\otimes$ abbreviating the initial product measure state.
\be{add122}
E[ \int  \CF^{++}_t d x^\otimes] =  E[ \int \CF^{++,<}_t d x^\otimes].\qad
\ee
\end{proposition}

The main idea for the analysis of $(\eta^\prime, \CF^{++,<})$ will
be to complete the specification of  the dynamic in such a way that the duality remains
true for the new modified function-valued dynamic  as well as
getting cancelations of terms by coupling the summands for
different $i$ in the representation in (\ref{ga8b}) by obtaining sums
combining to give factors 1.  What is in the way is that upon selection $1_A$
and $1-1_A$ have different ranks. We defined in Section \ref{basictool}  a  new
process denoted by $(\CG^{++}_t)_{t \geq 0}$ and  we define the enriched
$(\CG^{++,<}_t)_{t \geq 0}$ precisely
in the next Subsubsection \ref{sss.crittimesc}
by specifying the mutation and selection transitions  such that
they occur simultaneously at the same rank in each summand.

{\bf Step 6} {\em Simplifying the ordered dual to $(\eta^\prime, \wt \CF^{++}_t)$}.

Our main application of the dual is to determine the emergence time scale. For this reason
we make a few simplifications that are applicable for the calculations of the moments of the type $M+1$ when we begin with an initial measure on $(\Delta_M)^N$ of the form
\be{adg10}
\mu^{\otimes N} \mbox{ with } \mu(M+1)=0.
\ee

Note that for the moment we have included summands that when
integrated against $ \mu^{\otimes N} $ yield $0$, hence the dynamic can be
further simplified if we are interested only in the
integral on the r.h.s. of (\ref{ga9b2}). We observe that a summand
which contains a factor
\be{ga14b}
(0,\dots,0, \ast), \mbox{ with } \ast = 0 \mbox{ or } 1
\ee
becomes $0$ when integrated w.r.t. $\mu$ since we start
the original process with all the mass on types $1,\cdots,M$ so that a
summand with such a factor would not contribute taking the
expectation. Therefore we can and shall  agree to prune the
summands in the dual process by setting
\bea{ga14c}
&&\wt {\mathcal{F}}^{++,<}_{t,i}=0, \;\;\mbox{ summand } i \mbox{ contains a
factor } (0,\dots,0,\ast)\\&& \wt {\mathcal{F}}^{++,<}_{t,i}=
{\mathcal{F}}^{++,<}_{t,i},\; \mbox{ otherwise} .\nonumber \eea
 With this convention we have to deal only
with summands which have a contribution in the integral in
(\ref{add122}).

Similarly the evaluation of ranks for which  $f=(1,\dots,1,1)$ in all summands does not change
the integral on the r.h.s. of (\ref{add122}). However for the moment we do keep
factors which are effectively one, i.e. $(1,\dots,1, 1)$ since below we shall
see that this will help doing the bookkeeping once later on we couple
summands.

Furthermore note that the factor (0,\dots,0,1) can only gain further $1^\prime$'s by
a downward mutation and at the moment we consider
$m_{3,\rm{down}}=0$. (Even for $m_{3,\rm{down}}>0$ this would at
best become visible if the dual process reaches a state where we
have $O(N)$ summands with such a factor).

\beL{L.dual1} {(Simplified duality with ordered factors)}
For initial measures satisfying $\mu(M+1)=0$, we have the dual representation in terms of the simplified dual:
\be{ga16d2}
E [\int \CF^{++}_t d(\mu^\otimes)] = E [\int \wt
\CF^{++,<}_t d (\mu^\otimes)]. \qquad \square
\ee
\end{lemma}

\subsection{Modification of the dual 2:  $\mathcal{G}^{++,<}$ with ordered coupled summands}
\label{sss.crittimesc}

In the dual $(\eta^\prime, \CF^{++,<})$ factors and summands are ordered and part of the state
and summands are coupled,
since transitions in the particle system induce transitions
in $\CF^{++,<}_t$, which occur for certain factors in {\em all} the summands.
The strategy of this section is to to use the {\em new coupling} between
summands used by $\CG^{++}$, that took advantage of the structure of this dual
which introduces {\em cancellations} which will make  combinations of  some summands
in one variable equal to 1. This can be generated using the process $\CG^{++}$ in
combination with the enriched particle system $\eta^\prime$ and the order
introduced for $\CF^{++,<}_t$ to arrive at
\be{add201}
(\eta^\prime_t, \mathcal{G}^{++,<}).
\ee

This construction will be the basis for the second key idea that will be taken up
in the next Subsection \ref{ss.svdf}, which
is to introduce another dual object that
allows us to represent the process with coupled summands as a nice
{\em set-valued  process}
which allows to study the longtime behaviour. In our application
we can also calculate moments and with enrichments of it we
evaluate the Malthusian parameter $\beta$ that determines
the critical time scale.

The objective is to use the basic duality formula

\be{ga9b2}
E[\int\mathcal{F}_{0}^{++}d (x(t))^{\otimes}]    =E[\int\mathcal{F}_{t}
^{++}d (x(0))^{\otimes}] = E\left[\int {\CF}^{++,<}_td(\mu^\otimes)\right],
\ee
\bi

\noindent where $\mathcal{F}_t^{++}$ is the refined dual process
of Subsection \ref{ss.refine} and to rewrite the r.h.s. further,
in the sense that for every initial state $x^\otimes$:

\be{equiv}
 E\left[\int {\mathcal{F}}^{++,<}_td(x^\otimes)\right]
 = E\left[\int {\mathcal{G}}^{++,<}_td(x^\otimes)\right].
\ee

{\bf Step 1} {\em (Construction of the new couplings of summands in the ordered dual)}

We obtain the {\em new} coupling of the summands on the level of the tableau-valued process by
changing the rule of the selection
transition, parallel to the change in Section \ref{basictool} from $\CF^{++}$ to $\CG^{++}$, as follows.
The action of selection operator acting on
a particle with rank $\ell$ and associated factor $f$ and
location $j$ produces now at $j$ and the ranks $\ell$ and
$\ell +1$ the factors induced by

\be{ga16e}
f \la 1_A f \otimes 1 + (1-1_A) \otimes f.
\ee

This means that now the factors $1_A$ and $(1-1_A)$ have the
{\em same} variable, i.e. they are in the same column and the same rank.
In particular this means that under the further action of the migration,
coalescence or the selection
operator we  apply this simultaneously to the  $(1-1_A)$ and $1_A$ in the
product $(1-1_A) \otimes f$ and $1_A\cdot$  $f \otimes 1$,
that is,  to every factor in the corresponding rank (in the linear
order). This includes the transitions of these two factors in the
two summands arising from coalescence, migration, selection and
mutation, both normal and rare but for the moment we do not
consider the rare mutation.

Note that at this moment the variables of $1_A$ and $(1-1_A)$ are associated
by this modified mechanism, even though originally they where associated with different
individuals in $\eta$. The point is that taking the expectation of the sum
 w.r.t the dual dynamics  we
have {\em not} changed the value  of the expectation since the expectation depends only
on the marginals corresponding to the summands and not the joint law.
This means that we permute the order of the individuals in the $\eta$-process
depending on the summand, in order to make the cancellation effects
explicitly visible.

\begin{remark}
Note that $\CG_t^{++,<}$ have the same type of state as ${\CF}^{++,<}_t$
but the dynamic is now changed.
\end{remark}

\begin{remark}
A key point is that with this new coupling
the different (i.e. not repetitions of another one)
non-zero summands always correspond to {\em disjoint}
events. Also note that some summands can correspond to the empty
event which does not contribute in the duality.
\end{remark}

\begin{remark}
Note in the new process $(\CG^{++,<}_t)_{t \geq 0}$ can be viewed as a process,
where the  birth,
coalescence and migration transitions are now jumps associated with the ranks (recall
(\ref{AAG5})).
Recall that the {\em coupling} of summands in (\ref{ga16b}) was given by
the evolution rule that every transition occuring for an individual
in $\eta$, say $k$  of some particular rank  say $\ell$  occurs simultaneously
in each summand  that is for every   $ \wt{\CF}^{++,<}_{i,t}$ in (\ref{AAG5}).
 This induces a coupling of the rows in the tableau.
 Alternatively we can use  $\eta^\prime$ as the driving process
by permuting for each individual in $\eta^\prime$ the position in the ranking
to the factor it is associated with.
 \end{remark}

Recall that rows which contain factors $(0\dots 0)$ play no role in the duality
we do not need them in any calculation and we can completely remove them
which we did in passing from $\CF^{++,<}_t$ to $\wt \CF^{++,<}_t$.
Similarly a column which has only factors $(1 \cdots 1)$ plays no role
in the duality expression.
We call these ranks {\em inactive}.

\beD{D.nonzero} {(Non-zero summands and active ranks)}

(a) If we include only summands which are not 0,  we number these as

\be{dgr3} 1,\dots, \wt N^\ast_t.
\ee

(b) Columns (ranks) in the tableau which contain only $(1 1 \cdots 1)$
are called inactive. In particular sites which are only
colonized by inactive rows are called inactive.
$\qquad \square$
\end{definition}

\beD{D.G++} {(The coupled dual process $(\eta^\prime,\CG^{++,<})$)}

We shall denote the process arising from $(\eta^\prime,\mathfrak{t}((\eta^\prime \wt \CF^{++,<})))$
by introducing coupled summands via the rule (\ref{ga16e})
and by removing inactive ranks and zero rows by
\be{ga16e1}
(\eta_t^\prime, \CG^{++,<}_t)_{t \in \R}.\ee
The part $\CG^{++,<}_t$ of the  process consists of the collection of summands
of marked factors
\be{ga16e2}
\{\mathcal{G}^{++}_{i,t}:i=1,\dots,\wt N^\ast_t\},
\ee
each consisting of $N^\ast_t$ active ranks and the associated factors with their location,
thus specifying again a {\em marked tableau} of indicator functions
\be{adg11}
\mathfrak{t} (\CG^{++,<}). \qquad \square
\ee
\end{definition}
A key property of the new dual is that it allows us in a more transparent way
to keep track of how far our expression deviates from a product, which we
summarize in the next remark.

\begin{remark}\label{R.selec}
Recall from above that the sum can be organised with the help
of a tree, a splitting occurring with each operation of the
selection operator. Consider at time $t+s$ the two summands
corresponding to each of the two subtrees starting at a birth time
$s$ with selection operator  $ (1_A f \otimes 1 f)$,
respectively  $f \otimes (1-1_A)$  (where the order of the factors
is as indicated left to right). We have changed the dynamic such that we
couple the summands such that effectively we work instead of this with
\be{ga16ca}
1 \otimes (1_A f)  \mbox{ and } f \otimes (1-1_A),
\ee
 which is then a {\em decomposition in complementary terms } since
$1_A +(1-1_A)=1$. Furthermore if mutation acts on the partition $A$
the decompositions $(\emptyset,\I)$ respectively $(\I,\emptyset)$
are traps and if the mutation rates are strictly positive these traps are
actually reached. If the trap is reached (we call this resolution)
then only one of the two summands in (\ref{ga16ca}) remains.
\end{remark}

We have the identity (which we prove in Step 2 below):

\beP{P.dualcoupled} {(Tableau-valued duality with coupled summands)}

The process $(\eta^\prime, \CG^{++,<})$ gives the same expectation as the dual
$(\eta, \CF^{++,<})$, namely:

\be{ga16c1}
E\left[ \int \left(\suml^{\wt N_t}_{i=1}
 {\CF}^{++,<}_{i,t}\right)d(\mu^\otimes)\right] =
E \left[\int \suml^{\wt N^\ast_t}_{i=1} \CG^{++,<}_{i,t}
d(\mu^\otimes)\right] . \qquad \square
\ee
\end{proposition}

\begin{remark}\label{R.key}
We have a look at the behaviour of sites contributing in the duality.

In the two-type case we saw
new sites being colonized by particles was always equivalent of a
contribution of this site to the dual expression. This is not the
case for the case of $M>1$ and becomes very transparent with
the new dual. Consider the sequence of transitions

\be{xg61c0}\begin{array}{l}
f\\ \to \\ f\cdot1_B\otimes 1 +\\  1_{B^c}\otimes f
\\ \to \\
1_A 1_B\cdot f\otimes 1\otimes 1 +\\  1_{A^c} \otimes 1_B\cdot f\otimes 1+\\
1_A1_{B^c}\otimes f \otimes 1 +\\
1_{A^c}\otimes 1_{B^c}\otimes f.
\end{array}
\ee

If the second rank migrates and then the first rank
evolves by mutation via the jump $1_A1_B\to 1$ (hence
$1_A 1_{B^c}, 1_{A^c} 1_B, 1_{A^c} \to 0$),
then the newly occupied site becomes inactive (=1). Hence we see that factors which have colonized a new site need not
lead to a permanent colonization of that site as active site which is in contrast to the
case  $M=1$.

Note that the jump to 1 of a rank to the right of a migrating rank does not have this effect. For example consider
\be{xg61c}\begin{array}{l}
P \to f\cdot1_B\otimes 1 +  1_{B^c}\otimes f
\\ \to \\
1_A 1_B\cdot f\otimes 1\otimes 1 +\\  1_{A^c} \otimes 1_B\cdot f\otimes 1 +\\
1_A1_{B^c}\otimes f \otimes 1 +\\
1_{A^c}\otimes 1_{B^c}\otimes f
\end{array}
\ee
Now first let the first column migrate and then the second column resolve $1_B\to 1$. We then obtain

\be{xg61c2}\begin{array}{l}
1_A 1_B\cdot f\otimes 1\otimes 1 +\\  1_{A^c} \otimes  f\otimes 1+\\
1_A1_{B^c}\otimes f \otimes 1 \\=\\
1_A 1_B\cdot f\otimes 1  +\\  1_{(AB)^c} \otimes  f
\end{array}
\ee
but where now the first and second columns are located at different sites.  However the first column is not removed.
\end{remark}

{\bf Step 2} {\em (Alternative Proof of Proposition \ref{P.dualcoupled}).}

The statement can be proven viewing the process as enrichment of $\CG^{++}$ and
seeing that the latter satisfies the duality. We give an independent argument based
on $\CF^{++,<}$.

We have to recall how
we rewrote the state of  $\wt \CF^{++,<}_t$ (compare (\ref{AAG5}) and
(\ref{ga16c}) below). Consider the finite sum decomposition from
(\ref{AAG5})
\be{ga16b}
{\mathcal{F}}^{++,<}_t= \suml^{\wt N_t}_{i=1}
 \CF^{++,<}_{i,t}, \quad \mu^\otimes - a.s., \ee where each $\wt
\CF^{++}_{i,t}$ is an (ordered) product of indicator functions.

Given a realization of the particle system, we will build  a version of
$\wt{\mathcal{F}}^{++,<}_t$. Then in (\ref{ga9b2}) we can write for
the r.h.s.:

\be{ga16c0}
 E [ \int \left\{\mathcal{F}^{++,<}_{t}\right\}d(\mu^{\otimes})]
= \wh E\left[ \int \{ E_{{\wt N}_t}
{{\mathcal{F}}}^{++}_{t}\}d(\mu^{\otimes})\right], \ee where
$\widehat E$ denotes the expectation over the number of summands
process $(\wt N_t)_{t \geq 0}$ and $ E_{{N}_t}= E(\cdot| N_s:0\leq s\leq
t)$ the respective conditional expectation with respect to the mutation, coalescence and migration dynamics.

As in the case $M=2$ we have that $\mathcal{F}^{++,<}_t$
is decomposed into a sum  and then
\be{ga16c}
 E \int \left\{\sum^{\wt N_t}_{i=1} \mathcal{F}^{++,<}_{i,t}\right\}d(\mu^{\otimes})
= \wh E\left[  \int \{E_{\wt N_t}\left[\sum^{\wt N^\ast_t}_{i=1}
 {\mathcal{F}}^{++,<}_{i,t}\right]\}d(\mu^{\otimes})\right].
\ee
We note that we have the following upper bound on the number of  non-zero summands $\wt N^\ast_t\leq \wt N_t\leq 2^{N_t}$ (since selection acting on a rank can replace each summand by at most two non-zero summands).
We note that the sum over $\wt N_t$ or $\wt N^\ast_t$ results in the
same integrand for $\int$.

Therefore given $\{\wt N_s:0\leq s\leq t\}$ we have to keep track of
$\{{\CF}^{++,<}_{i,s}:{1\leq i\leq \wt N^\ast_s,\;0\leq s\leq t}\}$ and the way
they combine to form the sum. On the r.h.s. of (\ref{ga16c}) we see that the  conditional
expectation  depends only on the marginal distribution
of the vector of summands and {\em not} on the joint distribution for the given birth process.
This is where the coupling will come in.

Since in the dual representation we compute the expected value of
the sum this depends only on the marginal distributions of
the dynamics (conditioned on their starting points) of the different
summands after time $s$ (recall (\ref{ga16c}). Therefore we can {\em couple} the
evolution of the dynamics of the two   summands differently and still preserve
the duality relation.

Observe here that the different transitions of the dual $(\eta^\prime_t,
{\CF}^{++,<}_t)_{t \in \R^+}$ occur independently for each
individual in the dual process and the dynamic of the dual system
is Markov. Therefore if  our coupling of summands is constructed
by coupling transitions of factors in $\{{\CF}^{++,<}_{t,i}, \quad i =1,\cdots, N^\ast_t\}$,  which
correspond to different  leaves $i$  in the binary tree generated by the dual population the
basic duality relation  (\ref{ga16c1}) {\em will be
satisfied} automatically.

\subsection{The set-valued dual as functional of the ordered dual}
\label{ss.svdf}

We obtain now from $(\eta^\prime, \CG^{++})$ as a functional a process
$\wt \CG^{++}$ which  is set-valued and turns out to have a Markovian dynamic
which we exhibit.

\subsection{Dual $\wt \CG^{++}$: Examples for the special cases $M=2,3$}
\label{sss.examples}

To understand how this coupling construction of Subsubsection \ref{sss.crittimesc}
might be represented as set-valued process we consider
the example of $M=2,3$ and calculate what form the transitions take in
$\CG^{++,<}_t$ and indicate how we can represent the resulting process
which we  call $\wt \CG^{++}$
again in form of a  functional of a Markov process
and use this process as the new dual. These examples also suggest the concept
of {\em resolved} and {\em unresolved} factors, which is a key
to the growth of the number of factor in the dual expression.

{\bf The case $\mathbf{ M=2}$}

 In the case $M=2$ we construct a
Markov process with  dynamics looking as follows.
There are a finite number of types of factors, namely
$(100),(010),(110),(111)$. Starting with a factor of type
$(110)$ at a tagged site, new factors and a new summand are produced once a birth in $\eta^\prime$ occurs
by the action of the
selection operator. Factors can move to a new site,
can coalesce or change due to mutation. The mutation can result in the
removal of a summand as we shall see below, a phenomenon not occurring
with one type on the lower level.  We shall show below that
the factor-valued process  associated with
a given site is recurrent and returns either to the state
$(110)$ after a finite random time or to (000) and the whole summand is deleted.  During this sojourn some
factors emigrate to new sites and the founding factor at the new
site now (in contrast to  the case $M=1$)
may  become inactive  and thus a new site can be lost again.
The essential question is the rate of the {\em successful
emigration} (establishing a new permanent site)
which governs the growth of the number of the total number of (active) factors.

Here are the transitions in detail.

\begin{itemize}
\item \textbf{Selection}  We first consider selection.
Recall that the basic selection operator in $\CG^{++}$ is given by
\be{b81b}
f\to (1_A\cdot f) \otimes (111)+(1-1_A)\otimes f
\ee
where $1_A=(011)$ and the terms
$1_A\cdot f$ and $(1-1_A)$ have the same rank
in the order so that all other operations are performed
simultaneously on them.
Since $1_A, 1-1_A$ are of the form (011) respectively (100) they are not
traps under the mutation chain (different from (110) disregarding rare
mutation) and can be changed to $(00 \ast)$ or $(11\ast)$, which are
quasi-traps. We will therefore introduce the following concept:
\be{b82b}
(01 \ast), (10 \ast) \mbox{ are called {\em unresolved} factors}.
\ee

\item \textbf{Mutation}  We now include the action of mutation at
rate $(m_{12}+m_{21})$ which can {\em resolve} factors,
i.e. a transition of the selection factors (011) and (100)
to the quasi-trap (110) or to the quasi-trap (001) occurs:

\be{b9}(100)\to(110),\;(010)\to(000)\text{  with probability  }
\frac{m_{21}}{m_{12}+m_{21}},
\ee or
\be{b10}(100)\to(000),\;(010)\to(110) \text{  with probability  }
\frac{m_{12}}{m_{12}+m_{21}},
\ee
(thanks to the coupling we use).
This means that one summand is zero and the other has only
the factor (110), which is a quasi-trap.
Factors where the mutation has occurred are called ``{\em resolved}
factors''. These factors are important because they don't undergo
further changes by mutation (except rare mutation), where with unresolved factors in a row we do not
know yet if this summand remains or will become zero due to the action
of mutation.

\item \textbf{Coalescence}
   Now consider the coalescence of rank $k_1$ and $k_2$. If a rank $k_2$ looks-down to
    a position $k_1<k_2$, then the column $k_2$ is removed and the rows in which
   the combination $(010),\;(100)$ or $(100),\;(010)$ occur are removed (remember that here we are taking the product of indicator functions).
  Note that if subsequently the resolution
$(100) \to (000),\;(010)\to (110)$ occurs at the rank $k_1$,
then the offspring of both original ranks are deleted.

\item \textbf{Migration}  Any factor can undergo  migration,
that is, move to a new site, but factors of the same rank do this
together.  For the moment  we consider the case of infinitely many sites.  In this case a rank always moves to an unoccupied site (but note that the rank of the particles it contains does not change).

The importance of the location is that factors at different sites do not coalesce.
\end{itemize}

We now make a number of observations being crucial in the sequel.

Note that migrating factors can be one of three types (referred to as $\ast$-types)
\be{b83}
(110),(100)^\ast,(010)^\ast,
\ee
where for example $(100)^\ast$
denotes a factor that is currently given  as factor
$(100)$ later on the resolution $(100)\to (110)$ occurs,
which makes the summand with (010) equal to zero and the one
with $(100)$ (which became (110))  acts as the original factor (110).
Hence the $\ast$ indicates that the factor has not yet resolved in
$(110)$ or $(000)$.

Note that we can decide on the  final  type of a factor already at the time of migration with
probabilities $(\frac{m_{2,1}}{m_{1,2}+m_{2,1}})$,
$(\frac{m_{1,2}}{m_{1,2}+m_{2,1}})$.  Note furthermore that the choice of the
 final type  also determines the fate of the particles to the right of
that rank.

Note that it is possible for a newly occupied site to be deleted if the
founding rank is deleted due to resolution  $(100)\to (111)$
of an ancestral position (recall Remark \ref{R.key}).
In this case the descending newborn factors then have no effect
anymore, either a summand disappears or does not change the
founding father. Hence we have
to keep track of the ancestral relations in order to
determine which factors need to be deleted!

How do these transitions combine? Consider some examples.

Three successive applications of the
selection operator $1_A=(011)$  to a factor $f$, i.e. successively to rank 1 can be described in the form of
a tableau  with 8 rows and four columns.
Deleting now the zero rows we end up with the
simplified tableau of only four rows as follows:

\be{mtd0}
\begin{array}{llll}
   1_A^1\cdot f &\otimes 1 &\otimes 1 &\otimes 1 \\
  (1-1_A^1) & \otimes 1_A^2\cdot f & \otimes 1  &\otimes 1  \\
 (1-1_A^1) & \otimes(1-1_A^2) & \otimes 1_A^3\cdot f  & \otimes 1\\ (1-1_A^1) &\otimes (1-1_A^2) &\otimes (1-1_A^3)&\otimes f.
\end{array}
\ee
Here the first column corresponds to the initial
position of $f$ and the ``offspring'' of a selection operation is
placed immediately to its right. Each column consists of indicator functions which
induces a partition of the type space and the partition
elements evolve via the mutation process but in such a way that the
evolving partition elements at all time remain disjoint.

Applied to the case $f=(110)$, $1_A=(011)$,  (\ref{mtd0}) yields the following {\em tableau}
of indicators corresponding to subsets and inducing a partition
of $\{1,2,3\}^4$:
\be{mtd02}
\begin{array}{llll}
  (010) &\otimes 1 &\otimes 1 &\otimes 1 \\
  (100) & \otimes (010) & \otimes 1  &\otimes 1  \\
  (100) & \otimes(100) & \otimes (010)  & \otimes 1\\ (100) &\otimes (100) &\otimes (100)&\otimes (110).
\end{array}
\ee

Consider the case in which at the second rank in (\ref{mtd02}), we
have $(100)\to(000),\;(010)\to(110)$, by mutation. Then there are only two
summands and the third and fourth ranks  are inactive  (and therefore  columns)
are deleted.  This deletion of formerly  active ranks due to mutation is
the reason that the growth is slower than in the two type case.

{\bf Associated disjoint decomposition of product of type set.}

In (\ref{mtd02}) each row corresponds to one summand in the dual
$\mathcal{G}^{++}$ and each column defines a decomposition of the set of
the $M$ low level types. The collection of columns defines a
decomposition of the 4-fold product of the set of low-level types.
The number of factors is a pure birth
process with rate $s$ as in the two-type case.

More generally
each state of the $(k^\prime \times k)$-tableau induces a {\em decomposition of the
$k^\prime$-th product of the set of low level  types}, i.e. $\{1,2\}^{k^\prime}$
and thereby also of the set $\{1,2,3\}^{k^\prime}$.
The collection of {\em rows} of the resulting tableau corresponds for (\ref{mtd02})  to the
decomposition of $\{1,2\}^4$ given by the elements
$(x_1, x_2, x_3, x_4) \in \{1,2\}^4$:

\be{b8}\begin{array}{lll}
&&\{x_1 \in \{2\}\}\\&&\cup
\{x_1 \in \{1\}\}\cap\{x_2\in \{2\}\}\\
&&\cup\{x_1 \in \{1\}\}\cap\{x_2 \in \{1\}\}\cap\{x_3 \in \{2\}\}\\
&&\cup\{x_1 \in \{1\}\}\cap\{x_2 \in \{1\}\}\cap \{x_3 \in \{1\}\}\cap\{x_4 \in \{1,2\}\}.
\end{array}
\ee

{\bf The case $\mathbf{M=3}$}

We next consider the effect of having $M>2$ lower level types using the case
$M=3$ exhibiting the key features.
The main new feature is that with $M=2$ we had only one selection
operator, since we had two fitness levels. In general  $M-1$ of the $M$
lower level types  have different fitness levels from  1.
Having more such levels results in more  possible ways for  unresolved
factors to reach their final resolved state in a quasi-trap
called above $\ast$-types. Hence the essential difference
between the case $M=2$ and $M=3$ is the {\em resolution
mechanism} for the dual process which we shall describe in the case
$M=3$.

Suppose the fitness of the four types is
\be{xy2}
0,a,1,1 \mbox{ with } a \in (0,1). \ee Recall the upper level type
has downward mutation at rate $N^{-1}$ while the three others have
mutation rates $O(1)$.

We need to analyse the dual starting with one particle and with:
\be{du1} \CG^{++,<}_0 =(1110). \ee

\hspace{0.5cm}
Most of the things said for $M=2$ remain but what
we have to study in addition is what happens after the
creation of an unresolved factor by the action of a selection
operator, i.e. a factor not of the form $(111\ast)$ or $(000\ast)$
which are the quasi-traps, i.e. traps under mutation, excluding
rare mutation. On this unresolved factors various actions can take place,
mutations, further selection, coalescence
and migration. The new feature is that the selection can act at
different levels and  can be combined with intermediate mutation
steps. {\em Consequently the process of resolution  is no longer a 2 state
Markov chain} (with states $(10\ast)$ or $(01\ast)$,
but requires paths of intermediate steps.

We focus on this feature now and look at the transitions of
resolved and unresolved factors step by step.

\begin{itemize}

\item {\bf Coalescence and migration.} This is exactly as what we
described in the case $M=2$.

\item \textbf{Selection and mutation}

We focus first on the selection operator followed by
subsequent mutations.
We first note that we have two selection operators
corresponding to the sets of types
\be{du2}
A_2 \quad (\mbox{ fitness } \geq a) \mbox{ and
} A_3 \quad (\mbox{ fitness } \geq 1),
\ee
which leads at rate $s$ to the
multiplication by (0111) resp. (1000) and by (0011) respectively
(1100). Recall the choice between $A_2$ and $A_3$ occurs with
probability $a$ respectively $(1-a)$.
\bi

In particular we have for the selection operator with $A=A_3$:
\be{ga22c4} (1110)\to (0010) \otimes (1111)
+(1100)\otimes(1110)\quad \text{   at rates  } s\cdot (1-a).
\ee
After a finite
random time due to the action of mutation this yields either
$(1110) \otimes(1110)$ (mutation from type 3 to type 1 or 2) or
$(1110)\otimes (1111)$ (mutation from type  2 to type 3 and then
type 1 to 3) and this selection event is then in either case {\em
resolved} since we arrive in a quasi-trap under mutation.

If we use the selection operator $A=A_2$ we have the transition:
\be{ga22c5} (1110) \longrightarrow (0110) \otimes (1111) + (1000)
\otimes (1110)\quad \text{   at rates   } s \cdot a. \ee
In other words, selection occurs at rate $s$ and then the selection operator $(0111)$ is chosen with probability $a$ and the selection operator $(0011)$ is chosen with probability $1-a$.

Using the  coupling between summands in $\CG^{++}$,
eventually after a finite random time  and several
subsequent mutation steps (for example $3 \to 1, 2\to 1$ for the
$(0110)$ factor) the part of the sum due to the considered first action
of the selection operator collapses by mutations to either $(0000)$
or alternatively to $(1110) \otimes (1111) $ or $(1110)\otimes
(1110)$, and the factor on the r.h.s. of (\ref{ga22c5}) is
resolved.
\end{itemize}

What makes this dynamic difficult to handle?
The problem is to determine  what happens between creation of an
unresolved factor and its resolution by subsequent mutation steps.
As before when a position migrates we can decide in
advance, i.e. at creation, which state will occur at resolution.
But now to determine the probabilities of resolution in one of the
possible states, namely,  the two quasi-traps $(111{\ast})$ and
$(000\ast)$,  we need to follow  the {\em intermediate} states
which were not present for $M=2$.

To do
this we must work with a finite state Markov chain conditioned to
be absorbed by one of the two absorbing states. Having first
chosen this outcome we can then secondly replace the mutation Markov chain
by its {\em h-transform} corresponding to the condition to reach the
appropriate absorbing point.  The overall effect of this does not
change the law of the  resolution process.

For example the r.h.s. in (\ref{ga22c5}) can resolve to

\be{du3}
(1110)\otimes(1111)+ (0000)  \mbox{  mutations } 1 \to 3 \mbox{ or }1\to 2.
\ee
Similarly r.h.s. of (\ref{ga22c5}) can resolve to

\be{du4}
(000 \ast )+(1110)\otimes(1110) \mbox{ by mutations } 2\to 1,\; 3\to 1, \mbox{ etc.}.
\ee

The length of such paths to resolution is in principle {\em unbounded}
but is eventually absorbed in $(000\ast)$ or $(111 \ast)$
and can pass to a {\em finite} number of intermediate states.

 In other words between birth or coalescence events  each rank
involves the dynamics of a partition of the set of types
$\{1,2,3\}$. Moreover each partition element is associated to a
subset of the set of  columns  of the tableau of factors in the different
summands.  The dynamics of the
partition associated with a rank is given by a Markov chain induced by the mutation
process. Note that the number of partition elements can only
increase by a selection operation and can only decrease by coalescence
in each case by one element.  Otherwise the
partition process eventually ends in a quasi-trap and stays there
until the next selection event at this rank.

Note that a founding particle is not influenced by
coalescence at that new rank and its behaviour is given by a Markov
chain involving mutation and selection at that rank by the
look-down principle in carrying out the coalescence (this uses of
course the symmetry between individuals under the resampling
mechanism). We can use this as in the $M=1$ case to obtain a
simpler form of the sum. Namely we
 can determine its resolved value in advance and then follow the
h-transformed process as in the $M=2$ case.
As before we distinguish for each of the selection operators two
versions associated with the two principal states for absorption
but now instead of making the choice with probabilities
$m_{1,2}/(m_{1,2}+m_{2,1})$ respectively $m_{2,1}/(m_{1,2}+m_{2,1})$ we now use the two
absorption probabilities. The states of the process before
resolution is then described by the {\em $h$-transform}
corresponding to the absorbing state chosen.

\begin{remark}
In the new formulation of the dual process with coupled summands namely
$(\eta^\prime, \CG^{++})$ it is
important to  code the {\em new site} according to whether we have a
factor which survives or one which can be deleted after some
random time. In particular in the regime where we have no
collision we can do this easily since a site colonized will never
be hit by another migrant ever, so that these two events do not
depend on the evolution at other sites. For this purpose we shall later
introduce classes of permanent, removed, and transient
sites. A site is called permanent at the moment it is clear that  it
survives, it's called removed once it becomes deleted and transient
after the time of birth, until it will eventually be deleted.
Alternatively using the $h$-transform we can immediately at the
time of colonization decide to assign permanent or transient to the site.
\end{remark}

{\bf Example for coalescence of selection factors.}
In order to get a better feeling for the process, we discuss a particular effect.
Since there are more types we must also consider the result of
coalescence between two different partitions. As we has seen above
in the three type case it is necessary to determine the outcome of
collisions of unresolved factors of different types.

Now consider the case of selection operators corresponding to
$A=(0111),B=(0011)$ acting on two different ranks, say $i<j$, both
with factors corresponding to (1110):

\be{b11}
(1110)\rightarrow_{1_{A}}
\begin{tabular}
[c]{l}
$(0110)\otimes(1111)$\\
$(1000)\otimes(1110)$,
\end{tabular}
\;\;\;(1110)\rightarrow_{1_{B}}
\begin{tabular}
[c]{l}
$(0010)\otimes(1111)$\\
$(1100)\otimes(1110)$.
\end{tabular}
\ee
Then after coalescence of the first columns sitting and after insertion of two ranks
at ranks $i$ and $j+1$ we get the non-zero rows at rank $i$
\be{b12}
\\
\begin{array}{rl}
1_{A}1_{B} = (0110)(0010) &= (0010),\\
\;\;\;\;\;\;\;\;\\
(1-1_{A})(1-1_{B})  & =(1000)(1100) = (1000)\\
\;\;\;\;\;\;\;\; 1_A(1-1_B)&=(0111)(1100)=(0100).
\end{array}
\ee

The following illustrates how effect of  coalescing two columns followed by a  mutation.
Consider the tableau obtained starting with $(1110)$. Then let the  selection operator $1_B$
act on column 1 twice and then the selection operator $1_A$ acting on the third resulting column.
This yields the tableau on the left in (\ref{b13}).  Then the coalescence of the 3rd and 1st columns results in the
tableau on the right of (\ref{b13}).  Then mutation from 1 to 2 removes the last row and  mutation from 2 to 1
removes the third row, resulting in the total change.

\bigskip

\be{b13}
\begin{array}{|c|c|c|c|}
  0010 & 1111 & 1111 &  1111 \\
  1100 & 0010 & 1111 &  1111 \\
  1100 & 1100 & 0110 & 1111 \\
  1100 & 1100 & 1000 & 0110 \\
\end{array}
\qquad \mapsto \quad
\begin{array}{|c|c|c|c|}
  0010 & 1111 & 1111&  1111 \\
  1100 & 0010 & 1111&  1111 \\
  0100 & 1100 & 1111&  1111 \\
  1000 & 1100 & 1111&  0110 \\
\end{array}\ee


Then after resolution of the first column which would consist
either of (1110) factors or (0000) factors, either the third or
fourth rows  can be deleted since after
integration in the dual expression the result is not affected by
these terms anymore.
\bigskip

We note that very long strings on the way to resolution are very unlikely since
the mutation rates between the three lower types are strictly
positive and all other rates bounded above by some constant, so
that the quasi-traps for this mechanism are reached in finite
time. We also note that at each step we have a finite number  of
possible outcomes. Therefore by Markov chain theory the waiting
time needed to reach resolution and reach a quasi-trap  has
therefore {\em finite mean}.
\medskip

{\bf Conclusion from examples $M=2,3$:}

We can summarize the message of the two examples for the case of $M$-types
on the lower level as follows. Again as for $M=1$ we can simplify the
dual expression but instead of analysing strings of factors
we now work with {\em tableaus} of factors. In addition we have to
distinguish {\em resolved} and {\em unresolved} factors and hence resolved and unresolved sites.
In comparison with the case $M=2$, for the general case $M \geq 3$  there
are two new features in the resolution of the selection  factor:\\
(i) resolution does not occur at the first mutation, and\\
(ii) there are a finite number of  intermediate states, more precisely ($2^{M}-2)$ (here we have $M$ lower
level types)that can
be reached by mutation before resolution.

Nevertheless the observations on the form of $\CG^{++,<}$ in the two examples does
suggest that our {\em Markovian dynamics of marked tableaus}  corresponds to a {\em decomposition of subsets of $\{1,\cdots, M+1\}^\N$}. This we will now
formally introduce in Subsections \ref{ss.gensetup}-\ref{ss.svdualm} for the case of general $M$.

\subsection{The set-valued duality}
\label{ss.gensetup}

We are now ready to set up formally the representation of $(\eta^\prime, \CG^{++,<})$
by a marked {\em set}-valued Markov process and in particular specify its state space.
We make some definitions in a form, that it can be used
in fact on any geographic space, which is countable. All
the duality constructions we are carrying out in Subsections \ref{ss.gensetup}
-\ref{ss.svdualm} work for every finite type space.

Let
\be{setu1}
\I=\I_M:=\{1,\dots,M,M+1\},
\ee
\be{setu2}
\mathcal{I} :=
\mbox{  the algebra of subsets of } \I.
\ee
Let the geographic space be $S$ with
\be{setu3}
S=\{1,\dots,N\} \mbox{ or } S=\N.
\ee

The states of $(\eta^\prime, \CG^{++,<})$ can be coded as {\em marked} tableaus of columns of factors
of one variable marked by a location, where a row defines then the product of factors. This
defines a set-valued state
\be{setu4}
\wt \CG^{++}_t,
\ee
as follows.  The rows of the tableau (recall (\ref{b5bf})) correspond to the indicator function of subsets of
$(\I_M)^m$ for some $m$ and the product over the marks gives then a subset of $\prodl_{j \in S} (\I_M)^{m_j}$
such that  {\em distinct rows} correspond to {\em disjoint subsets}.
Therefore the sum of the rows is the {\em indicator function} of a subset of $\prodl_{j \in S}(\I_M)^{m_j}$.
It is often convenient to associate with this subset a subset of
$  \prod_{j\in S}(\I_M)^{\N}$ by considering
\be{add191}
\prodl_{j \in S} (A_j \times \I^\N) \quad , \quad A_j \subseteq (\I_M)^{m_j}.
\ee

\beD{D.codeG++}{(Marked set-valued process of $\wt \CG^{++}$)}

(a) The rows of the marked tableau $\mathfrak{t}(\CG^{++,<})$
define disjoint subsets of
\be{add123}
\prod_{j\in S}(\I_M)^{m_j}   \mbox{ for some } \{m_j<\infty\}_{j\in S}.
\ee
The union of these disjoint subsets is a subset of
$\Pi_{j \in S}({\I_M})^{m_j}$ whose indicator function is given by the sum of the rows
(where the factors in the row are multiplied corresponding to taking the intersection of the
sets associated with the different ranks)
of the tableau.

This subset of $\prod_{j\in S}(\I_M)^{m_j}$ does not characterizes the state
of $(\eta^\prime, \CG^{++,<})$ since in general
the ranks cannot be reconstructed from ranked subtableaus corresponding
to different sites in general.

(b) The set in (\ref{setu4}) defines a new process (the state space will be formally introduced
below) of subsets of $\prodl_{j\in S}(\mathbb{I}_M)^\N$ denoted:
\be{sub1a}
(\wt \CG^{++}_t)_{t \geq 0}.\qad
\ee
\end{definition}

\begin{remark}
Note that since the map giving $\wt \CG^{++}$ from the Markov process $(\eta^\prime, \CG^{++,<})$
is not bijective (at least if space contains more than one site, see above remark), the former need not be a priori Markov
in a spatial situation. However we will verify the Markov property below in Proposition \ref{P.MPsvd}.
\end{remark}

We now introduce some needed notation:
\be{tau.n}
 \mathcal{T}_{m} := \{\text{ the algebra of subsets of }(\I_M)^m\},
\ee
\be{tau.n2}
\mathcal{T}:= \text{ subsets of }(\I_M)^{\N}
\text{ of the form }A\times (\I_M)^{\N} \text{ with } A\in\mathcal{T}_m\text{ for some }m,
\ee
\be{tau.n3}
\mathcal{T}\subset \wh{\mathcal{T}}:=\sigma-\text{algebra of
subsets of }(\I_M)^{\mathbb{N}}.
\ee

The state space of $\wt \CG^{++}$ is contained in
\be{setu4b}   \textsf{I}_\ast= \text{Algebra }\bigotimes_{j\in S} \mathcal{I}_{j},
\mbox{ where for each } j,\;  \mathcal{I}_{j} = \mathcal{T}.
\ee

Given a non-empty set $G\in \textsf{I}_\ast$ let (recall we remove inactive ranks in $\CG^{++}$
and hence in $\wt \CG^{++}$ we get only finitely many factors not equal to $\I$)
\bea{supg}
|G|:= &&\min\{j: \exists S_j=\{s_1,\dots,s_j\}\subset S:
G= G_j\otimes (\I)^{S\backslash S_j} \text{ with }G_j\in \bigotimes_{i\in S_j} \mathcal{T}\}
\\&& :=+\infty \text{ if no such }S_j \text{ exists }.\nonumber\eea
If $|G|=j<\infty$, the {\em support} of $G$, $\rm{supp}(G)$  is defined to be the set $S_j$
that appears in (\ref{supg}).

Then we define
\beD{D.stspwtg}{(State space of $\wt \CG^{++}$)}

The state space of $(\wt \CG^{++}_t)_{t \geq 0}$ is defined as countable algebra of sets
\be{setu5}
\textsf{I}=\textsf{I}_S:=\{G\in\textsf{I}_\ast:|G|<\infty\}.
\qquad \square
\ee
\end{definition}

Given the finite type spatial Fleming-Viot process with selection
$ X^S_t=\{x^N(t,j)\}_{j\in S}\in (\mathcal{P}(\I_M))^S$, we define the expression
\be{setu6}
\wh X^S_t :=\prod_{j\in S}(x(t,j))^{\otimes\N}\in (\mathcal{P}((\I_M)^\N))^S,
\quad\text{ with  }S= \text{ the countable geographic space}.
\ee

The following dual representation is satisfied because of the very construction
of the process $\wt \CG^{++}$ as functional of the dual $\CG^{++,<}$:

\beP{P.msvdr}{(Set-valued duality relation)}

\bea{DGD}
E_{X_0} \left [<\wh X_t,\wt{\mathcal{G}}^{++}_0>\right ]
= E_{\wt{\mathcal{G}}^{++}_0} \left [<\wh X_0,\wt{\mathcal{G}}^{++}_t>\right ],
  \;\forall\;\wt{\mathcal{G}}^{++}_0\in \textsf{I}\text{  and  }X_0\in (\mathcal{P}(\I))^S.
  \qad
  \eea
  \end{proposition}

The point in the sequel will be to formulate a Markovian dynamic which generates
the process $\wt \CG^{++}$ in a transparent way and hence preserves  the
above duality relation. The Markovian dynamics of the process $\wt{\mathcal{G}}^{++}_t$
is constructed  in Subsubsections \ref{sss.prepsvd}-\ref{ss.svdualm}first without
migration, then incorporating the latter.
Then indeed we can  prove in:

\beP{P.MPsvd}{(Markov property of $\wt \CG^{++}$)}

The process $(\wt \CG^{++}_t)_{t \geq 0}$ is Markov. $\qquad \square$
\end{proposition}

In order to analyse the behaviour of the r.h.s. in (\ref{DGD})
for our questions on emergence and fixation we have
to classify factors of the geographic sites as follows.

\beD{D.resfac} {(Occupied, inactive and resolved sites)}

(a) A factor is called resolved if its indicator has the form $(1, \cdots, 1, \ast)$.

(b)
A site is occupied if there is at least one  factor with indicator not equal to $(1,\dots,1,\ast)$ at the site.
An occupied site  is called {\em resolved} if it contains only
one type of factor, namely with indicator $(1,\dots,1,\ast)$,  or {\em inactive}
if it contains only factors with indicators  $(1,1,\cdots,1)$. $\qad$

\end{definition}

\subsubsection{Preparations: Set-valued Markovian dynamics without migration}
\label{sss.prepsvd}

We now have to find a mechanism which turns $\wt \CG^{++}$ into a {\em Markov} process.
As a preparation for this analysis of the  set-valued dual  dynamic $\wt \CG^{++}_t$ we
begin by focussing on a single factor with subsets of $\I$ (corresponding in $\CG^{++,<}$
to a single column (rank)) since in this case we shall see that
$\wt \CG^{++}$ is in fact a Markov process. The difficulty arises due to migration.
First we analyse
the behaviour of $\wt \CG^{++}$ under the action of the mutation transition
at a single factor at a single location.
Later we consider also selection and coalescence, in particular
the interaction of mutation with selection.

Where we had before with $(\eta, \CF^{++})$ or $(\eta, \CG^{++})$
a collection of
birth and death processes with values in $\N$, we now get a dynamic on another
countable set, which can be viewed as subsets of
$\I_M$. Acting on each column we have a dynamic of subsets
of $\I_M$ if $j$ is the number of factors $\I$.
We study here these refined dynamics
and begin by looking at the dynamics at a given factor $\I$ (i.e. rank),
first under the effect of mutation, then only under selection and finally under the
combination of both. At last coalescence is added in a straightforward way.
At a single location we deal with subsets of $\I^\N$ as in (\ref{tau.n2}). Our subsets
typically do not have product form but are disjoint unions of subsets of products of
$\I$. So the disjoint union induces a partition of $\I$ in each factor of $\I$. Therefore
we will study the dynamics of these partitions of $\I$.
We introduce migration in the next subsubsection.
\bigskip

{\bf (1) Mutation}
Consider the mutation mechanism in the dual process
$(\wt \CG^{++}_t)$, ignoring for the moment all other transitions. This induces transitions on indicators of subsets
of $\I$ - see (\ref{b9}), (\ref{b10}), (\ref{du3}), (\ref{du4}) for examples.  In general, translating the transitions of  indicator functions into transitions on sets we obtain
for $j \in \I$:
\be{dgr62}
A \to A\cup\{j\}\text{  at rate  } \sum_{\ell\in A}m_{j,\ell},\quad A\to A\backslash \{j\}\text{  at rate  } \sum_{\ell\in A^c} m_{j,\ell}.
\ee

This defines a {\em Markov process}  on subsets of $\I$ which we denote by (here $m$ stands for the
matrix $(m_{i,j})_{i,j=1,\cdots,n}$):
\be{b18a}
(\zeta^m(t))_{t \geq 0}.
\ee

We want to define now a dynamic on a whole {\em partition} of $\I$.
Given a {\em partition of} $\I$  we can
let the mechanism in (\ref{dgr62}) act simultaneously on  {\em all} the indicator
functions corresponding to the elements of a given partition.

Note that for our model the order $O(1)$ transitions do not change the $(M+1)$ component so  that effectively this is described by a Markov process on the set of partitions of $\{1,\dots,M\}$ with
transition and transition rates given for every $k=1, \cdots, n$ by:

\be{add65}
(A_1,\dots,A_n)\to (A_1,\dots,A_k\backslash\{j\},\dots, A_\ell\cup\{j\},\dots,A_n)\text{   at rate}
\sum_{j^\prime\in A_\ell}m_{j,j^\prime}.\ee

This dynamic of partitions of the set $\I$ can be immediately generalized to a dynamics of
{\em tuples of sets} which are either disjoint or repetition
of another subset. (Note the union can be strictly contained in the type set).
Given a $k$ disjoint subsets with union $\{1,\dots,M\}$  the process from (\ref{dgr62}) induces a Markov jump process
\be{b18a2} (\wt \zeta^m (t))_{t \geq 0},
\ee
with values in $k$-tuples of disjoint possibly empty subsets with union $\{1,\dots,M\}$. Similarly we can start with the elements of any partition of $\{1,\dots,M\}$.

Then we can prove

\beL{mutationlemma} {(The set-valued mutation dynamics and h-transform)}

(a)
The Markov process $(\wt \zeta^m(t))_{t \geq 0}$ is a pure Markov jump process with state space given by the
set of partitions of $\{1,\dots,M\}$. The partition
\be{b18c}
(\I_M\backslash \{M+1\},\emptyset)
\mbox{ is a trap}
\ee
Assume that the mutation matrix is strictly positive
for all types $i,j \in \I_M\backslash\{M+1\}$. Then
starting with an initial partition, exactly one initial partition element
will grow to $\I_M\backslash \{M+1\}$ and all others will
go to the empty set.

(b)  Given the initial partition $(\{1\},\dots,\{M\}))$ of
$(\I_M\backslash\{M+1\})$ consider the process
 \be{b18d} (\wt\zeta^m_{t}(1),\dots,\wt\zeta^m_t(M)),
\ee
where the $\{\wt\zeta^m_t(i), i=1,\cdots,M\}$ are possibly empty
disjoint subsets  with union $\I_M\backslash\{M+1\})$. The quantity $\wt\zeta^m_t(i)$
denotes the set of initial types that have mutated to type $i$.
This is a Markov process.

If the $m_{k,\ell}>$ for $k,\ell\in \{1,\dots,M\}$, then for every
$k \in \N, 1 \leq j \leq k$ the $k$-tuples
\be{b18e}
A_j^\ast:=(\emptyset,\dots,\emptyset,\I_M\backslash\{M+1\},\emptyset,\dots,\emptyset),
\mbox{ with } \I_M\backslash \{M+1\} \mbox{ at } j\mbox{-th position}
\ee are absorbing points for the the $\wt\zeta^m_t$ process and with probability one the process reaches an absorbing
point at a finite time.

(c) Given $1\leq j\leq k$ the process conditioned to hit the absorbing point $A_j^\ast$ is
given by the $h$-transform  corresponding to $h_j$ which is
the Markov process  with transition rates \be{htransform}
q^{h_j}_{\zeta_1,\zeta_2}=
\frac{h_j(\zeta_2)}{h_j(\zeta_1)}q_{\zeta_1,\zeta_2}\quad\text{if
} h_j(\zeta_1)\ne 0, \ee where
\be{b20}
h_j((B_1,\dots,B_k))=P_{(B_1,\dots,B_k)}(\{\wt\zeta\text{ hits  }
A_j^\ast \mbox{ eventually}\}).
\ee

Denote the law of the $h$-transformed process by $P^h$.
Then the law of $\wt\zeta^m$ is given by \be{hjl}\sum_j h_j P^{h_j}. \qquad \square \ee

\end{lemma}

\begin{proof}{\bf of Lemma \ref{mutationlemma}}

(a) Note first that in a mutation transition two types $i,j$ are involved
and hence it suffices to consider the effect for a set $A$ and $A^c$.
Consider therefore the initial condition
\be{xg62}
\psi(0)=(A,A^c)$, $A\subset \I_M\backslash \{M+1\}.
\ee
Then a mutation from type $i$ to type
$j$ makes no change if $i,j\in A$ or $i,j\in A^c$.  However if
$i\in A,\;j\in A^c$, then after the transition at time $\tau $,
\be{xg63}
\wt \zeta(\tau+)=(A\backslash i,A^c\cup i).
\ee
Therefore after each
mutation, we have a  partition therefore resulting in a
partition-valued Markov process.

The process continues until it reaches
\be{xg64}
\wt \zeta^m=(\emptyset,\I_M\backslash \{M+1\}) \mbox{ or }
\wt \zeta^m=(\I_M\backslash \{M+1\},\emptyset).
\ee
Therefore given any partition exactly one partition elements grows to
$\I_M\backslash \{M+1\}$ and all the other go to the empty
set since we have assumed positive rates on $\I_M \backslash \{M+1\}$.

We can next consider the enriched process where we
start with a complete partition of $\I_M$. When type $i$ changes to type
$j$, $i$ is moved to the partition element containing $j$. This
process continues until there is only one non-empty partition element (these
are the traps).  Note that the number of different non-empty
partition elements is non-increasing and the process eventually
hits a trap with probability one.

(b)  and (c)
Relations (\ref{htransform}) and (\ref{hjl}) follow as a special case of Doob's
h-transform of a Markov process (\cite{RW}, III.45).  q.e.d.

If we now consider $j$-column, the mutation acts independently on each of
these components. Therefore we have now all
what we need to handle the mutation part of the evolution of the
set-valued tableau.
\end{proof}
\bigskip

\textbf{(2) Selection.} Given a  factors $k$,
the action of selection on this factor produces a new rank
(birth) placed at $k+1$ and the rank $k+\ell$ is
moved to rank $k+\ell+1$ for $\ell = 1,2,\dots $.

Now consider the action of a selection operator  on an indicator
function at a factor at a site.  Recall that selection acts on each factors with rate $s$
and when this occurs   the selection operators corresponding to ${A_i},\;i=2,\dots,M,$ is chosen
with probabilities $(e_i-e_{i-1})$ where $\{e_i\}$ are defined in (\ref{Gre3}).

Then the action of the chosen selection operators $1_{A_\ell}$ act by
multiplication  on the  active elements of
the factor and intersects with $A_\ell^c$ on the new elements of the factor as specified by (\ref{b3b}), (\ref{b3b2}).
On the level of a set-valued dynamics, if $B=\cup_{i=1}^m \cap^n_{k=1} B_{i,k}\in \mathcal{I}^n$,
then selection $1_{A_\ell}$ acting on factor  $j\in\{1,\dots,n\}$
produces a transition to the disjoint union  of ordered products

\bea{dgr6}
&&\left(\bigcup_{i=1}^m \left(\prod_{k=1}^{j-1}B_{i,k}\times (B_{i,j}\cap A_\ell)\times \I \times
\prod_{k=j+1}^n B_{i,k}\right) \right)
\\&& \qquad \qquad \bigcup \left(\bigcup_{i=1}^m \left(\prod_{k=1}^{j-1}B_{i,k}\times A_\ell^c \times \prod_{k=j}^n B_{i,k}
\right) \right)\in \mathcal{I}^{n+1}.\nonumber
\eea
We can represent such a subset in a more transparent way in a reduced tableau consisting of subsets of
$\I$ as entries, where the columns form a tuple of subsets of $\I$ taken from a partition of $\I$
and the rows define disjoint {\em nonempty} subsets of $\I^{n+1}$.

This means we can define uniquely the {\em evolution of a column} of a tableau since the indicators of sets, say
$A, B$, in a column satisfy either $A \cap B=\emptyset$ or $A=B$.
Given a specified rank in a tableau at a site the above mechanism
induces  the evolution  of the column  at this rank where of course in the case of
repetition of sets they follow all the same evolution.

\begin{remark}
Note that the corresponding tableau can be rewritten taking advantage of the simplifications $1_{A_i}+1_{A_i^c}=\I$,
$1_A1_C=0$ if $A$ and $C$ are disjoint, which may reduce the
number of rows, but of course not the number of columns.
\end{remark}

\begin{remark}
  Consider the result of two selection operations acting on
the same rank - we obtain four summands
\be{xg54}
(f 1_{A\cap B}) \otimes 1 \otimes 1,
\quad (1-1_A)1_B  \otimes f\otimes 1,
\quad (1-1_B) \otimes 1_Af \otimes 1,
\quad (1-1_B)  \otimes (1-1_A)\otimes f.
\ee
corresponding to
the decomposition of $\I_M$. If $f \equiv 1_\I$ the decomposition is given by:
\[({A\cap B})\cup ({A^c\cap B})\cup (B^c\cap A)\cup ( B^c\cap A^c).\]
The further evolution under mutation then determines  which of the four
decomposition elements takes over, that is, which of these will grow to $\I_M\backslash \{M+1\}$.
For general $f$ we get more decomposition elements.
\end{remark}

\textbf{(3) Combined effect of selection and mutation.}  The effect of selection is to produce new
nontrivial factors (real subsets of $(\I)$.
 Since mutation can change factors, it can change such a factor into
$\I$ or $\emptyset$  - when this happens the corresponding particle dies.
We can then think of the selection process as a ``proposal'' process and then mutation is acceptance-rejection process as follows (in the case $M=2$)
\be{add66}\begin{array}{lcl}
\{1,2\}& \longrightarrow & \{2\} \times \{1,2,3\} \cup \{2\} \times \{1,2,3\}\\
& \longrightarrow & \{1,2\}  \times \{1,2,3\} \cup  \{1,2\}  \times \{1,2\}.
 \end{array}
\ee
In the first case the new rank produced by selection if removed (rejected) by mutation.

\bigskip

\textbf{(4) Coalescence.} Factor $\ell$ and
$\ell^\prime>\ell$ which are at the same site coalesce by a look-down process: the  rank
$\ell^\prime$ is removed and the factor $\ell$ is modified by taking the intersections of the
corresponding partition elements - see (\ref{b3bcc}). (As a result in the representation
as irreducible tableau some rows and
columns can disappear. For example, in the case $M=2$ the row with
$(100)$ in the $\ell^\prime$th rank and $(010)$ in the $\ell^\prime$ th
rank is removed on coalescence so that one row and one column
is removed.)
\bi

\textbf{(5) Set-valued Markovian dynamic.}
Now we have defined all three mechanisms occuring in a single site model and
we can combine this now to a set-valued evolution  in $\CT$.
We denote the corresponding set-valued dynamics (respectively partition-valued)
combining mutation (rate $m$), selection (rate $s$),
coalescence (rate $d$) as
\be{b18a3}
(\zeta^{{\mathrm{m,s,d}}}_t)_{t \geq 0}, (\wt \zeta^{\mathrm{m,s,d}}_t)_{t \geq 0}.
\ee
The corresponding jump rates are denoted
\be{jr1}
\{q^{\mathrm{m,s,d}}_{\zeta,\zeta^\prime}\}_{\zeta,\zeta^\prime\in\mathcal{T}}.
\ee

By construction the evolution is  corresponding to the
evolution of a set-valued tableau.
Starting with a single factor or finite number of factors,  for example
$(1,\dots,1,0)^{\otimes k}$,   this determines  the
evolution at the given site in our indicator-function-tableau-valued  process $ \CG^{++,<}$.
This proves in particular Proposition \ref{P.MPsvd} for the case of a geographic space
with one single site.

The set-valued process has the following properties.

\beL{SMCM0}(Local set-valued dynamics with no migration)\\
(a) Consider the Markov process $(\zeta^{{\mathrm{m,s,d}}}_t)_{t \geq 0}$ with $d>0,\;  m>0$.
 This process is positive recurrent and returns to states $A \subseteq \I^\N$ where at most
 the first factor is different from $\I$. \\
 (b) In particular for our concrete model of $M$ lower and one upper type
 the process returns to the state $(1,1,\dots,1,0)$ (corresponding to
$(\I_M \backslash \{M+1\}) \times \I^\N$)
 after a recurrence  time called the {\em site resolution time} $\tau$  which has a
 distribution satisfying for some $\lambda>0$ and $k$ the initial number of variables
 \be{expmom}
E_k [e^{\lambda \tau}]<\infty. \qquad \square
\ee
\end{lemma}

\begin{proof}{\bf of Lemma \ref{SMCM0}}
 We first note that the number of active ranks at a site undergoes a birth and death process with linear birth rate and quadratic death rate.  Therefore the process will reach a state with only one active rank with probability one.  But then there is a positive probability that the mutation process at this rank will resolve, that is, reach a trap,  before the next selection event.

To prove the existence of a finite exponential moment, note that due to the quadratic death term there exists $k_0$ such that for $k\geq k_0$ the birth and death rates are dominated by a subcritical linear birth and death process.  Therefore if $\tau_{k_0-1}$ is the time to hit $k_0-1$, then
there exists some $\lambda^\prime >0$ such that for $k\geq k_0$,
\be{add67}
E_k[e^{\lambda^\prime \tau_{k_0-1}}]<\infty.\ee
We also note that (by recurrence) for each $k\leq k_0$, there is a positive probability of hitting  1 before $k+1$.
This means that the process returns to $k$ at most a finite number of times before hitting 1 and this random number has a geometric distribution. Therefore starting at $k\leq k_0$ the number of jumps before hitting 1 is the sum of a finite number of geometric random variables.  Since the time intervals between jumps are exponentially distributed, this implies that $\tau_1$ has a finite exponential moment. Noting that when the number of active ranks is reduced to 1 this corresponds to a column containing members of a partition of $\mathbb{I}\backslash \{M+1\}$ and therefore has union some $A\subset \mathbb{I}\backslash \{M+1\}$.   If this equals $\mathbb{I}\backslash \{M+1\}$ then $\tau=\tau_1$. Otherwise there is a positive probability that this rank will reach a trap $\I\backslash \{M+1\}$ or $\varnothing$ before the next selection event. This will then occur after a geometric number of trials. This implies the result.
q.e.d.
\end{proof}
\bi

\subsection{Set-valued dual with  migration}
\label{ss.svdualm}

We have sofar shown that we have for the non-spatial process a set-valued dual process with a Markovian dynamic.
We now turn to the spatial case, when migration can occur.
A transition due to mutation, selection and coalescence only effects the site at which the transition event occurs.
Moreover if $\wt{\mathcal{G}}^{++}_t$ {\em factors} between sites, then this property is {\em preserved} under these transitions.
Recall that we denote by $|\wt \CG^{++}_t|$ the number of sites with a factor of $\wt \CG^{++}_t$,
see (\ref{supg}) for a formal definition.

We now introduce the transition due to {\em migration} which is a transition changing both the originating (parent) site
and the site to which the migrant individual moves.
Recall that in $\CG^{++,<}$ this amounts to a column changing its mark from say $i$ to $j$
corresponding to a migration step from $i$ to $j$. More specifically
at rate $c a(i,j)$ a migration event from a specific factor at an {\em occupied} site $i$ to a site $j$ takes place.

There are two possibilities on migration at a time $t$, namely in the first case $j$ is an unoccupied site thus yielding
(1)
$|\wt{\mathcal{G}}^{++}_t|
=|\wt{\mathcal{G}}^{++}_{t-}|+1$, or $j$ is an occupied site in which case
we have (2) $|\wt \CG^{++}_t| = |\wt \CG^{++}_{t-}|$.

In case (1) the factor is removed at site $i$ and all factors to the right are shifted
one to the left and at site $j$ this migrant is placed as the first factor.

In case (2) we distinguish two subcases.

If prior to the migration the state of $\wt{\mathcal{G}}^{++}_t$ factored at these two involved sites,
then at the new site we add the individual at the first inactive rank at the new site and have the product form of the new individual and the old state at the new site.

If $\wt \CG^{++}_t$ does not factor at the site where the migrant appears, the
situation is more complicated. This case will not be dealt with here but is developed in \cite{DGsel}.


Altogether the transition for the set-valued process is well-defined by these finite rates and
well-defined transitions.

\medskip

We denote the process and  the resulting transition rates as (with $E$ denoting the state space):

\be{agdd68}
(\wt \zeta^{m,s,d,c}_t)_{t \geq 0}, \quad
\{q^{m,s,d,c}_{\zeta,\zeta^\prime}\}_{\zeta,\zeta^\prime\in {E}}.
\ee

\bi

{\bf Proof of Proposition \ref{P.MPsvd}} {\bf (Markov property)}\\

We have seen that the transitions that occur in the function-valued dual $\mathcal{G}^{++}_t$
(restricted to indicator functions)  translate into the  transitions on the corresponding set. Given any set $G\in \textsf{I}$, the set of possible transitions $G\to G^\prime$ due to selection, mutation, coalescence and migration are well defined and the transition rate is finite.  Therefore these transitions define a continuous time Markov chain $\wt{\mathcal{G}}^{++}_t$ on the countable  space $\textsf{I}$.
Therefore the functional of $\CG^{++,<}$ and the process $\wt \CG^{++}$ have the same law.
qed

\subsection{Calculating with $\wt \CG^{++}$}
\label{ss.migG++}

We first note that to carry out actual calculations it is often convenient to represent sets
by their indicators. For example (\ref{add66}) then reads:
\be{add68a}
(110)\to (010)\otimes 1 + (010)\otimes \to (110) \text{  or  } (110)\otimes (110).
\ee
Also we can then work with irreducible tableaus of indicator factors to represent subsets,
see (\ref{dgr6}) simple.

We demonstrate this at some examples. However one should have in mind that the objects
we deal with should be viewed as subsets. The art of using the duality is to use the
right representation at the right moment.

In order to make the effect of migration and its interplay with mutation more transparent,
we give some examples for the main effect that new
occupied sites are created by migration, but as a result of mutation  the parent site the offspring site can be deleted at a random time or be permanent.

\begin{example}\label{b13xx} {\em {The effects of migration for $\wt \CG^{++}$}

 We consider the case where we have two typs of low fitness and one of
higher fitness. Let us start with $(110)$ and assume that 4 births occur due to selection
acting at the ranks 1,2,3,4 and assume that we simplify the tableau by carrying out as much sums as possible
(so that each time here only one row is added), so we get:

\be{add68}
\left(
\begin{array}{lllll}
                        (010)_1&  &  &  &  \\
                        (100)_1 & (010)_1 &  &  &  \\
                        (100)_1 & (100)_1 & (010)_{1} &  &  \\
                        (100)_1 & (100)_1 & (100)_{1} & (010)_1 &  \\
                        (100)_1 & (100)_1 & (100)_{1} & (100)_1 & (110)_1
                      \end{array}
\right)
\ee
Here the missing entries are $(111)$.

If the individual at the third (local) rank below migrates to an empty site and then there are two selection operations at the new site (denoted by subscript 2) we obtain (representing the a set by its indicator function with each row given by a product of indicator functions and the different rows being the indicator function of disjoint sets):
{\scriptsize
\bea{dd0923}
&&\\&&\nonumber\\&&
\left(
\begin{array}{lllll}
                        (010)_1&  &  &  &  \\
                        (100)_1 & (010)_1 &  &  &  \\
                        (100)_1 & (100)_1 & (010)_{2} &  &  \\
                        (100)_1 & (100)_1 & (100)_{2} & (010)_1 &  \\
                        (100)_1 & (100)_1 & (100)_{2} & (100)_1 & (110)_1
                      \end{array}
\right)\bigodot \left(\begin{array}{lll}
                        (010)_2&  &    \\
                        (100)_2 & (010)_2& \\
                        (100)_2 & (100)_2 & (010)_2
                        \end{array}\right)\nonumber
\eea
}
Here $\bigodot$ indicates that the two sites are coupled.
This corresponds to the multisite tableau
\be{add69}
\begin{array}{ccccccc}
  (010)_1 &  &  &  &  &  &  \\
  (100)_1 & (010)_1 &  &  &  &  &  \\
  (100)_1 & (100)_1 & (010)_2 &  &  &  &  \\
  (100)_1 & (100)_1 & (100)_2 & (010)_2 &  &  &  \\
  (100)_1 & (100)_1 & (100)_2 & (100)_2 & (010)_2 &  &  \\
  (100)_1& (100)_1 & (100)_2 & (100)_2 & (100)_2 & (010)_1 &  \\
    (100)_1 & (100)_1 & (100)_2 & (100)_2 & (100)_2 & (100)_1 & (110)_1
\end{array}
\ee

If a mutation $1\to 2$ occurs in either column 1 or column 2, then the new site is removed, i.e.  becomes $(111)$
and therefore inactive.
On the other hand if the mutation $2\to 1$ occurs in both column 1 and column 2, and the last two columns coalesce with one of the first two columns, then the new site  is active and is decoupled from the parent site.

One additional complication arises if one of the last two columns migrates before this happens.  Then the future of this migrant site depends on the site 2 in the same way.
{\scriptsize
\bea{dd0923b}
&&\\&&\nonumber\\&&
\left(
\begin{array}{lllll}
                        (010)_1&  &  &  &  \\
                        (100)_1 & (010)_1 &  &  &  \\
                        (100)_1 & (100)_1 & (010)_{2^\ast} &  &  \\
                        (100)_1 & (100)_1 & (100)_{2^\ast} & (010)_{3^\ast} &  \\
                        (100)_1 & (100)_1 & (100)_{2^\ast} & (100)_{3^\ast} & (110)_1
                      \end{array}
\right)\bigodot \left(\begin{array}{lll}
                        (010)_2&  &    \\
                        (100)_2 & (010)_2& \\
                        (100)_2 & (100)_2 & (010)_2
                        \end{array}\right)
 \bigodot  \left(\begin{array}{l}
                        (010)_3    \\
                        (100)_3 \\
                        \end{array}\right)\nonumber
\eea
}

\be{add70}
\begin{array}{ccccccc}
  (010)_1 &  &  &  &  &  &  \\
  (100)_1 & (010)_1 &  &  &  &  &  \\
  (100)_1 & (100)_1 & (010)_2 &  &  &  &  \\
  (100)_1 & (100)_1 & (100)_2 & (010)_2 &  &  &  \\
  (100)_1 & (100)_1 & (100)_2 & (100)_2 & (010)_2 &  &  \\
  (100)_1& (100)_1 & (100)_2 & (100)_2 & (100)_2 & (010)_{3} &  \\
    (100)_1 & (100)_1 & (100)_2 & (100)_2 & (100)_2 & (100)_{3} & (110)_1
\end{array}
\ee

}
If the mutation $2\to 1$ occurs in both column 1 and column 2, and the last two columns coalesce with one of the first two columns,
then the new first site  is active.  Then if the mutation $1\to 2$ occurs in the first column at the second site,
then the third site is removed.
On the other hand if mutations $1\to 2$ occur at all three columns of the second site then the third site remains active and is decoupled.
\end{example}

Note that the
decomposition of the indicator function to the sum of products of indicator function is preserved but the constraint on which other individuals the migrating individual can coalesce with changes.

\section{Application: Ergodic theorem}
\label{s.appl}

\setcounter{secnum}{\value{section}} \setcounter{equation}{0}

In this section we demonstrate how to use the set-valued dual process to prove an
ergodic theorem for the case
\be{add196}
\I =\{1,\cdots,M\}.
\ee
We shall also consider the meanfield process $X^N$-model with geographic space $\{1,\cdots,N\}$
and migration kernel
$a(i,j) = \frac{1}{N}$ for $i,j \in \{1,2, \cdots, N\}$
and the limiting model as $N \to \infty$, the socalled McKean-Vlasov process.
(For more on this process see (\cite{DGsel}) with marginal law at time $t$
denoted $\CL_t$.)

\beT{MMergodic}{{(Ergodic theorem for $M$-type system)}}

We consider the type space $\{1,\dots,M\}$.

(a) Let $N<\infty$ and consider the exchangeably interacting system
\be{iacts1}\begin{array}{l}
X^N_t=(\mathbf{x}^N(1,t),\dots,\mathbf{x}^N(N,t))\in (\Delta_{M-1})^{N},\\
\hspace{4cm}
\mbox{ with } c>0, d\geq 0 \mbox{ and }  m_{i,j}>0  \mbox{ for all } (i,j)\in \{1,\dots,M\}.
\end{array}
\ee

Consider the distribution of $\ell$ tagged sites
$\CL[(\mathbf{x}^N(1,t),\dots,\mathbf{x}^N(\ell,t))]$.
Then for $\ell\leq N$
\be{dgr64}
\CL[(\mathbf{x}^N(1,t),\dots,\mathbf{x}^N(\ell,t))]
\Rightarrow \mu^{N,\ell}_{eq}\in\mathcal{P}((\Delta_{M-1})^\ell)  $ as $t\to\infty,
\ee
where $\mu^{N,\ell}_{eq}$   is an equilibrium state
and is independent of $X^N_0  \in (\Delta_{M-1})^{N}$.

(b) Consider the McKean-Vlasov dynamic $(\CL_t)_{t \geq 0}$ corresponding to above
set-up on $\I=\{1,\dots,M\}$.
Then
\be{dd915}\CL_t\Rightarrow \CL_\infty = \mu^{\infty}_{eq}\quad , \mbox{as  }t\to\infty,
\ee
where  $\mu^{\infty}_{eq}\in \mathcal{P}(\Delta_{M-1})$ does not depend on $\CL_0$
and is the marginal of the unique equilibrium of the McKean-Vlasov process.

(c) Now consider the $\ell$-dimensional marginal of the equilibrium measure
    $P^{0,N}_{\mathrm{eq}}, \{\mu^{N,{\ell}}_{eq}\}_{N\in\N}$ from Part (a).
    Then for every $\ell \in \N, (\ell \leq N)$, as $N\to\infty$,
    \be{dgr65}
    \mu^{N,{\ell}}_{eq}(dx_1,\dots,dx_\ell)\Rightarrow \prod_{i=1}^\ell(\mu^{\infty}_{eq}(dx_i)),
    \ee
    where $\mu^{\infty}_{eq}$ is the stationary distribution given in (\ref{dd915}).
     $\qquad \square$

(d)
Consider the system $X^N(t):= \{x^N_i(\xi,t)\}_{i\in\{1,\dots,M\},\xi\in\Omega_N}$  with $d>0$,
under the above assumptions on the $\{m(i,j)\}$  and $c_j>0$ for all $j$. Assume we
have a spatially homogeneous and shift ergodic initial condition.

In this case
\be{sh1}
\CL[X_N(t)] \Rightarrow \CL^{\Omega_N}_{eq}\quad, \mbox{  as  } t\to\infty,
 \ee
where $\CL_{eq}^{\Omega_N}\in\mathcal{P}((\Delta_{M-1})^{\Omega_N})$  is the law of a spatially
homogeneous shift-ergodic random field, which is an invariant measure of the evolution and which
is unique under all translation-invariant measures.

\end{theorem}	

\medskip

\begin{proof}{\bf of Theorem \ref{MMergodic}}
\sm

{\bf (a)}
It suffices that the joint moments of the field $\{\mathbf{x}^N(j,t),\;j=1,\dots,N\}$ converge. To show
this we use the set-valued dual $\wt{\mathcal{G}}^{++,N}_t$.  It suffices to take
the initial state for the dual given by
\be{dgr30}
\wt{\mathcal{G}}^{++,N}_0= \prod_{j=1}^{{\ell}}\prod_{i=1}^{n_j} 1_{A_{i,j}}(u_{i,j}),\; A_{i,j}\subsetneq E_0,
\ee
where $1\leq \sum_j n_j= n_\ast <\infty$ where $n_\ast$ is independent of $N$.
This means at a collection of sites $j \in \{1,\cdots,\ell\}$ we have $n_j$ factors.

{\em Case i:  $d=0, N=1$.}

In this deterministic case it suffices to consider only $n_\ast=1$ since different variables
in the dual do not interact by coalescence and hence evolve independently so that the product
structure is preserved.

First consider the case $M=2$ where type $2$ has fitness $1$ and type 1 fitness 0.

Start the dual with (01). (Note that this suffices in this case).
Consider the tableau which results from three selection events before the first mutation event
(here the special form results in the fact that selection creates only four nonzero rows):
\be{mtd02x}
\begin{array}{llll}
  (01) &\otimes 1 &\otimes 1 &\otimes 1 \\
  (10) & \otimes (01) & \otimes 1  &\otimes 1  \\
  (10) & \otimes(10) & \otimes (01)  & \otimes 1\\ (10) &\otimes (10) &\otimes (10)&\otimes (01).
\end{array}
\ee
Observe next that the mutation $2\to 1$ replaces $(01)$ by $(00)$ and this removes one row and column.
However if there is a single row, that is, the tableau $(01)$, then the result is a trap and the value is $0$.
On the other hand a mutation $1\to 2$ replaces $(01)$ by $(11)$ and this immediately produces  the value  $1$.
To see this note that  the first $m_{1,2}$ mutation terminates the process since  columns to the left of the column
$\ell$ at which the mutation occurs form $\{(01)_\ell\cup (10)_\ell\}$.

We can then view the number of active rows and columns as a birth and death process with birth rate $s$ and death rate $m_{2,1}$
and which is also terminated  forever at rate $m_{1,2}$.

We must now consider two cases $s>m_{2,1}$ and $s\leq m_{2,1}$.

In the latter case the birth and death process is recurrent and returns to $(01)$ infinitely often if no mutation occurs on each return.  Since on each such return there is a positive probability of reaching a
{\em trap} on the next transition,  the probability that the tableau has not reached a trap by time $t$ goes to $0$ as $t\to\infty$.

On the other hand if $s>m_{2,1}$ then there are an increasing number of rows each containing a factor $(01)$.  The probability that the next transition at this factor produces a $(11)$ and terminates the process is positive.
So again the probability that the tableau has not reached a trap by time $t$ goes to $0$ as $t\to\infty$.

Hence combining the cases the process reaches the trap with probability tending to $1$ as
$t\to\infty$. The value of the dual expressions is 0 or 1 irrespective of the initial value
of the process $X^N$. Hence we have convergence of moments to a value independently of the
initial state.

Now consider the case $M> 2$.

Again if the $d=0$ each birth factor (created by selection)  is eventually resolved to $(1,\dots,1)$ or $(0,\dots,0)$.  Since the probability of the latter is strictly positive there are at most finitely many which are not  $(0,\dots,0)$.
Also the initial column will finally by mutation reach the $(1,\cdots,1)$ or $(0,\cdots,0)$.
Therefore after a finite random time we again conclude that $\wt{\mathcal{G}}^{++,1}_0 = 0\mbox{  or  }1$ .
Denote the probabilities that the process reaches the trap  $(0,\dots,0)$, $(1,\dots,1)$ by $q^i_0,q^\ell_1$
respectively.
Then
\be{dgr31b}
\lim_{t\to\infty}  E[x^N_i(t)]= q^i_1\quad \text{independent of }\mathbf{x}^N(0)
\mbox{ for every } i \in E_0.
\ee

{\em Case ii: $d>0,\;N=1$.}

In this case
we must consider the initial condition (\ref{dgr30}) for $\ell =1$  and arbitrary $n_\ast$.

Again the selection process creates births (i.e. new columns and rows in the tableau).
One new feature is that a selection operation can produce more than one new row
since the number of rows doubles at each step in which $1_A f \otimes 1$
or $(1-1_A) \otimes f$ does not contain a $1_A f$ which is identically zero.

Since $d>0$, the process $|\wt{\mathcal{G}}^{++,N}_t|$, the number of active columns is always
{\em positive recurrent}. Each time the number of factors reduces to 1, there is a positive probability that the mutation process will hit a trap before the next selection event.
    Therefore $\wt{\mathcal{G}}^{++,N}_0 $ will reach a trap at a finite random time with probability 1 and then assume the value $\wt{\mathcal{G}}^{++,N}_0 = 0\mbox{  or  }1$ (see Lemma \ref{mutationlemma} and Lemma \ref{SMCM0}).
Again denote the probabilities that the process reaches the trap  $(0,\dots,0)$, $(1,\dots,1)$ by $q_0,q_1$
(note that $q_0=q_0(\ell, n_1, \cdots,n_\ell)$ similarly $q_1$).
Then
\be{dgr31}
\lim_{t\to\infty}  E[\prod_{j=1}^\ell  (x^N_j(t))^{n_j(t)}]= q_1\quad \text{independent of }\mathbf{x}^N(0).\ee
\medskip
This again shows convergence of all moments and the claim follows.
\sm

\begin{remark}
The mutation among the first $(M-1)$ types can kill these extra rows but not the primary row added by a selection event. We illustrate this as follows.
Let selection act with $1_A$ on the column number 1 with $A={3,4,5}$.
Then starting from the tableau
\be{mtd02xx}
\begin{array}{llll}
  (00001) &\otimes 1 &\otimes 1 &\otimes 1 \\
  (11100) & \otimes (00001) & \otimes 1  &\otimes 1  \\
\end{array}
\ee
we get the new tableau:

\be{mtd02xxx}
\begin{array}{llll}
  (00001) &\otimes 1 &\otimes 1 &\otimes 1 \\
  (11000) & \otimes (00001) & \otimes 1  &\otimes 1  \\
  (00100) & \otimes (00001) & \otimes 1  &\otimes 1  \\
  (11000) & \otimes(11100) & \otimes (00001)  & \otimes 1.\\
\end{array}
\ee
 The transition $(00001)\to (11111)$ no longer automatically terminates the process but does cut the tableau at the level of this factor.

\end{remark}

{\em Case iii: $d>0,\; 1<N<\infty$.}

In this case we must consider initial conditions given by (\ref{dgr30}) with initial factors at $\ell$ different sites and arbitrary $n_\ast$.  In this case the dynamics at each site are positive recurrent and the analysis at a site is as above.
However the new aspect is that before resolution at a site, migration can produce new occupied sites.

First consider the case $M=2$.

As above each site will eventually hit $(11)$ or $(00)$. If it resolves to $(11)$ the resulting value is $1$.  If it resolves to $(00)$ and produces no migrants   before resolution the value is $0$. If it resolves to $(00)$ during it lifetime it can produce a random number of migrant sites. Each of these can resolve to $(00)$ or $(11)$. If the second resolves to $(11)$ the resulting value is $1$. If it resolves to $(00)$ before producing a migrant site, the value is $0$,  otherwise it can produce a random number of offspring. As a result we have a growing number of visited sites. Since each new site can either resolve to $(11)$ or resolve to $(00)$ before producing an offspring with some probability $p_\ast>0$ the process will terminate after a random number of trials.
If one of the initial sites resolves to $(00)$ before producing a migrant the resulting value if $0$ and if it resolves to $(11)$ the value is $1$.
Hence again the process reaches the traps where $\wt \CG^{++}_t$ is the full set or is empty.
Hence as above we see that all moments converge again which gives the claim.

 Now consider the case $M>2$.

Here the principle is the same, the number of columns is a positive recurrent process
returning to 1. Then the process hits $(0,\cdots,0)$ or $(1,\cdots,1)$ with positive
probability by mutation before any selection event happens.
 Then  eventually the set-valued dual will hit a trap with all factors $\I$
 or the $\emptyset$-set (see Lemma \ref{mutationlemma} and Lemma \ref{SMCM0}) and all
 the mixed moments are computed as above giving convergence to a limit independently
 of the inital state.

{\bf (b)}
Now consider the McKean-Vlasov process in which we now use the set-valued dual with $S=\N$.
The new feature is that the dual is described by a set of occupied sites and the corresponding CMJ process describing the cloud of occupied (non-resolved sites) can be supercritical.
First of all all dual clouds starting
from different sites will become independent as $N \to \infty$
Furthermore since each site has a positive probability of resolving,  in a growing cloud of unresolved sites, eventually one of these will resolve to
$\emptyset$ and kill the whole expression.

Since the probability that this has not yet happened is monotone decreasing
and hence all moments converge to a limit which proves (\ref{dd915}).

{\bf (c)}
In order to prove the convergence fix $\ell$ and consider (\ref{dgr30}). As before $\wt{\mathcal{G}}^{++,N}$ will  hit a trap at a finite random
time.
The point is then that as $N\to\infty$ the time to hit the trap occurs before any collisions
occur tends to 1 and hence for very large $N$ the calculation is essentially the same as for the case $N=\infty$, that is, the McKean-Vlasov case
which identifies the factor in (\ref{dgr65}) as $\mu^\infty_{\mathrm{eq}}$.

{\bf (d)}
Now consider the ergodic theorem for the $M$-type interacting system on the hierarchical group $\Omega_N$ in the case $m_{i,j}>0$ for all $i,j\in \{1,\dots, M\}$. We must now prove the convergence of the joint moments. To do this we start $n_i$ factors at each of $\ell$ different sites.

The main point is that the dual system eventually hits a trap
and as before in b) this will guarantee convergence of all moments. The spatial
homogeneity follows from the one of the initial state if this has this property,
but since all initial states lead to the same limit we are done.
In fact we see that the system with collisions will tend to hit a trap faster that the one for the McKean-Vlasov dual.

Furthermore we note that two initial clouds of factors starting at two sites which are at hierarchical distance
distance $\ell$ will have a probability to ever meet before they resolve to 0 or 1 which
tends to 0 as $\ell \to \infty$.

More precisely, if $\lambda$ is the parameter in the exponential tail of the trapping time,
the probability of a collision between two clouds starting at two sites at hierarchical distance $\ell$  before
trapping is
\be{}
\leq\rm{const}\cdot  \frac{\frac{c_\ell}{N^{2\ell}}}{\lambda+\frac{c_\ell}{N^{2\ell}}}.\ee
Therefore the {\em all} the mixed moments connected with two sites at distance $\ell$ factor asymptotically as $\ell \to \infty$.
 Hence the limiting state is shift-ergodic. $\qquad\blacksquare$

\end{proof}

\bigskip

\newpage

\begin{center}

\end{center}

\newpage
\section* {Index of Notation}
\begin{itemize}
\item $\Z$ - the set of integers \item $\N = \{1,2,3,\dots \}$
\item $\Omega_N = \otimes_\N Z_N \quad , \quad Z_N$ cyclical group
of order $N$ \item $\mathcal{M}(E)$ denotes the space of finite
Borel measures on a Polish space $E$ \item $\CP(E)$ denotes the
space of probability measures on the Borel field on a Polish space
$E$.
\end{itemize}

\end{document}